\documentclass[10pt]{article}
\usepackage{amsmath}
\usepackage{amsfonts}
\usepackage{amscd}
\usepackage{amsthm}
\usepackage{a4}
\usepackage{graphicx, color, epsfig}
\usepackage[english]{babel}
\usepackage{titlesec}
\titleformat{\section}{\large\bfseries}{\thesection}{1em}{}
\newtheorem{theorem}{Theorem}[section]

\newtheorem{prop}[theorem]{Proposition}

\newenvironment{PROOF}[1][Proof]{\begin{trivlist}
\item[\hskip \labelsep {\bfseries #1}]}{\end{trivlist}}

\newcommand{\QED}{\hfill $\Box$}

\begin{document}

\numberwithin{figure}{section}
\numberwithin{equation}{section} 
\rmfamily
\setlength{\parindent}{15pt}
\title{The Kadomtsev-Petviashvili II Equation on the Half-Plane}
\author{
{D. Mantzavinos \& A. S. Fokas}\\
\\
{Department of Applied Mathematics and Theoretical Physics,}\\
{University of Cambridge, Cambridge CB3 0WA, UK.}}
%{E-mail:~t.fokas@damtp.cam.ac.uk}\\
\maketitle

\begin{abstract}
The KPII equation is an integrable nonlinear PDE in 2+1 dimensions (two spatial and one temporal), which arises in several physical circumstances, including fluid mechanics where it describes waves in shallow water. It provides a multidimensional generalisation of the renowned KdV equation. In this work, we employ a novel approach recently introduced by one of the authors in connection with the Davey-Stewartson equation \cite{FDS2009}, in order to analyse the initial-boundary value problem for the KPII equation formulated on the half-plane. The analysis makes crucial use of the so-called d-bar formalism, as well as of the so-called global relation. A novel feature of boundary as opposed to initial-value problems in 2+1 is that the d-bar formalism now involves a function in the complex plane which is discontinuous across the real axis. 
\end{abstract}

\normalsize

\section{Introduction}

The Kadomtsev-Petviashvili (KP) equation
\begin{equation}\label{kpintro}
q_t+6 q q_x+q_{xxx}+3\, \sigma\, \partial_{x}^{-1}q_{yy}=0, \quad \sigma=\pm 1,
\end{equation}
with the operator $\partial_{x}^{-1}$ defined by
\begin{equation*}
\partial_{x}^{-1}f(x)=\int_{-\infty}^{x}f(\xi)\, d\xi ,
\end{equation*}
is one of the most notable integrable nonlinear equations in 2+1 dimensions (i.e. evolution equations in two spatial dimensions). It is the natural generalisation of the celebrated Korteweg-de Vries (KdV) equation from one to two spatial dimensions and, as such, it appears in various physical problems. In particular, it was first derived in the study of waves of long wavelength in shallow water. In this specific application, if the surface tension dominates over the gravitational force, then $\sigma=-1$ and equation \eqref{kpintro} is called KPI \cite{FAbl1983}, whereas if the gravitational force is dominant, then $\sigma=1$ and equation \eqref{kpintro} is called KPII \cite{AblBarF1983}.

The initial-value problem for the KdV equation was solved via the so-called Inverse Scattering transform (IST) in 1967 \cite{GGKM}; KPI and KPII equations were formally solved in \cite{FAbl1983} and \cite{AblBarF1983} respectively, using a nonlocal Riemann-Hilbert formalism and a d-bar formalism respectively (see also \cite{BLP1989}, \cite{FSung1992}, \cite{Fkp2009} \cite{BC1984}). The solution of initial-boundary value (IBV) problems is considerably more complicated than the solution of pure initial-value problems. A unified transform method for solving linear and integrable nonlinear evolution PDEs in one spatial dimension was introduced in \cite{F1997}. In particular, interesting results for the linearised version of the KdV and for the KdV itself formulated on the half-line are presented in \cite{FPelloni1997}, \cite{F2002b}, \cite{Fbook} and \cite{FTreharne2008}. The generalisation of these results from one to two spatial dimensions for linear and for integrable nonlinear equations is presented in \cite{FPelloni1997} and \cite{FDS2009} respectively.

Here, we employ the general methodology of \cite{FDS2009} in order to analyse the KPII on the half-plane. It has been emphasised by one of the authors that the solution of the linearised version of a given nonlinear PDE via a Lax pair approach provides a useful starting point before analysing the nonlinear PDE itself. Hence, we will first solve the linearised version of the KPII equation on the half-plane,
\begin{equation}
q_t+q_{xxx}+3\partial_x^{-1}q_{yy}=0, \quad (x,y,t)\in \Omega, \label{lkpintro}
\end{equation}
where 
\begin{equation}\label{Omegaintro}
\Omega:=\{-\infty<x<\infty,\ 0<y<\infty,\ t>0\},
\end{equation}
before analysing the KPII equation,
\begin{equation}\label{kpiiintro}
q_t+6qq_x+q_{xxx}+3\partial^{-1}_x q_{yy}=0,  \quad (x,y,t)\in \Omega.
\end{equation} 

Our approach to both problems involves the following steps:
\begin{enumerate}
\item \textit{The formulation of the PDE in terms of a Lax pair.} For the spectral variable $k=k_R+i k_I,\ k_R,\, k_I \in \mathbb R$, and for some \textit{sectionally analytic} function $\mu=\mu(x,y,t,k_R,k_I)$, the linearised KP equation \eqref{lkpintro} admits the Lax pair: 
\begin{subequations}\label{Laxintro}
\begin{equation}
\mu_y-\mu_{xx}-2ik\mu_x=q,\label{Lax1intro}
\end{equation}
\begin{equation}
\mu_t+4\mu_{xxx}+12ik\mu_{xx}-12k^2\mu_x=-3\left(q_x+2ikq+\partial_x^{-1}q_y\right).\label{Lax2intro}
\end{equation}
\end{subequations}
Similarly, the KPII equation \eqref{kpiiintro} possesses the Lax pair:
\begin{subequations}\label{laxintro}
\begin{equation}\label{lax1intro}
\mu_y-\mu_{xx}-2ik\mu_x=q\mu,
\end{equation}
\begin{equation}\label{lax2intro}
\mu_t+4\mu_{xxx}+12ik\mu_{xx}-12k^2\mu_x=F\mu,
\end{equation}
\end{subequations}
where the operator $F$ is defined by
\begin{equation}\label{Fintro}
 F(x,y,t,k)= -6q\left(\partial_x+ik\right)-3\left(q_x+\partial^{-1}_xq_y\right).
\end{equation}

\item \textit{The direct problem.} By applying a Fourier transform in $x$, it is possible to analyse the two equations defining the Lax pair \textit{simultaneously} in order to obtain an expression for $\mu$ which is bounded for all $k \in \mathbb C$. It turns out that $\mu$ has different representations in different parts of the complex $k$-plane, namely
\begin{align}\label{spectralfunctionsintro}
\mu(x,y,t,k_R,k_I) = \left\lbrace \begin{array}{ll} \mu_{2}^{+}(x,y,t,k_R,k_I), &\quad k\in I, \\ \\ \mu_{1}^{+}(x,y,t,k_R,k_I), &\quad k\in II,  \\ \\ \mu_{1}^{-}(x,y,t,k_R,k_I), &\quad k\in III, \\ \\ \mu_{2}^{-}(x,y,t,k_R,k_I), &\quad k\in IV, \end{array} \right.
\end{align}
where $I-IV$ denote the four quadrants of the complex $k$-plane.

The above expression for $\mu_{1,2}^\pm$ depends on $q(x,y,t)$, $q(x,0,t)$ and $q_y(x,0,t)$.

\item \textit{The derivation of the global relation.} This relation is an algebraic equation coupling the so-called spectral functions. For the linearised KP these functions are appropriate transforms of
\begin{equation}\label{ibvintro}
 q_0(x,y)=q(x,y,0),\ g(x,t)=q(x,0,t),\ h(x,t)=q_y(x,0,t); 
\end{equation}
namely, the global relation couples $\hat q$, $\hat q_0$, $\tilde g_T$ and $\tilde h_T$, where
\begin{subequations}
\begin{align}
\hat q(k_1,k_2,t)&=\int_{-\infty}^{\infty}dx\int_{0}^{\infty}dy\, e^{-ik_1x-i k_2y} q(x,y,t)\label{qhatintro},\\
\hat q_0(k_1,k_2)&=\int_{-\infty}^{\infty}dx\int_{0}^{\infty}dy\, e^{-ik_1x-i k_2 y} q_0(x,y)\label{q0hatintro},\\
\tilde g_T(k_1,k_2)&=\int_{-\infty}^{\infty}dx\int_{0}^{T}d\tau\, e^{-ik_1x+\omega(k) \tau} g(x,\tau)\label{gttildeintro},\\
\tilde h_T(k_1,k_2)&=\int_{-\infty}^{\infty}dx\int_{0}^{T}d\tau\, e^{-ik_1x+\omega(k) \tau} h(x,\tau)\label{httildeintro}, \quad k_1\in \mathbb R,\ \mathrm{Im}k_2\leq 0
\end{align}
\end{subequations}
and 
\begin{equation}\label{omegaintro}
\omega(k):=-ik_1^3+3i\frac{k_2^2}{k_1}.
\end{equation}
For the KPII, the spectral functions can also be expressed in terms of \eqref{ibvintro} via linear integral equations.

\item \textit{The inverse problem.} By using the fact that the function $\mu$ defined in step 2 is bounded for all $k\in \mathbb C$, it is possible to obtain an alternative representation for this function using a d-bar formalism (or more precisely the so-called Pompeiu's formula). In order to achieve this, it is necessary to: (a) compute $\partial \mu/\partial \bar k$, (b) compute the jumps of $\mu_{1,2}^\pm$ across the real and imaginary $k$-axes; actually  the jump across the imaginary $k$-axis vanishes. For the linearised KP, $\partial \mu/\partial \bar k$ and $(\mu^+_{1,2}-\mu^-_{1,2})$ can be expressed in terms of $\hat q_0$, $\tilde g_T$ and $\tilde h_T$. Then, using Pompeiu's formula, it is possible to express $\mu$ in terms of $\hat q_0$, $\tilde g_T$ and $\tilde h_T$. In the nonlinear case, the d-bar derivatives and the jumps can be expressed in terms of the spectral functions. In summary, $\mu$ can be expressed via the spectral functions and hence via $q_0$, $g$ and $h$. After obtaining $\mu$ it is straightforward to obtain a formula for $q$.

\end{enumerate}

\paragraph{Notations and Assumptions}

\begin{itemize}
\item The complex variable $k$ is defined as
\begin{equation*}
 k = k_R + ik_I, \quad k_R,\, k_I \in \mathbb R.
\end{equation*}
\item A bar on top of a complex variable will denote the complex conjugate of this variable; in particular, $\bar k = k_R-ik_I$.
\item For the solution of the inverse problem we will make use of the so-called Pompeiu's formula: if $f(x,y)$ is a smooth function in some piece-wise smooth domain $ \mathcal D\subset \mathbb R^2$, then $f$ is related to its value on the boundary of $\mathcal D$ and to its d-bar derivative inside $\mathcal D$ via the equation
\begin{equation}\label{pompeiuintro}
f(x,y)=\frac{1}{2i\pi}\int_{\partial \mathcal D}\frac{d\zeta}{\zeta-z}\, f(\xi,\eta)+\frac{1}{2i\pi}\iint_{\mathcal D}\frac{d\zeta\wedge d\bar\zeta}{\zeta-z}\, \frac{\partial f}{\partial \overline \zeta}\,(\xi, \eta), \quad \zeta=\xi+i \eta,
\end{equation}
where 
\begin{equation*}
 d \zeta \wedge d \bar \zeta=-2i d \xi d \eta.
\end{equation*}
\item We will denote the initial value by
\begin{subequations}
\begin{equation}
q(x,y,0)=q_0(x,y), \label{icintro}
\end{equation}
and the boundary values $q(x,0,t)$ and $q_y(x,0,t)$ by
\begin{equation}
q(x,0,t)=g(x,t)\label{lkpbc1intro}
\end{equation}
and
\begin{equation}\label{lkpbc2intro}
q_y(x,0,t)=h(x,t).
\end{equation}
\end{subequations}
We will assume that $q_0\in \mathbb S(\mathbb R \times \mathbb R^+)$, where $\mathbb S$ denotes the space of Schwartz functions.
\item We will seek a solution which decays as $y\rightarrow \infty$ for all fixed $(x,t)$ and which also decays as $|x|\rightarrow \infty$ for all fixed $(y,t)$.
\item Throughout this paper we will {\it assume} that there exists a solution $q(x,y,t)$ with sufficient smoothness and decay in $\bar \Omega$, which denotes the closure of the domain $\Omega$.
\item A hat ``$\wedge$''above a function will denote the Fourier transform of this function in the variable $x$.
\item The index ``o'' on a function will denote evaluation at $t=0$.
\item Whenever we write $q$, we mean that $q$ depends on the physical variables $(x,y,t)$.
\item Whenever we write $\mu$, we mean that $\mu$ depends on the physical variables $(x,y,t)$ and on the spectral variables $k_R$, $k_I\in \mathbb R$.
\end{itemize}

\section{The linearised KP equation}\label{lkp}

We consider the linearised KP equation \eqref{lkpintro} with the associated Lax pair \eqref{Laxintro}. Differentiating equation \eqref{Lax1intro} with respect to $t$ and equation \eqref{Lax2intro} with respect to $y$, we find:
\begin{equation*}
\mu_{yt}-\mu_{xxt}-2ik\mu_{xt}=q_t  
\end{equation*}
and
\begin{equation*}
\mu_{ty}+4\mu_{xxxy}+12ik\mu_{xxy}-12k^2\mu_{xy}=-3\left(q_{xy}+2ikq_{y}+\partial_{x}^{-1}q_{yy}\right).
\end{equation*}
The requirement $\mu_{ty}=\mu_{yt}$ implies equation \eqref{lkpintro}, thus \eqref{Laxintro} is indeed a Lax pair for equation \eqref{lkpintro}.

\subsection{The direct problem}\label{lkpdp}

\begin{figure}[ht]
\begin{center}
\resizebox{4cm}{!}{\input{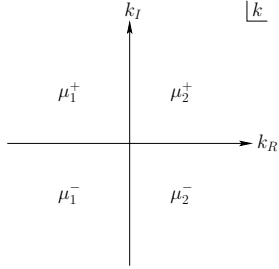}}
\end{center}
\caption{The function $\mu$ in the four quadrants of the complex $k$-plane.}
\label{eigenfunctions}
\end{figure}

\begin{prop}\label{lkpdpprop} 
Assume that there exists a solution $q(x,y,t), \ (x,y,t)\in \Omega,$ to an initial-boundary value problem for equation \eqref{lkpintro}. Then, there exists a solution $\mu$ of the Lax pair \eqref{Laxintro} which is bounded $\forall k \in \mathbb C$. This solution can be represented in the form
\begin{align}\label{spectralfunctions}
\mu(x,y,t,k_R,k_I) = \left\lbrace \begin{array}{ll} \mu_{1}^{+}(x,y,t,k_R,k_I), &\quad k_R\leq0,\ k_I\geq0, \\ \\ \mu_{1}^{-}(x,y,t,k_R,k_I), &\quad k_R\leq0,\ k_I\leq0,  \\ \\ \mu_{2}^{-}(x,y,t,k_R,k_I), &\quad k_R\geq0,\ k_I\leq0, \\ \\ \mu_{2}^{+}(x,y,t,k_R,k_I), &\quad k_R\geq0,\ k_I\geq0, \end{array} \right.
\end{align}
(see figure \ref{eigenfunctions}), where $\mu_{1,2}^\pm$ can be expressed in terms of $\{q(x,y,t), g(x,t), h(x,t)\}$ by the following formulae: 
{\small
\begin{align}\label{direct1}
&\mu_{1}^{\pm}=\frac{1}{2\pi}\mbox{\large $\Big($}\int_{-\infty}^{0}dl+\int_{-2k_R}^{\infty}dl\mbox{\large $\Big)$ }\int_{-\infty}^{\infty}d\xi\int_{0}^{y}d\eta\, e^{-il(\xi-x)+l(l+2k)(\eta-y)}q(\xi,\eta,t)\nonumber\\ 
&+ \frac{1}{2\pi}\mbox{\large $\Big($}\!\int_{-\infty}^{0}\!\!\!dl+\int_{-2k_R}^{\infty}\!\!\!dl\mbox{\large $\Big)$ }\!\!\!\!\!\int_{-\infty}^{\infty}\!\!\!d\xi\left\{\begin{array}{ll}\int_{0}^{t}d\tau \\ \\ -\int_{t}^{T}d\tau \end{array}\right.\!\!\!\!\! \!\!e^{-il(\xi-x)-l(l+2k)y-4il(l^2+3kl+3k^2)(\tau-t)} Q(\xi,\tau,k,l)\nonumber\\
&-\frac{1}{2\pi}\int_{0}^{-2k_R}dl \int_{-\infty}^{\infty}d\xi\int_{y}^{\infty}d\eta \, e^{-il(\xi-x)+l(l+2k)(\eta-y)}q(\xi,\eta,t),
\end{align}}
and
{\small
\begin{align}\label{direct2}
&\mu_{2}^{\pm}=\frac{1}{2\pi}\mbox{\large $\Big($}\int_{-\infty}^{-2k_R}dl+\int_{0}^{\infty}dl\mbox{\large $\Big)$ }\int_{-\infty}^{\infty}d\xi\int_{0}^{y}d\eta\, e^{-il(\xi-x)+l(l+2k)(\eta-y)}q(\xi,\eta,t)\nonumber\\ 
&+ \frac{1}{2\pi}\mbox{\large $\Big($}\!\int_{-\infty}^{-2k_R}\!\!\!dl+\int_{0}^{\infty}\!\!\!\!dl\mbox{\large $\Big)$ } \!\!\!\!\int_{-\infty}^{\infty}\!\!\!\!d\xi\left\{\!\!\begin{array}{ll}\int_{0}^{t}d\tau \\ \\ -\int_{t}^{T}d\tau \end{array}\right.\!\!\!\!\! e^{-il(\xi-x)-l(l+2k)y-4il(l^2+3kl+3k^2)(\tau-t)} Q(\xi,\tau,k,l)\nonumber\\
&-\frac{1}{2\pi}\int_{-2k_R}^{0}dl \int_{-\infty}^{\infty}d\xi\int_{y}^{\infty}d\eta \, e^{-il(\xi-x)+l(l+2k)(\eta-y)}q(\xi,\eta,t),
\end{align}}
where
\begin{equation}\label{Q}
Q(x,t,k,l)=-3\bigg[i(l+2k)g(x,t)+\partial_{x}^{-1}h(x,t)\bigg],
\end{equation}
and the functions $g$ and $h$ are defined by equations \eqref{lkpbc1intro} and \eqref{lkpbc2intro}.
\end{prop}

\begin{PROOF}

Let $\hat \mu$ denote the Fourier transform of $\mu$ with respect to $x$, i.e. 
\begin{subequations}
\begin{equation}\label{lkpft}
\hat{\mu}(l,y,t,k_R,k_I)=\int_{-\infty}^{\infty}dx\, e^{-ilx}\mu(x,y,t,k_R,k_I), \quad l \in \mathbb{R},\ k \in \mathbb C,
\end{equation}
with the inverse Fourier transform given by
\begin{equation}\label{lkpift}
\mu(x,y,t,k_R,k_I)=\frac{1}{2\pi}\int_{-\infty}^{\infty} dl\, e^{ilx}\hat{\mu}(l,y,t,k_R,k_I), \quad k \in \mathbb C.
\end{equation}
\end{subequations}

Then, the Fourier transform of equation \eqref{Lax1intro} gives:
\begin{equation*}
\bigg(\hat{\mu}(l,y,t,k_R,k_I)\, e^{l(l+2k)y}\bigg)_y=\hat{q}(l,y,t)e^{l(l+2k)y}.
\end{equation*}
We solve this equation by integrating either from $y=0$ or $y=-\infty$:
\begin{align}\label{if}
\hat{\mu}(l,y,t,k_R,k_I)=\left\{\begin{array}{ll} \int_0^y d\eta\, e^{l(l+2k)(\eta-y)}\hat{q}(l,\eta,t)+\hat{\mu}(l,0,t,k_R,k_I)\, e^{-l(l+2k)y} \\ \\ -\int_y^\infty d\eta\, e^{l(l+2k)(\eta-y)}\hat{q}(l,\eta,t). \end{array}\right.
\end{align}
We look for a solution which is bounded $\forall k \in \mathbb C$. In the first expression in equation \eqref{if}, $\eta-y\leq 0$, thus we require that $\mathrm{Re}\left[l(l+2k)\right]\geq 0$, i.e.
\begin{equation*}
l(l+2k_R)\geq 0\Leftrightarrow \left\lbrace\begin{array}{ll} l\in(-\infty,0]\cup[-2k_R,\infty),\ k_R\leq0 \\ \\ l\in(-\infty,-2k_R]\cup[0,\infty),\ k_R\geq0. \end{array} \right. 
\end{equation*}
Analogous considerations are valid for the second expression in equation \eqref{if}, where $\eta-y\geq 0$. 

Let 
\begin{align}\label{star}
\mu(x,y,t,k_R,k_I)&=\left\lbrace\begin{array}{ll}\mu_1(x,y,t,k_R,k_I), \quad k_R\leq0\\ \\ \mu_2(x,y,t,k_R,k_I), \quad k_R\geq0.\end{array}\right.
\end{align}
Inverting equation \eqref{if} we obtain:
\begin{align}\label{left}
\mu_1(x,y,t,k_R,k_I)&= \frac{1}{2\pi}\left(\int_{-\infty}^{0}dl+\int_{-2k_R}^{\infty}dl\right)e^{ilx-l(l+2k)y}\int_{0}^{y}d\eta\, e^{l(l+2k)\eta}\hat q(l,\eta,t)\nonumber\\
&+\frac{1}{2\pi}\left(\int_{-\infty}^{0}dl+\int_{-2k_R}^{\infty}dl\right)e^{ilx-l(l+2k)y}\hat\phi_1(l,t,k_R,k_I)\nonumber\\
&-\frac{1}{2\pi}\int_{0}^{-2k_R}dl\, e^{ilx-l(l+2k)y}\int_{y}^{\infty}d\eta \, e^{l(l+2k)\eta}\hat q(l,\eta,t), \quad k_R\leq0
\end{align}
and
\begin{align}\label{right}
\mu_2(x,y,t,k_R,k_I)&= \frac{1}{2\pi}\left(\int_{-\infty}^{-2k_R}dl+\int_{0}^{\infty}dl\right)e^{ilx-l(l+2k)y}\int_{0}^{y}d\eta\, e^{l(l+2k)\eta}\hat q(l,\eta,t)\nonumber\\
&+ \frac{1}{2\pi}\left(\int_{-\infty}^{-2k_R}dl+\int_{0}^{\infty}dl\right)e^{ilx-l(l+2k)y}\hat\phi_2(l,t,k_R,k_I)\nonumber\\
&-\frac{1}{2\pi}\int_{-2k_R}^{0}dl\, e^{ilx-l(l+2k)y}\int_{y}^{\infty}d\eta \, e^{l(l+2k)\eta}\hat q(l,\eta,t), \quad k_R\geq0,
\end{align}
where $\phi_{1,2}(x,t,k_R,k_I)$ denote the value of $\mu_{1,2}(x,y,t,k_R,k_I)$ at the boundary $y=0$, i.e.
\begin{equation}
\phi_j(x,t,k_R,k_I):=\mu_j(x,0,t,k_R,k_I),\quad j=1,2.
\end{equation}
Equation \eqref{left} evaluated at $y=0$ yields
\begin{align*}
\phi_1(x,t,k_R,k_I)&= \frac{1}{2\pi}\left(\int_{-\infty}^{0}dl+\int_{-2k_R}^{\infty}dl\right)\, e^{ilx}\hat\phi_1(l,t,k_R,k_I)\\
&-\frac{1}{2\pi}\int_{0}^{-2k_R}dl\, e^{ilx}\int_{0}^{\infty}d\eta\, e^{l(l+2k)\eta}\hat q(l,\eta,t).
\end{align*}
Using the definition \eqref{lkpift} with $y=0$, we find
\begin{align*}
\frac{1}{2\pi}\int_{-\infty}^{\infty}dl\, e^{ilx} \hat \phi_1(l,t,k_R,k_I)&= \frac{1}{2\pi}\left(\int_{-\infty}^{0}dl+\int_{-2k_R}^{\infty}dl\right)\, e^{ilx}\hat\phi_1(l,t,k_R,k_I)\\
&-\frac{1}{2\pi}\int_{0}^{-2k_R}dl\, e^{ilx}\int_{0}^{\infty}d\eta\, e^{l(l+2k)\eta}\hat q(l,\eta,t).
\end{align*}
Thus, $\hat \phi_1$ does not satisfy any constraints for $l\in (-\infty,0]\cup [-2k_R, \infty)$, whereas it does satisfy the following constraint for $l\in [0,-2k_R]$:
\begin{equation}\label{imposed}
\hat \phi_1(l,t,k_R,k_I)=-\int_{0}^{\infty}d\eta \, e^{l(l+2k)\eta}\hat q(l, \eta,t), \quad l\in [0,-2k_R].
\end{equation} 
The situation is similar for $\phi_2$, namely $\hat \phi_2$ does not satisfy any constraint for $l\in (-\infty,-2k_R)\cup (0, \infty)$, whereas it does satisfy the following constraint for $l\in [-2k_R,0]$:
\begin{equation*}
\hat \phi_2(l,t,k_R,k_I)=-\int_{0}^{\infty}d\eta \, e^{l(l+2k)\eta}\hat q(l, \eta,t) \quad l\in [-2k_R,0].
\end{equation*}
In order to determine $\hat\phi(l,t,k_R,k_I)$, we use the second of the equations defining the Lax pair, evaluated at $y=0$: 
\begin{equation}\label{mu0}
\phi_t+4\phi_{xxx}+12ik\phi_{xx}-12k^2\phi_x=F(x,t,k_R,k_I), 
\end{equation}
where 
\begin{align}\label{lF}
F(x,t,k_R,k_I):&=-3\left(q_x+2ikq+\partial_x^{-1}q_y\right)\bigg|_{y=0}\nonumber\\
&=-3\bigg[g_x(x,t)+2ikg(x,t)+\partial_x^{-1}h(x,t)\bigg].
\end{align}
The Fourier transform in $x$ of equation \eqref{mu0} yields
\begin{equation*}
\left(\hat \phi(l,t,k_R,k_I) e^{-4il(l^2+3kl+3k^2)t}\right)_t=\hat F(l,t,k_R,k_I) e^{-4il(l^2+3kl+3k^2)t}.
\end{equation*}
Thus, integrating with respect to $t$ either from $t=0$ or from $t=T$, we find
\begin{align*}
\hat \phi=\left\lbrace\begin{array}{ll}\int_{0}^{t}d\tau \, e^{-4il(l^2+3kl+3k^2)(\tau-t)} \hat F(l,\tau,k_R,k_I) + e^{4il(l^2+3kl+3k^2)t} \hat \phi_{0} \\ \\ \!\!\!\!-\!\int_{t}^{T}\!\!d\tau e^{-4il(l^2+3kl+3k^2)(\tau-t)} \hat F(l,\tau,k_R,k_I) + e^{-4il(l^2+3kl+3k^2)(T-t)} \hat \phi_{T},\  k \in \mathbb C. \end{array}\right.
\end{align*}
Moreover, integration by parts of $\hat F$ implies
\begin{align}\label{if'}
\hat \phi=\left\lbrace\begin{array}{ll}\int_{0}^{t}d\tau \, e^{-4il(l^2+3kl+3k^2)(\tau-t)} \hat Q(l,\tau,k_R,k_I) + e^{4il(l^2+3kl+3k^2)t}\hat \phi_0 \\ \\ \!\!\!\!-\!\int_{t}^{T}\!\!d\tau e^{-4il(l^2+3kl+3k^2)(\tau-t)} \hat Q(l,\tau,k_R,k_I) + e^{-4il(l^2+3kl+3k^2)(T-t)}\hat \phi_T,\ k \in \mathbb C.  \end{array}\right.
\end{align}
where $\hat Q$ is the Fourier transform \eqref{lkpft} of the function $Q$ defined by equation \eqref{Q} and $\phi_T$ denotes the evaluation of $\phi$ at $t=T$.

We look for a solution $\hat \phi$ which is bounded $\forall k \in \mathbb C$. Nothing that
\begin{equation*}
\mathrm{Re}\bigg[-4il(l^2+3kl+3k^2)(\tau-t)\bigg]=12lk_I(l+2k_R)(\tau-t),
\end{equation*} 
it follows that the exponential $e^{-4il(l^2+3kl+3k^2)(\tau-t)}$ is bounded, provided that
\begin{equation*}
lk_I(l+2k_R)(\tau-t)\leq0.
\end{equation*}
In order to analyse this inequality, we consider the following cases:
\begin{itemize}
\item If $k\in\{k\in \mathbb{C}: k_R\leq0,k_I\geq0\}$, then $l\in(-\infty,0)\cup (-2k_R, \infty)$ implies that $\tau-t\leq 0$, whereas $l\in(0,-2k_R)$ implies that $\tau-t\geq 0$. Recalling our earlier analysis for $\hat \phi_1$, we can choose 
\begin{equation*}
\hat \phi_1^+(l,0,k_R,k_I)=0 \quad \mathrm{for}\ l\in(-\infty, 0) \cup (-2k_R,\infty)  
\end{equation*}
and then relation \eqref{imposed} implies
\begin{equation*}
\hat \phi_1^+(l,T,k_R,k_I)=-\int_{0}^{\infty}d \eta \, e^{l(l+2k)\eta}\hat q(l,\eta,T) \quad \mathrm{for}\ l\in (0, -2k_R).
\end{equation*}
Consequently, equations \eqref{if'} yield
\begin{align}\label{choice1}
\hat \phi_1^+=\left\lbrace\begin{array}{ll}\int_{0}^{t}d\tau \, e^{-4il(l^2+3kl+3k^2)(\tau-t)} \hat Q(l,\tau,k), \quad l\in(-\infty,0)\cup(-2k_R, \infty)\\ \\ \\ -\int_{t}^{T}d\tau \, e^{-4il(l^2+3kl+3k^2)(\tau-t)} \hat Q(l,\tau,k)\\ \\-e^{-4il(l^2+3kl+3k^2)(T-t)}\int_{0}^{\infty}d \eta \, e^{l(l+2k)\eta}\hat q(l,\eta,T), \quad l\in(0,-2k_R).\end{array}\right.
\end{align} 
\item If $k\in\{k\in \mathbb C:k_R\leq0,k_I\leq0\}$, then $l\in(-\infty,0]\cup[-2k_R, \infty)$ implies that $\tau-t\geq 0$, whereas $l\in[0,-2k_R]$ implies that $\tau-t\leq 0$. Following steps similar to those used above, we find
\begin{equation*}
\hat \phi_1^-(l,0,k_R,k_I)=-\int_{0}^{\infty}d \eta \, e^{l(l+2k)\eta}\hat q_0(l,\eta) \quad \mathrm{for} \ l\in[0,-2k_R]
\end{equation*}
and
\begin{equation*}
\hat \phi_1^-(l,T,k_R,k_I)=0 \quad \mathrm{for} \ l\in(-\infty,0]\cup[-2k_R, \infty).
\end{equation*}
Then, equations \eqref{if'} imply
\begin{align}\label{choice2}
\hat \phi_1^-=\left\lbrace\begin{array}{ll}\!\!\!-\int_{t}^{T}d\tau \, e^{-4il(l^2+3kl+3k^2)(\tau-t)} \hat Q(l,\tau,k_R,k_I), \ l\in(-\infty,0]\cup[-2k_R, \infty)\\ \\ \\ \int_{0}^{t}d\tau \, e^{-4il(l^2+3kl+3k^2)(\tau-t)} \hat Q(l,\tau,k_R,k_I)\\ \\
-e^{4il(l^2+3kl+3k^2)t}\int_{0}^{\infty}d \eta \, e^{l(l+2k)\eta}\hat q_0(l,\eta), \quad l\in[0,-2k_R].\end{array}\right.
\end{align}
\item Similarly, we can show that for $k\in \{k\in \mathbb C: k_R\geq0, k_I\geq0\}$,
\begin{align}\label{choice3}
\hat \phi_2^+=\left\lbrace\begin{array}{ll}\int_{0}^{t}d\tau \, e^{-4il(l^2+3kl+3k^2)(\tau-t)} \hat Q(l,\tau,k_R,k_I), \quad l\in(-\infty,-2k_R]\cup[0, \infty)\\ \\ \\ -\int_{t}^{T}d\tau \, e^{-4il(l^2+3kl+3k^2)(\tau-t)} \hat Q(l,\tau,k_R,k_I)\\ \\-e^{-4il(l^2+3kl+3k^2)(T-t)}\int_{0}^{\infty}d \eta \, e^{l(l+2k)\eta}\hat q(l,\eta,T), \quad l\in[-2k_R,0].\end{array}\right.
\end{align} 
\item For $k\in\{k\in \mathbb C:k_R\geq0, k_I\leq0\}$,
\begin{align}\label{choice4}
\hat \phi_2^-=\left\lbrace\begin{array}{ll}\!\!\!-\int_{t}^{T}d\tau \, e^{-4il(l^2+3kl+3k^2)(\tau-t)} \hat Q(l,\tau,k_R,k_I), \ l\in(-\infty,-2k_R]\cup[0, \infty)\\ \\ \\ \int_{0}^{t}d\tau \, e^{-4il(l^2+3kl+3k^2)(\tau-t)} \hat Q(l,\tau,k_R,k_I)\\ \\
-e^{4il(l^2+3kl+3k^2)t}\int_{0}^{\infty}d \eta \, e^{l(l+2k)\eta}\hat q(l,\eta,0), \quad l\in[-2k_R, 0].\end{array}\right.
\end{align}
\end{itemize}
Inserting equations \eqref{choice1}-\eqref{choice4} in equations \eqref{left} and \eqref{right}, we find equations \eqref{direct1} and \eqref{direct2} respectively.
\QED 
\end{PROOF}

\subsection{The global relation}\label{lkpgr}

\begin{prop}\label{GRprop}
Define $\hat q(k_1,k_2,t)$, $\tilde g_t(k_1,k_2)$ and $\tilde h_t(k_1,k_2)$ by
\begin{subequations}\label{transfs}
\begin{align}
\hat q(k_1,k_2,t)&=\int_{-\infty}^{\infty}dx\int_{0}^{\infty}dy\, e^{-ik_1x-i k_2y} q(x,y,t)\label{qhatgr},\\
\hat q_0(k_1,k_2)&=\int_{-\infty}^{\infty}dx\int_{0}^{\infty}dy\, e^{-ik_1x-i k_2y} q_0(x,y)\label{q0hatgr},\\
\tilde g_t(k_1,k_2)&=\int_{-\infty}^{\infty}dx\int_{0}^{t}d\tau\, e^{-ik_1x+\omega(k) \tau} g(x,\tau)\label{gttilde},\\
\tilde h_t(k_1,k_2)&=\int_{-\infty}^{\infty}dx\int_{0}^{t}d\tau\, e^{-ik_1x+\omega(k) \tau} h(x,\tau)\label{httilde},
\end{align}
\end{subequations}
where
\begin{equation}\label{omega}
\omega(k):=-ik_1^3+3i\frac{k_2^2}{k_1}.
\end{equation}
The functions defined in \eqref{transfs} satisfy the so-called global relation:
\begin{equation}\label{GR1}
e^{\omega(k)t}\hat q(k_1,k_2,t)=\hat q_0(k_1,k_2)+3\left[\frac{k_2}{k_1}\tilde g_t(k_1,k_2)-\frac{i}{k_1} \tilde h_t(k_1,k_2)\right], \ \ k_1\in \mathbb R, \ \mathrm{Im} k_2\leq 0.
\end{equation}

\end{prop}

\begin{PROOF}

The linearised KP equation \eqref{lkpintro} can be written in divergence form as: 
\begin{align}\label{divform}
&\bigg[-ik_1e^{-ik_1x-ik_2y+\omega(k)t}q_x\bigg]_t+3\bigg[e^{-ik_1x-ik_2y+\omega(k)t}\bigg(-ik_1q_y-ik_2q_x\bigg)\bigg]_y\nonumber\\
&+\bigg[e^{-ik_1x-ik_2y+\omega(k)t}\bigg(-ik_1q_{xxx}+ik_1^3q_x+k_1^2q_{xx}+3ik_2q_y\bigg)\bigg]_x=0.
\end{align}
Green's theorem in the plane implies
\begin{align}\label{GR}
\int_{-\infty}^{\infty}dx\int_{0}^{\infty}dy\, \left(e^{-ik_1x-ik_2y+\omega(k)t}q\right)_t\!=\!&\int_{-\infty}^{\infty}dx\, e^{-ik_1x-ik_2y+\omega(k)t}\ 3\left(\partial_{x}^{-1}h+\frac{k_2}{k_1}g\right)\nonumber\\
=&\int_{-\infty}^{\infty}dx\, e^{-ik_1x-ik_2y+\omega(k)t}\ 3\left(\frac{k_2}{k_1}\,g-\frac{i}{k_1}\,h\right).
\end{align}
Integrating equation \eqref{GR}, we obtain the identity
\begin{align}
&\int_{-\infty}^{\infty}d\xi\int_{0}^{\infty}d\eta\, e^{-ik_1\xi-ik_2\eta+\omega(k)t}q=\nonumber\\
=&\int_{-\infty}^{\infty}d\xi\int_{0}^{\infty}d\eta\, e^{-ik_1\xi-ik_2\eta}q_0+\int_{-\infty}^{\infty}d\xi\int_{0}^{t}d\tau\, e^{-ik_1\xi+\omega(k)\tau}3\left(\frac{k_2}{k_1}\,g-\frac{i}{k_1}\,h\right),
\end{align}
valid for $k_1 \in \mathbb R$ and $\mathrm{Im} k_2\leq0$, which is equation \eqref{GR1}.
\QED
\vskip 3mm

\noindent\textbf{Remark 2.1} We also note that if $t$ is replaced by $T\geq t$, equation \eqref{GR1} becomes
\begin{equation}\label{GR2}
e^{\omega(k)T}\hat q(k_1,k_2,T)=\hat q_0(k_1,k_2)+3\!\left[\frac{k_2}{k_1}\tilde g_T(k_1,k_2)-\frac{i}{k_1}\tilde h_T(k_1,k_2)\right]\!, \ k_1\in \mathbb R, \ \mathrm{Im} k_2\leq 0,
\end{equation}
with $\tilde g_T$ and $\tilde h_T$ defined by equations analogous to \eqref{gttilde} and \eqref{httilde}.
\end{PROOF}

\noindent\textbf{Remark 2.2} The global relation can also be obtained for the Lax pair \eqref{Laxintro}. In fact, we have that
\begin{equation*}
\left(\hat \mu e^{l(l+2k)y-4il(l^2+3kl+3k^2)t}\right)_y=e^{l(l+2k)y-4il(l^2+3kl+3k^2)t}\int_{-\infty}^{\infty}dx\, e^{-ilx}q
\end{equation*}
and
\begin{equation*}
\begin{split}
\left(e^{l(l+2k)y-4il(l^2+3kl+3k^2)t}\hat \mu\right)_t&=e^{l(l+2k)y-4il(l^2+3kl+3k^2)t}\\
&\times\int_{-\infty}^{\infty}dx\, e^{-ilx}\Big[-3\left(q_x+2ikq+\partial_x^{-1}h\right)\Big].
\end{split}
\end{equation*}

\begin{figure}[ht]
\begin{center}
\resizebox{4cm}{!}{\input{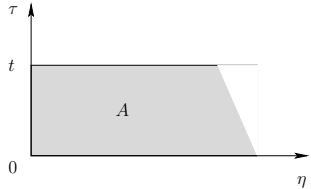}}
\end{center}
\caption{The domain $A$ for the global relation.}
\label{lA1}
\end{figure}

These two equations are compatible iff
\begin{align*}
&\left(e^{l(l+2k)y-4il(l^2+3kl+3k^2)t}\int_{-\infty}^{\infty}dx\, e^{-ilx}q\right)_t=\nonumber\\
&=-\left(e^{l(l+2k)y-4il(l^2+3kl+3k^2)t}\int_{-\infty}^{\infty}dx\, e^{-ilx}\, 3\Big[q_x+2ikq+\partial_x^{-1}h\Big]\right)_y.
\end{align*}
Hence, Green's theorem over a domain $A$ yields
\begin{equation*}
\begin{split}
\int_{\partial A}&\left[e^{l(l+2k)y-4il(l^2+3kl+3k^2)t}\int_{-\infty}^{\infty}dx\, e^{-ilx}q\right]dy=\\
&=\left\{e^{l(l+2k)y-4il(l^2+3kl+3k^2)t}\int_{-\infty}^{\infty}dx\, e^{-ilx}\, 3\bigg[g_\xi+2ikg+\partial_\xi^{-1}h\bigg]\right\}dt.
\end{split}
\end{equation*}
In particular, for the domain depicted in figure \ref{lA1} we find the identity

\begin{align}\label{lgr}
&\int_{-\infty}^{\infty}d\xi\int_{0}^{\infty}d\eta\, e^{-il\xi+l(l+2k)\eta-4il(l^2+3kl+3k^2)t}q=\int_{-\infty}^{\infty}d\xi\int_{0}^{\infty}d\eta\, e^{-il\xi+l(l+2k)\eta}q_0\nonumber\\
&+\int_{-\infty}^{\infty}d\xi\int_{0}^{t}d\tau\, e^{-il\xi-4il(l^2+3kl+3k^2)\tau}\, 3\bigg[i(l+2k)g+\partial_\xi^{-1}h\bigg],\quad l(l+2k_R)\leq 0.
\end{align}
The constraint on $l$ is imposed in order for the exponential $e^{l(l+2k)\eta}$ to be bounded.

\subsection{The inverse problem}\label{lkpip}
We will use Pompeiu's formula in order to reconstruct $\mu$ via a d-bar problem. Since the functions $\mu_1^{\pm}$ and $\mu_{2}^{\pm}$ depend explicitly both on $k$ and on $k_R$, we need to compute the jumps of $\mu$ across the real and the imaginary $k$-axes, as well as the $k$-bar derivatives of $\mu$ in each of the four quadrants. 

\begin{prop}\label{lkpipprop}
Let $\mu$ be defined by equation \eqref{spectralfunctions}. This function admits the following alternative representation in terms of appropriate transforms of $q_0(x,y)$, $g(x,t)$, $h(x,t)$:
\begin{align}\label{newvariables}
\mu&=\frac{1}{2\pi^2}\left(\int_{0}^{\infty}\!\!\!d\nu_R-\int_{-\infty}^{0}\!\!\!d\nu_R\right)\!\int_{-\infty}^{\infty} \frac{d\nu_I}{\nu-k}\, e^{-2i\nu_Rx+4i\nu_R\nu_Iy+8i\nu_R(3\nu_I^2-\nu_R^2)t}\tilde q_0(\nu_R, \nu_I)\nonumber\\
\nonumber\\
&+\frac{1}{2\pi^2}\int_{0}^\infty d\nu_R\int_{\partial D_1^+}\frac{d\nu_I}{\nu-k}\, e^{-2i\nu_R x+4i\nu_R\nu_I y+8i\nu_R(3\nu_I^2-\nu_R^2)t}\tilde S_T(\nu_R, \nu_I)\nonumber\\
\nonumber\\
&-\frac{1}{2\pi^2}\int_{-\infty}^0 d\nu_R\int_{\partial D_1^-}\frac{d\nu_I}{\nu-k}\, e^{-2i\nu_R x+4i\nu_R\nu_I y+8i\nu_R(3\nu_I^2-\nu_R^2)t}\tilde S_T(\nu_R, \nu_I),
\end{align}

\begin{figure}[ht]
\begin{center}
\resizebox{5cm}{!}{\input{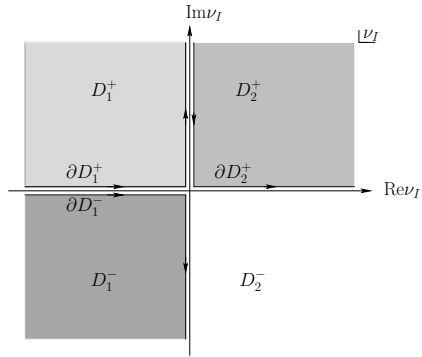}}
\end{center}
\caption{The regions $D^\pm_{1,2}$.}
\label{D1D2}
\end{figure}

where $\tilde q_0$ and $\tilde S_{T}$ are defined by
\begin{align}
\tilde q_0(\nu_R, \nu_I)&= \int_{-\infty}^{\infty}d\xi\int_{0}^{\infty}d\eta\, e^{2i\nu_R\xi-4i\nu_R\nu_I\eta}q_0(\xi, \eta)\label{q0tilde},\\
\tilde S_T(\nu_R, \nu_I)&= \int_{-\infty}^{\infty}d\xi\int_{0}^{T}d\tau \, e^{2i\nu_R\xi-8i\nu_R(3\nu_I^2-\nu_R^2)\tau}\, 3\bigg[\partial^{-1}_{\xi}h(\xi,\tau)-2\nu_Ig(\xi,\tau)\bigg],\label{STtilde}
\end{align}
and the contours $\partial D_1^+$ and $\partial D_1^-$ are the boundaries of the second and the third quadrant, with the orientation as shown in figure \ref{D1D2}.
\end{prop}

\begin{PROOF}
The definitions \eqref{direct1} and \eqref{direct2} of $\mu_1^{\pm}$ and $\mu_2^{\pm}$ yield the following equations:
{\small
\begin{align}\label{direct0+}
&\mu^{\pm}_1\mbox{\large $\Big|$}_{k_R=0}=\frac{1}{2\pi}\mbox{\large $\Big($}\int_{-\infty}^{0}dl+\int_{0}^{\infty}dl\mbox{\large $\Big)$ }\int_{-\infty}^{\infty}d\xi\int_{0}^{y}d\eta\, e^{-il(\xi-x)+l(l+2ik_I)(\eta-y)} q(\xi,\eta,t)\nonumber\\ 
&+\frac{1}{2\pi}\mbox{\large $\Big($}\int_{-\infty}^{0}\!\!\!\!\!dl+\!\!\int_{0}^{\infty}\!\!\!\!\!dl\mbox{\large $\Big)$ }\!\!\!\!\int_{-\infty}^{\infty}\!\!\!\!d\xi \left\lbrace\!\!\!\!\begin{array}{ll}\int_{0}^{t}d\tau \\ \\ -\int_{t}^{T}d\tau \end{array}\right.\!\!\!\!\!\!\! e^{-il(\xi-x)-l(l+2ik_I)y-4il(l^2+3ik_Il-3k_I^2)(\tau-t)} Q(\xi,\tau,ik_I,l),
\end{align}}
{\small
\begin{align}\label{direct0-}
&\mu^{\pm}_2\mbox{\large $\Big|$}_{k_R=0}=\frac{1}{2\pi}\mbox{\large $\Big($}\int_{-\infty}^{0}dl+\int_{0}^{\infty}dl\mbox{\large $\Big)$ }\int_{-\infty}^{\infty}d\xi\int_{0}^{y}d\eta\, e^{-il(\xi-x)+l(l+2ik_I)(\eta-y)}q(\xi,\eta,t)\nonumber\\ 
&+\frac{1}{2\pi}\mbox{\large $\Big($}\int_{-\infty}^{0}\!\!\!\!\!dl+\!\!\int_{0}^{\infty}\!\!\!\!\!dl\mbox{\large $\Big)$ }\!\!\!\!\int_{-\infty}^{\infty}\!\!\!\!d\xi \left\lbrace\!\!\!\!\begin{array}{ll}\int_{0}^{t}d\tau \\ \\ -\int_{t}^{T}d\tau \end{array}\right.\!\!\!\!\!\!\! e^{-il(\xi-x)-l(l+2ik_I)y-4il(l^2+3ik_Il-3k_I^2)(\tau-t)} Q(\xi,\tau,ik_I,l),
\end{align}}
{\small
\begin{align}\label{direct01}
&\mu^{\pm}_1\mbox{\large $\Big|$}_{k_I=0}=\frac{1}{2\pi}\mbox{\large $\Big($}\int_{-\infty}^{0}dl+\int_{-2k_R}^{\infty}dl\mbox{\large $\Big)$ }\int_{-\infty}^{\infty}d\xi\int_{0}^{y}d\eta\, e^{-il(\xi-x)+l(l+2k_R)(\eta-y)}  q(\xi,\eta,t)\nonumber\\ 
&+\frac{1}{2\pi}\mbox{\large $\Big($}\!\int_{-\infty}^{0}\!\!\!\!\!dl+\!\!\int_{-2k_R}^{\infty}\!\!\!\!\!dl\mbox{\large $\Big)$ }\!\!\!\!\!\int_{-\infty}^{\infty}\!\!\!\!d\xi\left\lbrace\!\!\!\!\begin{array}{ll}\int_{0}^{t}d\tau \\ \\ -\int_{t}^{T}d\tau \end{array}\right.\!\!\!\!\!\!\! e^{-il(\xi-x)-l(l+2k_R)y-4il(l^2+3k_Rl+3k_R^2)(\tau-t)} Q(\xi,\tau,k_R,l)\nonumber\\
&-\frac{1}{2\pi}\int_{0}^{-2k_R}dl\int_{-\infty}^{\infty}d\xi\int_{y}^{\infty}d\eta \, e^{-il(\xi-x)+l(l+2k_R)(\eta-y)} q(\xi,\eta,t), \quad k_R\leq0
\end{align}}
and
{\small
\begin{align}\label{direct02}
&\mu^{\pm}_2\mbox{\large $\Big|$}_{k_I=0}= \frac{1}{2\pi}\mbox{\large $\Big($}\int_{-\infty}^{-2k_R}dl+\int_{0}^{\infty}dl\mbox{\large $\Big)$ }\int_{-\infty}^{\infty}d\xi\int_{0}^{y}d\eta\, e^{-il(\xi-x)+l(l+2k_R)(\eta-y)} q(\xi,\eta,t)\nonumber\\ 
&+\frac{1}{2\pi}\mbox{\large $\Big($}\!\int_{-\infty}^{-2k_R}\!\!\!\!\!dl\!\!+\!\!\int_{0}^{\infty}\!\!\!\!\!dl\mbox{\large $\Big)$ }\!\!\!\!\! \int_{-\infty}^{\infty}\!\!\!\!d\xi\left\lbrace\!\!\!\!\begin{array}{ll}\int_{0}^{t}d\tau \\ \\ -\int_{t}^{T}d\tau \end{array}\right.\!\!\!\!\!\!\! e^{-il(\xi-x)-l(l+2k_R)y-4il(l^2+3k_Rl+3k_R^2)(\tau-t)}Q(\xi,\tau,k_R,l)\nonumber\\
&-\frac{1}{2\pi}\int_{-2k_R}^{0}dl\int_{-\infty}^{\infty}d\xi\int_{y}^{\infty}d\eta \, e^{-il(\xi-x)+l(l+2k_R)(\eta-y)} q(\xi,\eta,t),\quad k_R\geq0.
\end{align}}
Equations \eqref{direct0+} and \eqref{direct0-} imply
\begin{equation*}
\left(\mu_1^+-\mu_2^+\right)\bigg|_{k_R=0}=0 \quad \mathrm{and} \quad \left(\mu_1^--\mu_2^-\right)\bigg|_{k_R=0}=0,
\end{equation*}
thus $\mu$ does not have a discontinuity across the imaginary $k$-axis. 

On the other hand, equations \eqref{direct01} and \eqref{direct02} imply 
\begin{equation*}
\left(\mu_1^+-\mu_1^-\right)\bigg|_{k_I=0}= \frac{1}{2\pi}\left(\int_{-\infty}^{0}dl+\int_{-2k_R}^{\infty}dl\right)e^{ilx-l(l+2k_R)y}\, \tilde Q_T(k_R,l)
\end{equation*}
and
\begin{equation*}
\left(\mu_2^+-\mu_2^-\right)\bigg|_{k_I=0}= \frac{1}{2\pi}\left(\int_{-\infty}^{-2k_R}dl+\int_{0}^{\infty}dl\right)e^{ilx-l(l+2k_R)y}\, \tilde Q_T(k_R,l),
\end{equation*}
where the function $\tilde Q_T$ is given by
\begin{equation}\label{QTtilde}
\tilde Q_T(k,l)=\int_{-\infty}^{\infty}d\xi\int_{0}^{T}d\tau\, e^{-il\xi-4il(l^2+3kl+3k^2)\tau}Q(\xi,\tau,k,l),
\end{equation}
with $Q$ defined by equation \eqref{Q}. Hence, $\mu$ is discontinuous across the real $k$-axis. 

Pompeiu's formula for the contour depicted in figure \ref{pompeiu} yields
\begin{figure}[ht]
\begin{center}
\resizebox{4.5cm}{!}{\input{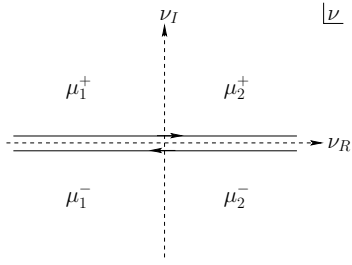}}
\end{center}
\caption{The contour of integration for Pompeiu's formula.}
\label{pompeiu}
\end{figure}
\begin{align}\label{Pomp}
\mu&= \frac{1}{2i\pi}\int_{-\infty}^0 \frac{d\nu_R}{\nu_R-k} \,\left(\mu_1^+-\mu_1^-\right)\!\Big|_{\nu_I=0}+\frac{1}{2i\pi} \int_{0}^{\infty} \frac{d\nu_R}{\nu_R-k}\, \left(\mu_2^+-\mu_2^-\right)\!\Big|_{\nu_I=0}\nonumber\\
&-\frac{1}{\pi}\int_{0}^{\infty}d\nu_R\int_{0}^{\infty} \frac{d\nu_I}{\nu-k}\frac{\partial \mu_2^+}{\partial \bar \nu}-\frac{1}{\pi}\int_{-\infty}^{0}d\nu_R\int_{0}^{\infty}\frac{d\nu_I}{\nu-k}\frac{\partial \mu_1^+}{\partial \bar \nu}\nonumber\\
&-\frac{1}{\pi}\int_{0}^{\infty}d\nu_R\int_{-\infty}^{0}\frac{d\nu_I}{\nu-k}\frac{\partial \mu_2^-}{\partial \bar \nu}-\frac{1}{\pi}\int_{-\infty}^{0}d\nu_R\int_{-\infty}^{0}\frac{d\nu_I}{\nu-k}\frac{\partial \mu_1^-}{\partial \bar \nu}.
\end{align}
In order to compute the associated d-bar derivatives, we recall the definition
\begin{equation}\label{d-bardef}
\frac{\partial \mu}{\partial \bar \nu}=\, \frac{1}{2}\left(\frac{\partial}{\partial \nu_R}+i\frac{\partial}{\partial \nu_I}\right)\mu. 
\end{equation}
Thus, by differentiating equation \eqref{direct1}, we find
\begin{align*}
\frac{\partial \mu_{1}^{\pm}}{\partial \bar \nu}&= \frac{1}{2\pi}\int_{-\infty}^{\infty}d\xi\int_0^{\infty}d\eta\, e^{2i\nu_R(\xi-x)-4i\nu_R\nu_I(\eta-y)}q(\xi,\eta,t)\nonumber\\
&+\frac{1}{2\pi}\int_{-\infty}^{\infty}d\xi\left\{\begin{array}{ll} \int_{0}^{t}d\tau \\ \\ -\int_{t}^{T}d\tau \end{array}\right.\!\!\! e^{2i\nu_R (\xi-x)+4i\nu_R \nu_I y-8i\nu_R(3\nu_I^2-\nu_R^2)(\tau-t)}Q(\xi,\tau,\nu,-2\nu_R).
\end{align*}
Similarly, using equation \eqref{direct2}, we find 
\begin{align*}
\frac{\partial \mu_{2}^{\pm}}{\partial \bar \nu}&=-\frac{1}{2\pi}\int_{-\infty}^{\infty}d\xi\int_0^{\infty}d\eta\, e^{2i\nu_R(\xi-x)-4i\nu_R\nu_I(\eta-y)}q(\xi,\eta,t)\nonumber\\
&-\frac{1}{2\pi}\int_{-\infty}^{\infty}d\xi\left\{\begin{array}{ll} \int_{0}^{t}d\tau \\ \\ -\int_{t}^{T}d\tau \end{array}\right.\!\!\! e^{2i\nu_R (\xi-x)+4i\nu_R \nu_I y-8i\nu_R(3\nu_I^2-\nu_R^2)(\tau-t)}Q(\xi,\tau,\nu,-2\nu_R).
\end{align*}
Thus, the d-bar derivatives of $\mu_1^\pm$ and $\mu_2^\pm$ are equal, i.e.,
\begin{equation*}
\frac{\partial \mu_{2}^{\pm}}{\partial \bar \nu}=-\frac{\partial \mu_{1}^{\pm}}{\partial \bar \nu}.
\end{equation*}
Moreover, employing the global relation \eqref{lgr} together with the identity
\begin{equation}\label{tT}
 -\int_{t}^{T}=\int_{0}^{t}-\int_{0}^{T},
\end{equation}
we can express the above d-bar derivatives in terms of the initial and boundary values of $q$ and its derivatives:
\begin{subequations} 
\begin{align*}
\frac{\partial \mu_{1}^+}{\partial \bar \nu}&= \frac{1}{2\pi}\, e^{-2i\nu_Rx+4i\nu_R\nu_Iy+8i\nu_R(3\nu_I^2-\nu_R^2)t}\tilde q_0(\nu_R, \nu_I),\\
\frac{\partial \mu_{2}^+}{\partial \bar \nu}&=-\frac{1}{2\pi}\, e^{-2i\nu_Rx+4i\nu_R\nu_Iy+8i\nu_R(3\nu_I^2-\nu_R^2)t}\tilde q_0(\nu_R, \nu_I)
\end{align*}
and
\begin{align*}
\frac{\partial \mu_{1}^-}{\partial \bar \nu}&=\frac{1}{2\pi}\, e^{-2i\nu_Rx+4i\nu_R\nu_Iy+8i\nu_R(3\nu_I^2-\nu_R^2)t}\left[\tilde q_0(\nu_R, \nu_I)+\tilde S_T(\nu_R, \nu_I)\right],\\
\frac{\partial \mu_{2}^-}{\partial \bar \nu}&=-\frac{1}{2\pi}\, e^{-2i\nu_Rx+4i\nu_R\nu_Iy+8i\nu_R(3\nu_I^2-\nu_R^2)t}\left[\tilde q_0(\nu_R, \nu_I)+\tilde S_T(\nu_R, \nu_I)\right],
\end{align*}
\end{subequations}
where $q_0$ and $\tilde S_T$ are defined by equations \eqref{q0tilde} and \eqref{STtilde} respectively.

Thus, inserting the above expressions in Pompeiu's formula \eqref{Pomp}, we obtain the following expression for $\mu$:
\begin{align}\label{oldvariables}
\mu&= \frac{1}{4i\pi^2}\int_{-\infty}^0 \! \frac{d\nu_R}{\nu_R-k}\left(\int_{-\infty}^{0}\!dl+\!\int_{-2\nu_R}^{\infty}dl\right)e^{ilx-l(l+2\nu_R)y+4il(l^2+3\nu_R l+3\nu_R^2)t} \tilde Q_T(\nu_R,l)\nonumber\\
&+\, \frac{1}{4i\pi^2}\int_{0}^{\infty}\!  \frac{d\nu_R}{\nu_R-k}\left(\int_{-\infty}^{-2\nu_R}\!dl+\!\int_{0}^{\infty}dl\right)e^{ilx-l(l+2\nu_R)y+4il(l^2+3\nu_R l+3\nu_R^2)t} \tilde Q_T(\nu_R,l)\nonumber\\
&+\frac{1}{2\pi^2}\left(\int_{0}^{\infty}d\nu_R-\int_{-\infty}^{0}d\nu_R\right)\int_{-\infty}^{\infty} \frac{d\nu_I}{\nu-k}\, e^{-2i\nu_Rx+4i\nu_R\nu_Iy+8i\nu_R(3\nu_I^2-\nu_R^2)t}\tilde q_0(\nu_R, \nu_I)\nonumber\\
&+\frac{1}{2\pi^2}\left(\int_{0}^{\infty}\!d\nu_R-\int_{-\infty}^{0}\!d\nu_R\right)\!\!\int_{-\infty}^{0}\frac{d\nu_I}{\nu-k}\, e^{-2i\nu_Rx+4i\nu_R\nu_Iy+8i\nu_R(3\nu_I^2-\nu_R^2)t}\tilde S_T(\nu_R, \nu_I).
\end{align}
In the first two integrals on the RHS of equation \eqref{oldvariables}, we make the change of variables
\begin{align}\label{change1}
\left.\begin{array}{rl}l&=2\alpha\\ \nu_R&=\beta-\alpha\end{array}\right\}\Leftrightarrow\left\{\begin{array}{rl}\alpha&=\frac{l}{2}\\ \beta&=\nu_R+\frac{l}{2}.\end{array}\right.
\end{align}
The Jacobian of this transformation is given by 
\begin{equation*}
\frac{\partial (\nu_R, l)}{\partial(\beta, \alpha)}=\frac{\partial \nu_R}{\partial \beta}\frac{\partial l}{\partial \alpha}-\frac{\partial \nu_R}{\partial \alpha}\frac{\partial l}{\partial \beta}=2.
\end{equation*}
The domains of integration are transformed as shown in figure \ref{transform1}; the first two terms on the RHS of equation \eqref{oldvariables} become
\begin{figure}[ht]
\begin{center}
\resizebox{7.5cm}{!}{\input{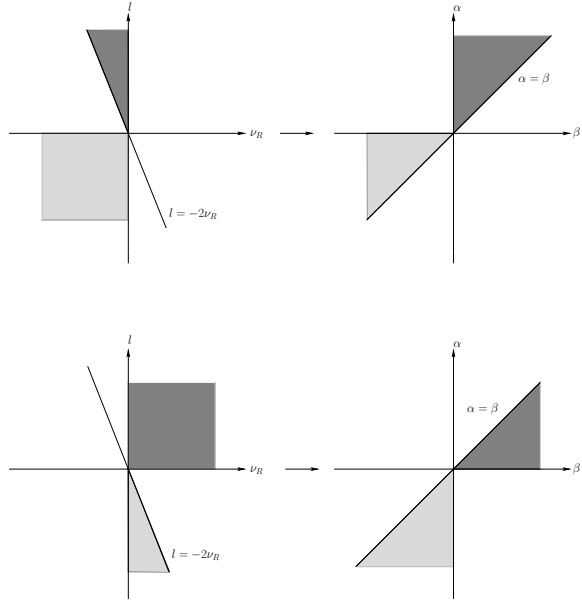}}
\end{center}
\caption{The regions of integration before and after the change of variables \eqref{change1}.}
\label{transform1}
\end{figure}
{\small
\begin{align*}
&2\ \mbox{\large $\Big($}\int_{-\infty}^0\!\!\!d\beta\!\int_{\beta}^{0}\!\frac{d\alpha}{\beta-\alpha-k}+\!\int_{0}^\infty \!\!\!d\beta\!\int_{\beta}^{\infty} \!\!\!\frac{d\alpha}{\beta-\alpha-k}\mbox{\large $\Big)$}e^{2i\alpha x-4\alpha \beta y+8i\alpha (\alpha^2+3\beta^2)t}\tilde Q_T(2\alpha,\beta-\alpha)\\
+\ &2\ \mbox{\large $\Big($}\int_{-\infty}^{0}\!\!\!d\beta\!\int_{-\infty}^{\beta}\!\frac{d\alpha}{\beta-\alpha-k}+\!\int_{0}^{\infty}\!\!\!d\beta\!\int_{0}^{\beta}\!\!\!\!\frac{d\alpha}{\beta-\alpha-k} \mbox{\large $\Big)$}e^{2i\alpha x-4\alpha \beta y+8i\alpha(\alpha^2+3\beta^2)t}\tilde Q_T(2\alpha,\beta-\alpha).
\end{align*}}
We next introduce the transformation 
\begin{align*}
\left.\begin{array}{rl}\beta&=i\nu_I\\ \alpha&=-\nu_R\end{array}\right\}\Leftrightarrow\left\{\begin{array}{rl}\nu_I&=-i\beta\\ \nu_R&=-\alpha. \end{array}\right.
\end{align*}
Then, integrating by parts and recalling definitions \eqref{STtilde} and \eqref{QTtilde}, we find
\begin{align*}
\tilde Q_T(\beta-\alpha, 2\alpha)&=\int_{-\infty}^{\infty}d\xi \int_{0}^{T}d\tau\, e^{-2i\alpha\xi-8i\alpha(\alpha^2+3\beta^2)\tau}Q(\xi,\tau,\beta-\alpha,2\alpha)\\
&=\int_{-\infty}^{\infty}d\xi\int_{0}^{T}d\tau\, e^{2i\nu_R\xi-8i\nu_R(3\nu_I^2-\nu_R^2)\tau}\, 3\bigg[2\nu_I g(\xi, \tau)-\partial_{\xi}^{-1}h(\xi,\tau)\bigg]\nonumber\\
&=-\tilde S_T(\nu_R,\nu_I).
\end{align*}
Thus, the first two integrals on the RHS of equation \eqref{oldvariables} are equal to the following expression:
{\small
\begin{align*}
&-2i\mbox{\large $\Big($}\!\int_0^{i\infty}\!\!\!d\nu_I\!\int_{-i\nu_I}^{0}\!\frac{d\nu_R}{\nu-k}+\!\!\int_{-i\infty}^0 \!\!\!d\nu_I\!\int_{-i\nu_I}^{-\infty}\!\frac{d\nu_R}{\nu-k} \mbox{\large $\Big)$}e^{-2i\nu_R x+4i\nu_R\nu_I y+8i\nu_R (3\nu_I^2-\nu_R^2)t}\tilde S_T(\nu_R,\nu_I)\\
&-2i\mbox{\large $\Big($}\!\int_0^{i\infty}\!\!\!d\nu_I\!\!\int_{\infty}^{-i\nu_I}\!\!\frac{d\nu_R}{\nu-k}+\!\!\int_{-i\infty}^0 \!\!\!d\nu_I\!\int_{0}^{-i\nu_I}\!\!\!\frac{d\nu_R}{\nu-k} \mbox{\large $\Big)$}e^{-2i\nu_R x+4i\nu_R\nu_I y+8i\nu_R (3\nu_I^2-\nu_R^2)t}\tilde S_T(\nu_R,\nu_I)\\
&= 2i\mbox{\large $\Big($}\int_0^{i\infty}\!\!\!d\nu_I\int_{0}^{\infty}\!\!\!\frac{d\nu_R}{\nu-k}-\!\!\int_0^{-i\infty}\!\!\!\!d\nu_I\int_{-\infty}^0\frac{d\nu_R}{\nu-k}\mbox{\large $\Big)$}\, e^{-2i\nu_R x+4i\nu_R\nu_I y+8i\nu_R (3\nu_I^2-\nu_R^2)t}\tilde S_T(\nu_R,\nu_I).
\end{align*}}
Using this expression into equation \eqref{oldvariables}, we find equation \eqref{newvariables}.
\QED
\end{PROOF}

\noindent\textbf{Remark 2.3} In equation \eqref{newvariables}, $\tilde S_T$ can be replaced by $\tilde S_t$, where
\begin{equation}\label{Sttilde}
\tilde S_t(\nu_R, \nu_I)= \int_{-\infty}^{\infty}d\xi\int_{0}^{t}d\tau \, e^{2i\nu_R\xi-8i\nu_R(3\nu_I^2-\nu_R^2)\tau}\, 3\bigg[\partial^{-1}_{\xi}h(\xi,\tau)-2\nu_Ig(\xi,\tau)\bigg].
\end{equation}
Indeed, we have that
\begin{align*}
\mathrm{Re}(4i\nu_R\nu_Iy)\leq 0 &\Leftrightarrow\left\{\begin{array}{ll} \arg \nu_I \in \left[0,\pi\right], \ \nu_R\geq 0 \\ \\ \arg \nu_I \in \left[\pi,2\pi\right], \ \nu_R\leq 0  \end{array} \right.
\end{align*}
and
\begin{align*}
\mathrm{Re}\bigg[-8i\nu_R(3\nu_I^2-\nu_R^2)(t-\tau)\bigg]\leq 0 &\Leftrightarrow\left\{\begin{array}{ll} \arg \nu_I\in \left[\dfrac{\pi}{2},\pi\right]\cup \left[\dfrac{3\pi}{2},2\pi\right], \ \nu_R\geq 0 \\ \\ \arg \nu_I\in \left[0,\dfrac{\pi}{2}\right]\cup \left[\pi,\dfrac{3\pi}{2}\right], \ \nu_R\leq 0.  \end{array} \right.
\end{align*}
Hence, the exponential involved in the integrals along the contours $\partial D_1^+$ and $\partial D_1^-$ is bounded and goes to zero as $|\nu|\rightarrow \infty$ for $\nu\in\{\nu_R\geq 0,\, \nu_I \in D_1^+\}$, and for $\nu\in\{\nu_R\leq 0,\, \nu_I\in D_1^-\}$. Hence, by Cauchy's theorem and Jordan's lemma in the regions $D_1^+$ and $D_1^-$ the term which involves the integral $\int_{t}^{T}$ vanishes.

\subsection{An integral representation for $q(x,y,t)$}\label{lkpir}

The representation \eqref{lkpipprop} for $\mu$ suggests that 
\begin{equation*}
\mu(x,y,t,k_R,k_I)=O\left(\frac{1}{k}\right),\quad k \rightarrow \infty.
\end{equation*}
Consequently, equation \eqref{Lax1intro} allows us to write $q$ in terms of $\mu$:
\begin{equation*}
q(x,y,t)=-2i\lim_{k\rightarrow \infty}\bigg[k\mu_x(x,y,t,k_R,k_I)\bigg].
\end{equation*}
Hence, by differentiating equations \eqref{lkpipprop} we find:
{\small
\begin{align}\label{finalrepT}
q(x,y,t)&= \frac{2}{\pi^2}\mbox{\large $\Big($}\int_{0}^{\infty}d\nu_R-\int_{-\infty}^{0}d\nu_R\mbox{\large $\Big)$}\!\!\int_{-\infty}^{\infty}d\nu_I\, e^{-2i\nu_R x+4i\nu_R\nu_I y+8i\nu_R(3\nu_I^2-\nu_R^2)t}\nu_R\, \tilde q_0(\nu_R, \nu_I)\nonumber\\
&+\frac{2}{\pi^2}\int_{0}^\infty d\nu_R\int_{\partial D_1^+}d\nu_I\, e^{-2i\nu_R x+4i\nu_R\nu_I y+8i\nu_R(3\nu_I^2-\nu_R^2)t}\nu_R\, \tilde S_T(\nu_R, \nu_I)\nonumber\\
&-\frac{2}{\pi^2}\int_{-\infty}^0 d\nu_R\int_{\partial D_1^-}d\nu_I\, e^{-2i\nu_R x+4i\nu_R\nu_I y+8i\nu_R(3\nu_I^2-\nu_R^2)t}\nu_R\, \tilde S_T(\nu_R, \nu_I),
\end{align}}
where $\tilde q_0$ and $\tilde S_T$ are defined by equations \eqref{q0tilde} and \eqref{STtilde} respectively.

According to remark 2.3 we can replace $\tilde S_T$ by $\tilde S_t$, thus an alternative representation for $q$ is given by
{\small
\begin{align}\label{finalrept}
q(x,y,t)&= \frac{2}{\pi^2}\mbox{\large $\Big($}\int_{0}^{\infty}d\nu_R-\int_{-\infty}^{0}d\nu_R\mbox{\large $\Big)$}\!\!\int_{-\infty}^{\infty}d\nu_I\, e^{-2i\nu_R x+4i\nu_R\nu_I y+8i\nu_R(3\nu_I^2-\nu_R^2)t}\nu_R\, \tilde q_0(\nu_R, \nu_I)\nonumber\\
&+\frac{2}{\pi^2}\int_{0}^\infty d\nu_R\int_{\partial D_1^+}d\nu_I\, e^{-2i\nu_R x+4i\nu_R\nu_I y+8i\nu_R(3\nu_I^2-\nu_R^2)t}\nu_R\, \tilde S_t(\nu_R, \nu_I)\nonumber\\
&-\frac{2}{\pi^2}\int_{-\infty}^0 d\nu_R\int_{\partial D_1^-}d\nu_I\, e^{-2i\nu_R x+4i\nu_R\nu_I y+8i\nu_R(3\nu_I^2-\nu_R^2)t}\nu_R\,\tilde S_t(\nu_R, \nu_I),
\end{align}}
where $\tilde q_0$ and and $\tilde S_t$ are defined by equations \eqref{q0tilde} and \eqref{Sttilde} respectively. 

Let us introduce the variables $k_1, k_2$ by the equations
\begin{align*}
\left.\begin{array}{rl}k_1=-2\nu_R\\ \\ k_2=4\nu_R\nu_I\end{array}\right\}\Leftrightarrow \left\{\begin{array}{rl}\nu_R=-\dfrac{k_1}{2}\\ \\ \nu_I=-\dfrac{k_2}{2k_1}.\end{array}\right.
\end{align*}
The Jacobian of this transformation is given by
\begin{equation*}
\left|\frac{\partial (\nu_R, \nu_I)}{\partial(k_1,k_2)}\right|=\left|\frac{\partial \nu_R}{\partial k_1}\frac{\partial \nu_I}{\partial k_2}-\frac{\partial \nu_R}{\partial k_2}\frac{\partial \nu_I}{\partial k_1}\right|=\frac{1}{4|k_1|}.
\end{equation*} 
Thus, equation \eqref{finalrept} becomes
\begin{align}\label{ktransform2}
q(x,y,t)&= \frac{1}{(2\pi)^2}\left(\int_{-\infty}^0\!dk_1\int_{-\infty}^{\infty}\!dk_2+\int_0^{\infty}\!dk_1\int_{-\infty}^{\infty}\!dk_2\right)\, e^{ik_1 x+ik_2 y-\omega(k)t}\hat q_0(k_1,k_2)\nonumber\\
&+\frac{1}{(2\pi)^2}\int_{-\infty}^0 dk_1\int_{\partial D_1^+}dk_2\, e^{ik_1 x+ik_2 y-\omega(k)t}\, \tilde S_t(-\frac{k_1}{2}, -\frac{k_2}{2k_1})\nonumber\\
&+\frac{1}{(2\pi)^2}\int_0^{\infty} dk_1\int_{\partial D_1^-}dk_2\, e^{ik_1 x+ik_2 y-\omega(k)t}\, \tilde S_t(-\frac{k_1}{2}, -\frac{k_2}{2k_1}),
\end{align}
where $\hat q_0$ and $\tilde S_t$ are defined by equations \eqref{q0hatgr} and \eqref{Sttilde}, $\omega(k)$ is given by equation \eqref{omega} and the contours of integration are shown in figure \ref{D1D2}.

The real part of $\omega(k)$ is equal to $-\mathrm{Re}k_2\, \mathrm{Im}k_2/k_1$ and, since $\tau-t\leq 0$, the exponential term in equation \eqref{ktransform2} is bounded in the domains $\{k_1\geq 0,\ k_2 \in D_1^+\cup D_2^-\}$ and $\{k_1\leq 0,\ k_2 \in D_2^+\cup D_1^-\}$. Hence, the contour $\partial D_1^-$ can be deformed onto $\partial D_2^+$ (see figure \ref{D1D2}).

Then, equation \eqref{ktransform2} implies
\begin{align*}
q(x,y,t)&= \frac{1}{(2\pi)^2}\left(\int_{-\infty}^0\!dk_1\int_{-\infty}^{\infty}dk_2+\int_0^{\infty}\!dk_1\int_{-\infty}^{\infty}dk_2\right)\, e^{ik_1 x+ik_2 y-\omega(k)t}\hat q_0(k_1,k_2)\nonumber\\
&+\frac{1}{(2\pi)^2}\int_{-\infty}^0 dk_1\int_{\partial D_1^+}dk_2\, e^{ik_1 x+ik_2 y-\omega(k)t}\, \tilde S_t(-\frac{k_1}{2}, -\frac{k_2}{2k_1})\nonumber\\
&+\frac{1}{(2\pi)^2}\int_0^{\infty} dk_1\int_{\partial D_2^+}dk_2\, e^{ik_1 x+ik_2 y-\omega(k)t}\, \tilde S_t(-\frac{k_1}{2}, -\frac{k_2}{2k_1}).
\end{align*}
Integration by parts in the definition \eqref{Sttilde} of $\tilde S_t$ yields that
\begin{align*}
\tilde S_t(-\frac{k_1}{2}, -\frac{k_2}{2k_1})=\int_{-\infty}^{\infty}d\xi\int_{0}^{t}d\tau\, e^{-ik_1\xi+\omega(k)\tau} 3\left[\frac{k_2}{k_1}\,g(\xi,\tau)-\frac{i}{k_1}\,h(\xi,\tau)\right],
\end{align*}
hence, recalling the definitions \eqref{gttilde} and \eqref{httilde}, we find the following expression:
\begin{align}\label{Fsolutiont}
q(x,y,t)&= \frac{1}{(2\pi)^2}\int_{-\infty}^{\infty}dk_1\int_{-\infty}^{\infty}dk_2\, e^{ik_1 x+ik_2 y-\omega(k)t}\hat q_0(k_1,k_2)\nonumber\\
&+\frac{1}{(2\pi)^2}\int_{-\infty}^0\!\!\! dk_1\int_{\partial D_1^+}\!dk_2\, e^{ik_1 x+ik_2 y-\omega(k)t}\,3\left[\frac{k_2}{k_1}\,\tilde g_t(k_1,k_2)-\frac{i}{k_1}\,\tilde h_t(k_1,k_2)\right]\nonumber\\
&+\frac{1}{(2\pi)^2}\int_0^{\infty}\!\!\! dk_1\int_{\partial D_2^+}dk_2\, e^{ik_1 x+ik_2 y-\omega(k)t}\,3\left[\frac{k_2}{k_1}\,\tilde g_t(k_1,k_2)-\frac{i}{k_1}\,\tilde h_t(k_1,k_2)\right],
\end{align}
where $\omega$ is defined by equation \eqref{Omegaintro}, $\hat q_0$, $\tilde g_t$ and $\tilde h_t$ are defined by equations \eqref{q0hatgr}, \eqref{gttilde} and \eqref{httilde} respectively, and the contours $\partial D_1^+$ and $\partial D_2^+$ are shown in figure \ref{D1D2}.

The solution $q$ of an IBV problem for the linearised KP in the domain $\Omega$ is given by equation \eqref{Fsolutiont}.

According to remark 2.3, an equivalent formula for $q$ is given by
\begin{align}\label{FsolutionT}
q(x,y,t)&= \frac{1}{(2\pi)^2}\int_{-\infty}^{\infty}dk_1\int_{-\infty}^{\infty}dk_2\, e^{ik_1 x+ik_2 y-\omega(k)t}\hat q_0(k_1,k_2)\nonumber\\
&+\frac{1}{(2\pi)^2}\int_{-\infty}^0\!\!\!\! dk_1\int_{\partial D_1^+}\!dk_2\, e^{ik_1 x+ik_2 y-\omega(k)t}\,3\left[\frac{k_2}{k_1}\,\tilde g_T(k_1,k_2)-\frac{i}{k_1}\,\tilde h_T(k_1,k_2)\right]\nonumber\\
&+\frac{1}{(2\pi)^2}\int_0^{\infty}\!\!\! dk_1\int_{\partial D_2^+}\!dk_2\, e^{ik_1 x+ik_2 y-\omega(k)t}\,3\left[\frac{k_2}{k_1}\,\tilde g_T(k_1,k_2)-\frac{i}{k_1}\,\tilde h_T(k_1,k_2)\right]\!\!,
\end{align}
with $\tilde g_T$ and $\tilde h_T$ defined analogously to $\tilde g_t$ and $\tilde h_t$.

\begin{prop}
The formula defined by equation \eqref{FsolutionT} in terms of $\{q_0(x,y)$, $q(x,0,t)$, $q_y(x,0,t)\}$ represents the solution to the linearised KP equation \eqref{lkpintro} in $\Omega$, provided that the functions $\hat q_0$, $\tilde g_T$, $\tilde h_T$ satisfy the global relation \eqref{GR2}. An equivalent representation is given by equation \eqref{Fsolutiont}. 
\end{prop}
\begin{PROOF}
Since equations \eqref{Fsolutiont} and \eqref{FsolutionT} are equivalent, we can use either of them in the calculations that follow. 

By differentiating equation \eqref{FsolutionT}, it is straightforward to show that $q(x,y,t)$ satisfies the linearised KP equation.

\underline{Initial Condition}: Evaluating the expression \eqref{Fsolutiont} at $t=0$ and using the fact that
\begin{equation}\label{delta}
 \int_{\mathcal{D}}d\lambda\, e^{i\lambda(x-\xi)}=2\pi \delta(x-\xi),\quad \lambda \in \mathcal D,
\end{equation}
we find
\begin{align*}
q(x,y,0)&= \frac{1}{(2\pi)^2}\int_{-\infty}^{\infty}dk_1\int_{-\infty}^{\infty}dk_2\, e^{ik_1x+ik_2y}\hat q_0(k_1,k_2)\nonumber\\
&= \frac{1}{(2\pi)^2}\int_{-\infty}^{\infty}dk_1\int_{-\infty}^{\infty}dk_2\, \int_{-\infty}^{\infty}d\xi\int_{0}^{\infty}d\eta\, e^{ik_1(x-\xi)+ik_2(y-\eta)}q_0(\xi, \eta)\nonumber\\
&= \frac{1}{(2\pi)^2}\int_{-\infty}^{\infty}d\xi\int_{0}^{\infty}d\eta\, (2\pi)^2\delta(x-\xi)\, \delta(y-\eta)\,q_0(\xi, \eta)\nonumber\\
&= q_0(x,y),
\end{align*}
thus the initial condition \eqref{icintro} is satisfied.

\underline{Boundary Values}: At $y=0$, equation \eqref{FsolutionT} becomes
\begin{align}\label{FsolutionT0}
q(x,0,t)&= \frac{1}{(2\pi)^2}\int_{-\infty}^{\infty}dk_1\int_{-\infty}^{\infty}dk_2\, e^{ik_1 x-\omega(k)t}\hat q_0(k_1,k_2)\nonumber\\
&+\frac{1}{(2\pi)^2}\int_{-\infty}^0\!\!\!\! dk_1\int_{\partial D_1^+}\!dk_2\, e^{ik_1 x-\omega(k)t}\,3\left[\frac{k_2}{k_1}\,\tilde g_T(k_1,k_2)-\frac{i}{k_1}\,\tilde h_T(k_1,k_2)\right]\nonumber\\
&+\frac{1}{(2\pi)^2}\int_0^{\infty}\!\!\! dk_1\int_{\partial D_2^+}\!dk_2\, e^{ik_1 x-\omega(k)t}\,3\left[\frac{k_2}{k_1}\,\tilde g_T(k_1,k_2)-\frac{i}{k_1}\,\tilde h_T(k_1,k_2)\right].
\end{align}
Under the summetry transformation $k_2\mapsto-k_2$ of $\omega(k)$, the global relation \eqref{GR2} assumes the form:
{\small\begin{equation}\label{GR1sym}
e^{\omega(k)t}\hat q(k_1,-k_2,t)=\hat q_0(k_1,-k_2)-3\,\frac{k_2}{k_1}\tilde g_T(k_1,k_2)-3\,\frac{i}{k_1} \tilde h_T(k_1,k_2), \ \ k_1\in \mathbb R, \ \mathrm{Im} k_2\geq 0,
\end{equation}}
while we notice that the functions $\tilde g_T$ and $\tilde h_T$, which depend on $k_2$ only through $\omega(k)$, remain invariant. Solving this expression for $h_T$ and substituting into equation \eqref{FsolutionT0}, we find
\begin{align}\label{FsolutionT0b}
q(x,0,t)&= \frac{1}{(2\pi)^2}\int_{-\infty}^{\infty}dk_1\int_{-\infty}^{\infty}dk_2\, e^{ik_1 x-\omega(k)t}\hat q_0(k_1,k_2)\nonumber\\
&+\frac{1}{(2\pi)^2}\left(\int_{-\infty}^0\!\!\!\! dk_1\int_{\partial D_1^+}\!\!\!dk_2+\int_0^{\infty}\!\!\! dk_1\int_{\partial D_2^+}\!\!\!dk_2\right) e^{ik_1 x-\omega(k)t}\,6\,\frac{k_2}{k_1}\,\tilde g_T(k_1,k_2)\nonumber\\
&-\frac{1}{(2\pi)^2}\left(\int_{-\infty}^0\!\!\!\! dk_1\int_{\partial D_1^+}\!dk_2+\int_0^{\infty}\!\!\! dk_1\int_{\partial D_2^+}\!dk_2\right) e^{ik_1 x-\omega(k)t}\hat q_0(k_1,-k_2)\nonumber\\
&+\frac{1}{(2\pi)^2}\left(\int_{-\infty}^0\!\!\!\! dk_1\int_{\partial D_1^+}\!dk_2+\int_0^{\infty}\!\!\! dk_1\int_{\partial D_2^+}\!dk_2\right) e^{ik_1 x}\hat q(k_1,-k_2,t).
\end{align}
Recalling the definition \eqref{qhatgr} of $\hat q$, we note that the last term on the RHS of the above equation vanishes due to Jordan's lemma on the complex $k_2$-plane. Also, recalling the definition \eqref{q0hatgr} of $q_0$, equation \eqref{FsolutionT0b} becomes
\begin{align}\label{FsolutionT0c}
&q(x,0,t)= \frac{1}{(2\pi)^2}\int_{-\infty}^{\infty}dk_1\int_{-\infty}^{\infty}dk_2\, e^{ik_1 x-\omega(k)t}\int_{-\infty}^{\infty}d\xi\int_{0}^{\infty}d\eta\, e^{-ik_1\xi-i k_2\eta} q_0(\xi,\eta)\nonumber\\
&+\frac{1}{(2\pi)^2}\left(\!\int_{-\infty}^0\!\!\!\!\!\! dk_1\!\!\int_{\partial D_1^+}\!\!\!\!\!\!dk_2+\!\!\int_0^{\infty}\!\!\!\!\! dk_1\!\!\int_{\partial D_2^+}\!\!\!\!\!\!dk_2\right)\!\! e^{ik_1 x-\omega(k)t}\,6\,\frac{k_2}{k_1}\!\!\int_{-\infty}^{\infty}\!\!\!\!\!\!d\xi\!\int_{0}^{T}\!\!\!\!\!\!d\tau\, e^{-ik_1\xi+\omega(k) \tau} g(\xi,\tau)\nonumber\\
&-\frac{1}{(2\pi)^2}\left(\int_{-\infty}^0\!\!\!\!\!\!dk_1\int_{\partial D_1^+}\!\!\!\!\!dk_2+\int_0^{\infty}\!\!\!\!\! dk_1\int_{\partial D_2^+}\!\!\!\!\!dk_2\right) e^{ik_1 x-\omega(k)t}\!\!\int_{-\infty}^{\infty}\!\!\!\!\!\!d\xi\!\int_{0}^{\infty}\!\!\!\!\!\!d\eta\, e^{-ik_1\xi+i k_2\eta} q_0(\xi,\eta).
\end{align}
The definition \eqref{omega} of $\omega$ implies
\begin{equation*}
 \mathrm{Re}\,\omega(k)= -\frac{6}{k_1}\, \mathrm{Re}k_2\,\mathrm{Im}k_2,
\end{equation*}
hence we can deform the contours $\partial D_1^+$ and $\partial D_2^+$ onto the real $k$-axis in the third term on the RHS of equation \eqref{FsolutionT0c}, i.e.
\begin{align*}
&q(x,0,t)= \frac{1}{(2\pi)^2}\int_{-\infty}^{\infty}dk_1\int_{-\infty}^{\infty}dk_2\, e^{ik_1 x-\omega(k)t}\int_{-\infty}^{\infty}d\xi\int_{0}^{\infty}d\eta\, e^{-ik_1\xi-i k_2\eta} q_0(\xi,\eta)\nonumber\\
&+\frac{1}{(2\pi)^2}\left(\!\int_{-\infty}^0\!\!\!\!\!\! dk_1\!\!\int_{\partial D_1^+}\!\!\!\!\!\!dk_2+\!\!\int_0^{\infty}\!\!\!\!\! dk_1\!\!\int_{\partial D_2^+}\!\!\!\!\!\!dk_2\right)\!\! e^{ik_1 x-\omega(k)t}\,6\,\frac{k_2}{k_1}\!\!\int_{-\infty}^{\infty}\!\!\!\!\!\!d\xi\!\int_{0}^{T}\!\!\!\!\!\!d\tau\, e^{-ik_1\xi+\omega(k) \tau} g(\xi,\tau)\nonumber\\
&-\frac{1}{(2\pi)^2}\int_{-\infty}^\infty dk_1\int_{-\infty}^{\infty}dk_2\,e^{ik_1 x-\omega(k)t}\int_{-\infty}^{\infty}d\xi\int_{0}^{\infty}d\eta\, e^{-ik_1\xi+i k_2\eta} q_0(\xi,\eta).
\end{align*}
The first and the third term on the RHS of this equation cancel, while the change of variables
\begin{equation*}
 l=-3 \frac{k_2^2}{k_1} \Rightarrow dl=-6\frac{k_2}{k_1}
\end{equation*}
on the second term yields
\begin{align}\label{FsolutionT0d}
&q(x,0,t)= \frac{1}{(2\pi)^2}\int_{-\infty}^{\infty}\!\!\!\!dk_1\int_{-\infty}^{\infty}\!\!\!\!dl\, e^{ik_1x-(-ik_1^3+il)t}\int_{-\infty}^{\infty}\!\!\!\!\!d\xi\int_{0}^{T}\!\!\!\!\!d\tau\, e^{-ik_1\xi+(-ik_1^3+il) \tau} g(\xi,\tau).
\end{align}
Using the identity \eqref{delta}, we find
\begin{align}\label{FsolutionT0e}
q(x,0,t)&= \frac{1}{(2\pi)^2}\int_{-\infty}^{\infty}\!\!\!d\xi\int_{-\infty}^{\infty}\!\!dk_1e^{-ik_1(\xi-x)}\int_{0}^{T}\!\!\!d\tau\, e^{-ik_1^3(\tau-t)}\int_{-\infty}^{\infty}\!\!dl\, e^{il( \tau-t)} g(\xi,\tau)
\nonumber\\
&= \frac{1}{2\pi}\int_{-\infty}^{\infty}\!\!\!d\xi\int_{-\infty}^{\infty}\!\!dk_1e^{-ik_1(\xi-x)}\int_{0}^{T}\!\!\!d\tau\, e^{-ik_1^3(\tau-t)} g(\xi,\tau)\,\delta(\tau-t)\nonumber\\
&= \frac{1}{2\pi}\int_{-\infty}^{\infty}\!\!\!d\xi\int_{-\infty}^{\infty}\!\!dk_1e^{-ik_1(\xi-x)} g(\xi,t)\nonumber\\
&=\int_{-\infty}^{\infty}d\xi \, g(\xi, t)\,\delta(\xi-x)\nonumber\\
&=g(x,t).
\end{align}
In a similar way, it can be shown that $q_y(x,0,t)=h(x,t)$. 
\QED
\end{PROOF}

\section{The KPII equation}\label{kpiispect}
The first Lax equation \eqref{lax1intro} suggests that there exists a solution $\mu$ with the asymptotic behaviour 
\begin{equation}\label{asymptotic}
 \mu=1+O\left(\frac{1}{k}\right), \quad k \rightarrow \infty.
\end{equation}

\subsection{The direct problem}\label{kpiidp}

\begin{prop}\label{kpiidpprop}
Assume that there exists a solution $q(x,y,t), \ (x,y,t)\in \Omega$, for an initial-boundary value problem of equation \eqref{kpiiintro}. Then, the function $\mu$, which is bounded $\forall k \in \mathbb C$ and is defined as
\begin{align}\label{spectralfunctionsii}
\mu(x,y,t,k_R,k_I) = \left\lbrace \begin{array}{ll} \mu_{1}^{+}(x,y,t,k_R,k_I), &\quad k_R\leq0,\ k_I\geq0, \\ \\ \mu_{1}^{-}(x,y,t,k_R,k_I), &\quad k_R\leq0,\ k_I\leq0,  \\ \\ \mu_{2}^{-}(x,y,t,k_R,k_I), &\quad k_R\geq0,\ k_I\leq0, \\ \\ \mu_{2}^{+}(x,y,t,k_R,k_I), &\quad k_R\geq0,\ k_I\geq0, \end{array} \right.
\end{align}
in the different quadrants of the complex $k$-plane (see figure \ref{eigenfunctions}), admits the following representations in terms of $q(x,y,t)$, $q(x,0,t)$ and $q_y(x,0,t)$:
{\small
\begin{align}\label{direct1ii}
&\mu^{\pm}_1=1+\frac{1}{2\pi}\mbox{\large $\Big($}\int_{-\infty}^{0}dl+\int_{-2k_R}^{\infty}dl\mbox{\large $\Big)$}\int_{-\infty}^{\infty}d\xi \int_{0}^{y}d\eta\, e^{-il(\xi-x)+l(l+2k)(\eta-y)}q(\xi,\eta,t)\mu_1^{\pm}\nonumber\\ 
&+\frac{1}{2\pi}\mbox{\large $\Big($}\int_{-\infty}^{0}\!\!\!dl+\!\!\int_{-2k_R}^{\infty}\!\!\!\!\!dl\mbox{\large $\Big)$}\!\!\int_{-\infty}^{\infty}\!\!d\xi\left\lbrace \!\!\!\begin{array}{ll}\int_{0}^{t}d\tau \\ \\ \!-\!\int_{t}^{T}d\tau \end{array}\right.\!\!\!\!\!\!\!\! e^{-il(\xi-x)-l(l+2k)y-4il(l^2+3kl+3k^2)(\tau-t)} H(\xi,\tau,k,l)\phi_1^{\pm}\nonumber\\
&-\frac{1}{2\pi}\int_{0}^{-2k_R}dl\int_{-\infty}^{\infty}d\xi\int_{y}^{\infty}d\eta \, e^{-il(\xi-x)+l(l+2k)(\eta-y)}q(\xi,\eta,t)\mu_1^{\pm},
\end{align}}
{\small
\begin{align}\label{direct2ii}
&\mu^{\pm}_2=1+\frac{1}{2\pi}\mbox{\large $\Big($}\int_{-\infty}^{-2k_R}dl+\int_{0}^{\infty}dl\mbox{\large $\Big)$}\int_{-\infty}^{\infty}d\xi \int_{0}^{y}d\eta\, e^{-il(\xi-x)+l(l+2k)(\eta-y)}q(\xi,\eta,t)\mu_2^{\pm}\nonumber\\ 
&+\frac{1}{2\pi}\mbox{\large $\Big($}\int_{-\infty}^{-2k_R}\!\!\!dl+\!\!\int_{0}^{\infty}\!\!\!\!\!dl\mbox{\large $\Big)$}\!\!\int_{-\infty}^{\infty}\!\!\!\!d\xi\left\lbrace \!\!\!\begin{array}{ll}\int_{0}^{t}d\tau \\ \\ \!-\!\int_{t}^{T}d\tau \end{array}\right.\!\!\!\!\!\!\!\! e^{-il(\xi-x)-l(l+2k)y-4il(l^2+3kl+3k^2)(\tau-t)} H(\xi,\tau,k,l)\phi_2^{\pm}\nonumber\\
&-\frac{1}{2\pi}\int_{-2k_R}^{0}dl\int_{-\infty}^{\infty}d\xi\int_{y}^{\infty}d\eta \, e^{-il(\xi-x)+l(l+2k)(\eta-y)}q(\xi,\eta,t)\mu_2^{\pm},
\end{align}}
where
\begin{equation}\label{phiold}
\phi^\pm_{1,2}(x,t,k_R,k_I)=\mu^\pm_{1,2}(x,0,t,k_R,k_I),
\end{equation}
the function $H$ is given by
\begin{equation}\label{H}
H(x,t,k,l):=3\bigg[g_x(x,t)-2i(l+k)g(x,t)-\partial_{x}^{-1}h(x,t)\bigg]
\end{equation}
and 
\begin{equation*}
 g(x,t)=q(x,0,t), \quad h(x,t)=q_y(x,0,t),
\end{equation*}
as defined by equations \eqref{lkpbc1intro} and \eqref{lkpbc2intro}.
\end{prop}

\begin{PROOF}  
 
Define the Fourier transform of $\mu$ with respect to $x$ by
\begin{subequations}\label{ftpair}
\begin{equation}\label{ft}
 \hat \mu(l,y,t,k_R,k_I):=\int^{\infty}_{-\infty}dx\, e^{-ilx}\Big[\mu(x,y,t,k_R,k_I)-1\Big], \quad l \in \mathbb R, \ k\in \mathbb C,
\end{equation}
with the inverse Fourier transform given by 
\begin{equation}\label{ift}
 \mu(x,y,t,k_R,k_I):=1+\frac{1}{2\pi}\int^{\infty}_{-\infty}dl\, e^{ilx}\hat \mu(l,y,t,k_R,k_I), \quad l \in \mathbb R, \ k\in \mathbb C.
\end{equation}
\end{subequations}
Applying the Fourier transform on the first Lax equation yields an ODE in $y$:
\begin{equation}\label{ftonlax}
\hat \mu_y(l,y,t,k_R,k_I)+l\left(l+2k\right)\hat \mu(l,y,t,k_R,k_I)=\widehat{q\mu}(l,y,t,k_R,k_I),
\end{equation}
where 
\begin{equation*}
\widehat{q\mu}(l,y,t,k_R,k_I):=\int_{-\infty}^{\infty}dx \, e^{-ilx}\left(q\mu\right)(x,y,t,k_R,k_I).
\end{equation*}

We solve this equation by integrating with respect to $y$ either from $y=0$ or from $y=\infty$:
\begin{align}\label{ifii}
\hat{\mu}(l,y,t,k_R,k_I)=\left\{\begin{array}{ll} \!\!\!\!\int_0^y d\eta\, e^{l(l+2k)(\eta-y)}\widehat{q\mu}(l, \eta, t,k_R,k_I)+\hat{\mu}(l,0,t,k_R,k_I)\, e^{-l(l+2k)y} \\ \\ -\int_y^\infty d\eta\, e^{l(l+2k)(\eta-y)}\widehat{q\mu}(l, \eta,t,k_R,k_I), \end{array}\right.
\end{align}
where we have assumed that $\hat \mu\rightarrow 0$ as $y\rightarrow \infty$.

We look for a solution which is bounded $\forall k \in \mathbb C$. In the first expression in equation \eqref{ifii}, $\eta-y\leq 0$, thus we require that $\mathrm{Re}\left[l(l+2k)\right]\geq0$, i.e.
\begin{equation*}
l(l+2k_R)\geq 0\Leftrightarrow \left\lbrace\begin{array}{ll} l\in(-\infty,0]\cup[-2k_R,\infty),\ k_R\leq 0 \\ \\ l\in(-\infty,-2k_R]\cup[0,\infty),\ k_R\geq 0. \end{array} \right. 
\end{equation*}
Analogous considerations are valid for the second expression in equation \eqref{ifii}: \newline $\eta-y\geq 0$ implies that $\mathrm{Re}\left[l(l+2k)\right]\leq0$, i.e.
\begin{equation*}
l(l+2k_R)\leq 0\Leftrightarrow \left\lbrace\begin{array}{ll} l\in[0,-2k_R],\ k_R\leq 0 \\ \\ l\in[-2k_R, 0],\ k_R\geq 0. \end{array} \right. 
\end{equation*}

Applying the inverse Fourier transform \eqref{ift} on equations \eqref{ifii} and writing $\mu$ according to the definition \eqref{spectralfunctionsii}, we obtain
\begin{align}\label{leftii}
&\mu_1(x,y,t,k_R,k_I)=\nonumber\\
&=1+\frac{1}{2\pi}\left(\int_{-\infty}^{0}dl+\int_{-2k_R}^{\infty}dl\right)e^{ilx-l(l+2k)y}\int_{0}^{y}d\eta\, e^{l(l+2k)\eta}\widehat{q\mu}(l,\eta,t,k_R,k_I)\nonumber\\
&+\frac{1}{2\pi}\left(\int_{-\infty}^{0}dl+\int_{-2k_R}^{\infty}dl\right)e^{ilx-l(l+2k)y}\hat\phi_1(l,t,k_R,k_I)\nonumber\\
&-\frac{1}{2\pi}\int_{0}^{-2k_R}dl\, e^{ilx-l(l+2k)y}\int_{y}^{\infty}d\eta \, e^{l(l+2k)\eta}\widehat{q\mu}(l,\eta,t,k_R,k_I), \quad k_R\leq 0,
\end{align}
\begin{align}\label{rightii}
&\mu_2(x,y,t,k_R,k_I)=\nonumber\\
&=1+\frac{1}{2\pi}\left(\int_{-\infty}^{-2k_R}dl+\int_{0}^{\infty}dl\right)e^{ilx-l(l+2k)y}\int_{0}^{y}d\eta\, e^{l(l+2k)\eta}\widehat{q\mu}(l,\eta,t,k_R,k_I)\nonumber\\
&+\frac{1}{2\pi}\left(\int_{-\infty}^{-2k_R}dl+\int_{0}^{\infty}dl\right)e^{ilx-l(l+2k)y}\hat\phi_2(l,t,k_R,k_I)\nonumber\\
&-\frac{1}{2\pi}\int_{-2k_R}^{0}dl\, e^{ilx-l(l+2k)y}\int_{y}^{\infty}d\eta \, e^{l(l+2k)\eta}\widehat{q\mu}(l,\eta,t,k_R,k_I), \quad k_R\geq 0.
\end{align}
Evaluating equation \eqref{leftii} at $y=0$ and recalling the definition \eqref{phiold} we find
\begin{align}
\phi_1(x,t,k_R,k_I)&=1+\, \frac{1}{2\pi}\left(\int_{-\infty}^{0}dl+\int_{-2k_R}^{\infty}dl\right)e^{ilx}\hat \phi_1(l,t,k_R,k_I)\nonumber\\
&-\frac{1}{2\pi}\int_{0}^{-2k_R}dl\int_{0}^{\infty}d\eta \, e^{ilx+l(l+2k)\eta}\widehat{q\mu}(l,\eta,t,k_R,k_I).
\end{align}
The definition \eqref{ift} yields
\begin{align*}
\frac{1}{2\pi}\int_{-\infty}^{\infty}dl\,e^{ilx} \hat \phi_1(l,t,k_R,k_I)&= \frac{1}{2\pi}\left(\int_{-\infty}^{0}dl+\int_{-2k_R}^{\infty}dl\right)e^{ilx}\hat \phi_1(l,t,k_R,k_I)\nonumber\\
&-\frac{1}{2\pi}\int_{0}^{-2k_R}dl\int_{0}^{\infty}d\eta \, e^{ilx+l(l+2k)\eta}\widehat{q\mu}(l,\eta,t,k_R,k_I).
\end{align*}
This equation implies that $\hat \phi_1$ does not satisfy any constraints for $l\in(-\infty,0]\cup[-2k_R, \infty)$, whereas it does satisfy the following constraint for $l\in [0,-2k_R]$: 
\begin{equation}\label{phi1restriction}
 \hat\phi_1(l,t,k_R,k_I)=-\int_{0}^{\infty}d\eta \, e^{l(l+2k)\eta}\widehat{q\mu}(l,\eta,t,k_R,k_I), \quad l\in[0,-2k_R].
\end{equation}

The second Lax equation \eqref{lax2intro} evaluated at $y=0$ yields
\begin{equation}\label{phi}
 \phi_t+4\phi_{xxx}+12ik\phi_{xx}-12k^2\phi_x-\Lambda\phi=0,
\end{equation}
where $\Lambda$ is defined as
\begin{align*}
 \Lambda(x,t,k_R,k_I)&=F(x,0,t,k_R,k_I)\nonumber\\
&=-6g(x,t)\left(\partial_x+ik\right)-3\left[g_x(x,t)+\partial_{x}^{-1}h(x,t)\right].
\end{align*}
By applying the Fourier transform \eqref{ft} on equation \eqref{phi}, we obtain the following ODE in $t$:
\begin{equation*}
 \hat \phi_t(l,t,k_R,k_I)-4il\left(l^2+3kl+3k^2\right) \hat \phi(l,t,k_R,k_I)=\widehat {\Lambda \phi}(l,t,k_R,k_I).
\end{equation*}
Integrating this expression with respect to $t$ either from $t=0$ or from $t=T$, we find:
\begin{align*}
\hat \phi(l,t,k_R,k_I)\!=\!\!\left\lbrace\begin{array}{ll}\int_{0}^{t}d\tau \, e^{-4il(l^2+3kl+3k^2)(\tau-t)}\, \widehat {\Lambda\phi} +  e^{4il(l^2+3kl+3k^2)t}\, \hat \phi\big|_{t=0}\\ \\ \!\!\!\!\!-\!\!\int_{t}^{T}\!\!d\tau e^{-4il(l^2+3kl+3k^2)(\tau-t)}\widehat {\Lambda \phi}+\! e^{-4il(l^2+3kl+3k^2)(T-t)}\hat \phi\big|_{t=T}, k \!\in \!\mathbb C. \end{array}\right.
\end{align*}
Integration by parts of the integral in $\widehat {\Lambda\phi}$ implies
\begin{align}\label{if'ii}
\hat \phi(l,t,k_R,k_I)\!=\!\!\left\lbrace\begin{array}{ll}\int_{0}^{t}d\tau \, e^{-4il(l^2+3kl+3k^2)(\tau-t)} \widehat {H\phi} + e^{4il(l^2+3kl+3k^2)t}\, \hat \phi\big|_{t=0} \\ \\ \!\!\!\!\!\!-\!\!\int_{t}^{T}\!\!d\tau e^{-4il(l^2+3kl+3k^2)(\tau-t)}\widehat {H \phi}+\! e^{-4il(l^2+3kl+3k^2)(T-t)}\hat \phi\big|_{t=T}, k \!\in \!\mathbb C. \end{array}\right.
\end{align}
where the quantity $H$ is defined by equation \eqref{H}, and 
\begin{equation*}
\widehat {H\phi}(l,t,k_R,k_I):= \int_{-\infty}^{\infty}dx\, e^{-ilx} \left(H\phi\right)(x,\tau,k_R,k_I).
\end{equation*}
In order for $\phi$ to be bounded $\forall k \in \mathbb C$, we require that the exponential involved in equation \eqref{if'ii} is bounded. Noting that
\begin{equation*}
\mathrm{Re}\left[-4il(l^2+3kl+3k^2)\right]=12k_Il\left(l+2k_R\right),
\end{equation*}
we are led to consider the following cases (see also figure \ref{eigenfunctions}):
\begin{itemize}
\item If $k\in\{k\in \mathbb{C}: k_R\leq 0,k_I\geq 0\}$, then $l\in(-\infty,0]\cup[-2k_R, \infty)$ implies that $\tau-t\leq 0$, whereas $l\in[0,-2k_R]$ implies that $\tau-t\geq 0$. Thus, recalling our analysis for $\hat \phi_1$, we can choose 
\begin{equation*}
 \hat \phi_1^+(l,0,k_R,k_I)=0, \quad l\in(-\infty,0]\cup [-2k_R, \infty),
\end{equation*}
Equation \eqref{phi1restriction} imposes
\begin{equation*}
 \hat \phi_1^+(l,T,k_R,k_I)=-\int_{0}^{\infty}d\eta \, e^{l(l+2k)\eta}\widehat{q\mu_1^+}(l,\eta,T,k_R,k_I), \quad l\in[0,-2k_R].
\end{equation*}
Therefore, using equation \eqref{if'ii} we find:
\begin{align}\label{phi1hat}
\hat \phi_1^+(l,t,k_R,k_I)\!=\!\left\lbrace\begin{array}{ll}\!\!\int_{0}^{t}d\tau \, e^{-4il(l^2+3kl+3k^2)(\tau-t)} \widehat {H\phi_1^+},\ l\in\!(-\infty,0)\cup (-2k_R, \infty) \\ \\\!\!\!\!-\int_{t}^{T}d\tau \, e^{-4il(l^2+3kl+3k^2)(\tau-t)} \widehat {H\phi_1^+} \vspace{2mm}\\ \!\!\!\!-e^{-4il(l^2+3kl+3k^2)(T-t)}\!\int_{0}^{\infty}\!\!d\eta \, e^{l(l+2k)\eta}\widehat{q\mu_1^+}\Big|_{t=T},\ l\in\![0,-2k_R]. \end{array}\right.
\end{align}
Similar calculations for $k$ in each of the other three quandrants yield the following results:
\item for $k\in\{k\in \mathbb{C}: k_R\leq 0,k_I\leq 0\}$, 
\begin{align}\label{phi1-hat}
\hat \phi_1^-(l,t,k_R,k_I)=\left\lbrace\begin{array}{ll}\!\!\!\!\!-\!\int_{t}^{T}\!\!d\tau \, e^{-4il(l^2+3kl+3k^2)(\tau-t)} \widehat {H\phi_1^-},\ l\in\!(-\infty,0]\!\cup\![-2k_R, \infty)\\  \\ \int_{0}^{t}d\tau \, e^{-4il(l^2+3kl+3k^2)(\tau-t)} \widehat {H\phi_1^-}\vspace{2mm}\\
-e^{4il(l^2+3kl+3k^2)t}\int_{0}^{\infty}d \eta \, e^{l(l+2k)\eta}\widehat{q\mu_1^-}\Big|_{t=0},\ l\in[0,-2k_R].\end{array}\right.
\end{align}
\item for $k\in\{k\in \mathbb{C}: k_R\geq 0,k_I\geq 0\}$, 
\begin{align}\label{phi2+hat}
\hat \phi_2^+(l,t,k_R,k_I)\!=\!\left\lbrace\begin{array}{ll}\!\!\!\int_{0}^{t}d\tau \, e^{-4il(l^2+3kl+3k^2)(\tau-t)} \widehat {H\phi_2^+},\ l\in(-\infty,-2k_R]\cup[0, \infty)\\  \\ \!\!\!\!-\int_{t}^{T}d\tau \, e^{-4il(l^2+3kl+3k^2)(\tau-t)} \widehat {H\phi_2^+}\vspace{2mm}\\ \!\!\!\!\!-e^{-4il(l^2+3kl+3k^2)(T-t)}\!\int_{0}^{\infty}\!\!d \eta \, e^{l(l+2k)\eta}\widehat{q\mu_2^+}\Big|_{t=T},\ l\in[-2k_R,0].\end{array}\right.
\end{align} 
\item and for $k\in\{k\in \mathbb{C}: k_R\geq 0,k_I\leq 0\}$, 
\begin{align}\label{phi2-hat}
\hat \phi_2^-(l,t,k_R,k_I)=\left\lbrace\begin{array}{ll}\!\!\!\!\!-\!\int_{t}^{T}\!\!d\tau \, e^{-4il(l^2+3kl+3k^2)(\tau-t)} \widehat {H\phi_2^-},\ l\in\!(-\infty,-2k_R]\!\cup\![0, \infty)\\ \\ \int_{0}^{t}d\tau \, e^{-4il(l^2+3kl+3k^2)(\tau-t)} \widehat {H\phi_2^-}\vspace{2mm}\\ 
-e^{4il(l^2+3kl+3k^2)t}\int_{0}^{\infty}d \eta \, e^{l(l+2k)\eta}\widehat{q\mu_2^-}\Big|_{t=0},\ l\in[-2k_R, 0].\end{array}\right.
\end{align}
\end{itemize}
Inserting equations \eqref{phi1hat}-\eqref{phi2-hat} into equations \eqref{leftii} and \eqref{rightii}, we obtain for $\mu_1$ and $\mu_2$ the representations \eqref{direct1ii} and \eqref{direct2ii} respectively.
\QED
\end{PROOF}

\subsection{The global relation}\label{kpiigr}

\begin{prop}\label{GRpropkpii}
Define $q_0$, $g$ and $h$ by equations \eqref{icintro}, \eqref{lkpbc1intro} and \eqref{lkpbc2intro} respectively. Also, define the functions $\phi$ and $H$ by equations \eqref{phiold} and \eqref{H}.

Then, these functions satisfy the so-called global relation:
{\small\begin{align}
&\int_{-\infty}^{\infty}d\xi\int_{0}^\infty d\eta\, e^{-il\xi+l(l+2k)\eta} q_0(\xi,\eta)\mu_0-\int_{-\infty}^{\infty}d\xi\int_{0}^{t}d\tau\, e^{-il\xi-4il(l^2+3kl+3k^2)\tau}H(\xi,\tau,k,l)\phi\nonumber\\
&=\int_{-\infty}^{\infty}d\xi\int_{0}^\infty d\eta\, e^{-il\xi+l(l+2k)\eta-4il(l^2+3kl+3k^2)t} q(\xi,\eta,t)\mu,\quad k \in \mathbb C,\quad l(l+2k_R)\leq0. \label{GRjump}
\end{align}}

\end{prop}

\begin{PROOF}

Equation \eqref{ftonlax} yields
\begin{equation*}
\left(\hat \mu\, e^{l(l+2k)y-4il(l^2+3kl+3k^2)t}\right)_y=e^{l(l+2k)y-4il(l^2+3kl+3k^2)t}\int_{-\infty}^{\infty}dx\, e^{-ilx}q\mu.
\end{equation*}
Also, the second equation of the Lax pair implies
\begin{equation*}
\left(\hat \mu\, e^{l(l+2k)y-4il(l^2+3kl+3k^2)t}\right)_t=e^{l(l+2k)y-4il(l^2+3kl+3k^2)t}\int_{-\infty}^{\infty}dx\, e^{-ilx}H\mu.
\end{equation*}
The above equations imply
\begin{align*}
&\left(\!e^{l(l+2k)y-4il(l^2+3kl+3k^2)t}\!\!\int_{-\infty}^{\infty}\!\!\!\!\!\!dx\, e^{-ilx}q\mu\!\right)_{t}\!\!\!=\left(\!e^{l(l+2k)y-4il(l^2+3kl+3k^2)t}\!\!\int_{-\infty}^{\infty}\!\!\!\!\!\!dx\, e^{-ilx}H\mu\!\right)_y\!\!.
\end{align*}
Using Green's theorem over the domain depicted in figure \ref{lA1}, we obtain the global relation \eqref{GRjump}:
{\small\begin{align*}
&\int_{-\infty}^{\infty}d\xi\int_{0}^\infty d\eta\, e^{-il\xi+l(l+2k)\eta} q_0(\xi,\eta)\mu_0-\int_{-\infty}^{\infty}d\xi\int_{0}^{t}d\tau\, e^{-il\xi-4il(l^2+3kl+3k^2)\tau}H(\xi,\tau,k,l)\phi\nonumber\\
&=\int_{-\infty}^{\infty}d\xi\int_{0}^\infty d\eta\, e^{-il\xi+l(l+2k)\eta-4il(l^2+3kl+3k^2)t} q(\xi,\eta,t)\mu,\quad k \in \mathbb C,\ \ l(l+2k_R)\leq0. 
\end{align*}}
Note that the constraint on $l$ is needed in order for the exponential $e^{l(l+2k)\eta}$ to be bounded.
\QED
\end{PROOF}

\subsection{The inverse problem}\label{invprob}

It will be shown later (see equation \eqref{imagcont}) that 
\begin{equation*}
\left(\mu_1^\pm-\mu_2^\pm\right)\big|_{k_R=0}=0.
\end{equation*} 
Assuming that $\mu$ is continuous across the imaginary $k$-axis, Pompeiu's formula \eqref{pompeiuintro}
\begin{equation*}
\mu(k_R,k_I)=\frac{1}{2i\pi}\int_{\partial \mathcal D}\frac{d\nu}{\nu-k}\, \mu(\nu_R,\nu_I)+\frac{1}{2i\pi}\iint_{\mathcal D}\frac{d\nu\wedge d\bar\nu}{\nu-k}\, \frac{\partial \mu}{\partial \bar \nu}(\nu_R, \nu_I),
\end{equation*}
applied in the simple domain $\mathcal D$ (see figure \ref{pompeiu}), yields
\begin{align}\label{Pompeiu}
\mu&=1+ \frac{1}{2i\pi} \int_{-\infty}^0 \frac{d\nu_R}{\nu_R-k}\, \Delta \mu_1 +\frac{1}{2i\pi}\int_{0}^{\infty} \frac{d\nu_R }{\nu_R-k}\, \Delta \mu_2\nonumber\\
&-\frac{1}{\pi}\int_{0}^{\infty}d\nu_R\int_{0}^{\infty}\frac{d\nu_I}{\nu-k}\, \frac{\partial \mu_2^+}{\partial \bar \nu}-\frac{1}{\pi}\int_{-\infty}^{0}d\nu_R\int_{0}^{\infty}\frac{d\nu_I}{\nu-k}\, \frac{\partial \mu_1^+}{\partial \bar \nu}\nonumber\\
&-\frac{1}{\pi}\int_{0}^{\infty}d\nu_R\int_{-\infty}^{0} \frac{d\nu_I}{\nu-k}\,\frac{\partial \mu_2^-}{\partial \bar \nu}-\frac{1}{\pi}\int_{-\infty}^{0}d\nu_R\int_{-\infty}^{0} \frac{d\nu_I}{\nu-k}\,\frac{\partial \mu_1^-}{\partial \bar \nu},
\end{align}
where 
\begin{equation}\label{Deltaj}
\Delta \mu_j(x,y,t,k_R):=\left(\mu_j^+-\mu_j^-\right)(x,y,t,k_R,0),\quad j=1,2.
\end{equation}
Therefore, we need to derive expressions for the discontinuities of $\mu$ across the real $k$-axis and for the d-bar derivatives of $\mu$ in the complex $k-$plane.

\begin{prop}(The d-bar derivatives)\label{kpiid-barprop}
Define the functions $\mu_1^\pm$ and $\mu_2^\pm$ by equations \eqref{spectralfunctionsii}, \eqref{direct1ii} and \eqref{direct2ii}. Then,
\begin{equation}\label{kbar1finalprop}
\frac{\partial \mu_1^\pm}{\partial \bar k }(x,y,t,k_R,k_I)= e^{-2ik_Rx+4ik_Rk_Iy-8ik_R(k_R^2-3k_I^2)t}\, \gamma_1^{\pm}(k_R,k_I)\mu_2^\pm(x,y,t,-k_R, k_I)
\end{equation}
\textit{and}
\begin{equation}\label{kbar2finalprop}
\frac{\partial \mu_2^\pm}{\partial \bar k }(x,y,t,k_R,k_I)\!=\!-e^{-2ik_Rx+4ik_Rk_Iy-8ik_R(k_R^2-3k_I^2)t} \gamma_2^{\pm}(k_R,k_I)\mu_1^\pm(x,y,t,-k_R, k_I), 
\end{equation}
where 
\begin{align}
\gamma_1^+
(k_R,k_I)&=\beta_1^+(k_R,k_I), & k_R\leq0,\, k_I\geq0,\label{gamma1pprop}\\
\nonumber\\
\gamma_1^-
(k_R,k_I)&=\beta_1^-(k_R,k_I)-\alpha^-_1(k_R,k_I), & k_R\leq0,\, k_I\leq0,\label{gamma1mprop}\\
\nonumber\\
\gamma_2^+
(k_R,k_I)&=\beta_2^+(k_R,k_I),& k_R\geq0,\, k_I\geq0,\label{gamma2pprop}\\
\nonumber\\
\gamma_2^-
(k_R,k_I)&=\beta_2^-(k_R,k_I)-\alpha^-_2(k_R,k_I),& k_R\geq0,\, k_I\leq0\label{gamma2mprop}
\end{align}
and $\alpha_{1,2}^\pm$, $\beta_{1,2}^\pm$ are defined in the appropriate quadrant of the complex $k$-plane by the following formulae:
\begin{align}
\alpha_1^-(k_R,k_I)&=\frac{1}{2\pi}\int_{-\infty}^{\infty}\!d\xi\int_0^T\!\!\!d\tau\, e^{2ik_R\xi+8ik_R(k_R^2-3k_I^2)\tau} H(\xi,\tau,k,-2k_R)\phi_1^-(\xi,\tau,k_R,k_I)\label{alpha1-prop},\\
\beta_1^{\pm}(k_R,k_I)&=\frac{1}{2\pi} \int_{-\infty}^{\infty}d\xi\int_0^\infty d\eta\, e^{2ik_R\xi-4ik_Rk_I\eta}q_0(\xi,\eta)\mu_{1}^\pm(\xi, \eta,0,k_R,k_I)\label{beta1prop},\\
\alpha_2^-(k_R,k_I)&=\frac{1}{2\pi}\int_{-\infty}^{\infty}\!d\xi\int_0^T\!\!\!d\tau\, e^{2ik_R\xi+8ik_R(k_R^2-3k_I^2)\tau}H(\xi,\tau,k,-2k_R)\phi_2^-(\xi,\tau,k_R,k_I)\label{alpha2-prop},\\
\beta_2^{\pm}(k_R,k_I)&=\frac{1}{2\pi} \int_{-\infty}^{\infty}d\xi\int_0^\infty d\eta\ e^{2ik_R\xi-4ik_Rk_I\eta}q_0(\xi,\eta)\mu_{2}^\pm(\xi, \eta,0,k_R,k_I)\label{beta2prop},
\end{align}
where the functions $\phi_{1,2}^-$ and $H$ are defined by equations \eqref{phiold} and \eqref{H} respectively.
\end{prop}

\begin{PROOF} 
 
Differentiating equation \eqref{direct1ii} according to the definition \eqref{d-bardef} of the d-bar derivative, we find the following expression:
{\small\begin{align*}
&\frac{\partial \mu_1^\pm}{\partial \bar k }=\frac{1}{2\pi}\mbox{\large $\Big($}\int_{-\infty}^0 dl + \int_{-2k_R}^{\infty} dl\mbox{\large $\Big)$}\int_{-\infty}^{\infty}d\xi\int_0^y d\eta\, e^{-il(\xi-x)+l(l+2k)(\eta-y)} q \frac{\partial \mu_1^\pm}{\partial \bar k }\nonumber\\
&+\frac{1}{2\pi}\mbox{\large $\Big($}\!\int_{-\infty}^0\!\!\!\!\!dl +\!\! \int_{-2k_R}^{\infty}\!\!\!\!\!dl\mbox{\large $\Big)$}\!\!\!\int_{-\infty}^{\infty}\!\!\!\!d\xi \left\{\!\!\!\!\begin{array}{ll}\int_0^t d\tau \\ \\ -\int_t^T d\tau\end{array}\!\!\!\!\!\!e^{-il(\xi-x)-l(l+2k)y-4il(l^2+3kl+3k^2)(\tau-t)}\!H(\xi,\tau,k,l)\frac{\partial \phi_1^\pm}{\partial \bar k}\right.\nonumber\\
&-\frac{1}{2\pi}\int_{0}^{-2k_R} dl \int_{-\infty}^{\infty}d\xi\int_y^\infty d\eta \, e^{-il(\xi-x)+l(l+2k)(\eta-y)}q \frac{\partial \mu_1^\pm}{\partial \bar k }\nonumber\\
&+ \frac{1}{2\pi}\int_{-\infty}^{\infty}d\xi \int_0^y d\eta \, e^{2ik_R(\xi-x)-4ik_Rk_I(\eta-y)} q\mu_1^\pm  \nonumber\\
&+\frac{1}{2\pi} \int_{-\infty}^{\infty}d\xi\left\{\begin{array}{ll}\int_0^t d\tau \\ \\ -\int_t^T d\tau\end{array}e^{2ik_R(\xi-x)+4ik_Rk_Iy+8ik_R(k_R^2-3k_I^2)(\tau-t)} H(\xi,\tau,k,-2k_R)\frac{\partial \phi_1^\pm}{\partial \bar k}\right.\nonumber\\
&+\frac{1}{2\pi} \int_{-\infty}^{\infty}d\xi\int_y^\infty d\eta  \, e^{2ik_R(\xi-x)-4ik_Rk_I(\eta-y)} q\mu_1^\pm, \quad k_R\leq 0.
\end{align*}}
This expression can be written as an integral equation for $\partial \mu_1^\pm/\partial \bar k$:
{\small\begin{align}
&\frac{\partial \mu_1^\pm}{\partial \bar k }= \frac{1}{2\pi}\mbox{\large $\Big($}\int_{-\infty}^0 dl + \int_{-2k_R}^{\infty} dl\mbox{\large $\Big)$} \int_{-\infty}^{\infty}d\xi\int_0^y d\eta\, e^{-il(\xi-x)+l(l+2k)(\eta-y)} q \frac{\partial \mu_1^\pm}{\partial \bar k }\nonumber\\
&+\frac{1}{2\pi}\mbox{\large $\Big($}\!\int_{-\infty}^0\!\!\!\!\!dl +\!\! \int_{-2k_R}^{\infty}\!\!\!\!\!dl\mbox{\large $\Big)$}\!\!\!\int_{-\infty}^{\infty}\!\!\!\!d\xi \left\{\!\!\!\!\begin{array}{ll}\int_0^t d\tau \\ \\ -\int_t^T d\tau\end{array}\!\!\!\!\!\!e^{-il(\xi-x)-l(l+2k)y-4il(l^2+3kl+3k^2)(\tau-t)}\!H(\xi,\tau,k,l)\frac{\partial \phi_1^\pm}{\partial \bar k}\right.\nonumber\\
&-\frac{1}{2\pi}\int_{0}^{-2k_R} dl \int_{-\infty}^{\infty}d\xi\int_y^\infty d\eta \, e^{-il(\xi-x)+l(l+2k)(\eta-y)}q \frac{\partial \mu_1^\pm}{\partial \bar k }\nonumber\\
&+e^{-2ik_Rx+4ik_Rk_Iy}\, \Gamma^\pm_1(t,k_R, k_I),\label{kbar1}
\end{align}}
where $\Gamma^\pm_1$ is defined by
\begin{align}
\Gamma^\pm_1(t,k_R,k_I)&= \frac{1}{2\pi}\int_{-\infty}^{\infty}d\xi\int_0^\infty d\eta\, e^{2ik_R\xi-4ik_Rk_I\eta} q\mu_1^\pm \nonumber\\
&+\frac{1}{2\pi}\int_{-\infty}^{\infty}d\xi\left\{\begin{array}{ll}\int_0^t d\tau \\ \\ -\int_t^T d\tau\end{array}e^{2ik_R\xi+8ik_R(k_R^2-3k_I^2)(\tau-t)}  H(\xi,\tau,k,-2k_R) \phi_1^\pm\right..\label{gammadeftemp}
\end{align}

Using the identity \eqref{tT}, we can write equation \eqref{gammadeftemp} in the following form:
\begin{align}
&\Gamma^\pm_1(t,k_R,k_I)=\, \frac{1}{2\pi} \int_{-\infty}^{\infty}d\xi\int_0^\infty d\eta\, e^{2ik_R\xi-4ik_Rk_I\eta} q\mu_1^\pm \nonumber\\
&+\frac{1}{2\pi} \int_{-\infty}^{\infty}d\xi \int_0^t d\tau \, e^{2ik_R\xi+8ik_R(k_R^2-3k_I^2)(\tau-t)} H(\xi,\tau,k,-2k_R)\phi_1^\pm\nonumber\\
&-\frac{1}{2\pi}\int_{-\infty}^{\infty}d\xi\int_0^T d\tau \, e^{2ik_R\xi+8ik_R(k_R^2-3k_I^2)(\tau-t)} H(\xi,\tau,k,-2k_R)\phi_1^\pm\ \mathcal H(-k_I),\label{gammadef}
\end{align}
where $\mathcal H$ denotes the Heaviside function.

The global relation \eqref{GRjump} for $l=-2k_R$ yields
\begin{align*}
&\int_{-\infty}^{\infty}d\xi\int_{0}^{\infty}d\eta\ e^{2ik_R\xi-4ik_Rk_I\eta}q\mu_1^\pm=\nonumber\\
&=\int_{-\infty}^{\infty}d\xi\int_{0}^{\infty}d\eta\ e^{2ik_R\xi-4ik_Rk_I\eta-8ik_R(k_R^2-3k_I^2)t}q_0\mu_{1_0}^\pm\nonumber\\
&-\int_{-\infty}^{\infty}d\xi\int_0^td\tau\ e^{2ik_R\xi+8ik_R(k_R^2-3k_I^2)(\tau-t)}H(\xi,\tau,k,-2k_R)\phi_1^\pm.
\end{align*}
Hence, equation \eqref{gammadef} implies that
\begin{equation*}
\Gamma_1^\pm(t,k_R,k_I)=e^{-8ik_R(k_R^2-3k_I^2)t}\,\gamma_1^\pm(k_R,k_I),
\end{equation*}
where $\gamma_1^\pm$ are defined in proposition \ref{kpiidpprop}. Then, equation \eqref{kbar1} becomes
{\small\begin{align}
&\frac{\partial \mu_1^\pm}{\partial \bar k }= \frac{1}{2\pi}\mbox{\large $\Big($}\int_{-\infty}^0 dl + \int_{-2k_R}^{\infty} dl\mbox{\large $\Big)$}\int_{-\infty}^{\infty}d\xi\int_0^y d\eta\, e^{-il(\xi-x)+l(l+2k)(\eta-y)} q \frac{\partial \mu_1^\pm}{\partial \bar k }\nonumber\\
&+\frac{1}{2\pi}\mbox{\large $\Big($}\!\int_{-\infty}^0\!\!\!\!\!dl +\!\! \int_{-2k_R}^{\infty}\!\!\!\!\!dl\mbox{\large $\Big)$}\!\!\!\int_{-\infty}^{\infty}\!\!\!\!d\xi \left\{\!\!\!\!\begin{array}{ll}\int_0^t d\tau \\ \\ -\int_t^T d\tau\end{array}\!\!\!\!\!\!e^{-il(\xi-x)-l(l+2k)y-4il(l^2+3kl+3k^2)(\tau-t)}\!H(\xi,\tau,k,l)\frac{\partial \phi_1^\pm}{\partial \bar k}\right.\nonumber\\
&-\frac{1}{2\pi}\int_{0}^{-2k_R} dl \int_{-\infty}^{\infty}d\xi\int_y^\infty d\eta\, e^{-il(\xi-x)+l(l+2k)(\eta-y)}q \frac{\partial \mu_1^\pm}{\partial \bar k }\nonumber\\
&+e^{-2ik_Rx+4ik_Rk_Iy-8ik_R(k_R^2-3k_I^2)t}\, \gamma^\pm_1(k_R, k_I).\label{kbar1b}
\end{align}}
Let us introduce the notations
\begin{equation}
\mathcal M_1^\pm(x,y,t,k_R,k_I)=\frac{e^{2ik_Rx-4ik_Rk_Iy+8ik_R(k_R^2-3k_I^2)t}}{\gamma^\pm_1(k_R, k_I)}\ \frac{\partial \mu_1^\pm}{\partial \bar k }(x,y,t,k_R,k_I)\label{M1pm}
\end{equation}
and 
\begin{equation*}
\mathcal F_1^\pm(x,t,k_R,k_I)=\mathcal M_1^\pm(x,0,t,k_R,k_I). 
\end{equation*}
Then, equation \eqref{kbar1b} can be written in the form
{\small\begin{align}
&\mathcal M_1^\pm= 1+\frac{1}{2\pi}\mbox{\large $\Big($}\int_{-\infty}^0\!\!\!\!dl + \int_{-2k_R}^{\infty}\!\!\!\!dl\mbox{\large $\Big)$}\mbox{\large $\Big[$}\!\int_{-\infty}^{\infty}\!\!d\xi\int_0^y \!\!d\eta \, e^{-i(l+2k_R)(\xi-x)+(l(l+2k)+4ik_Rk_I)(\eta-y)} q \mathcal M_1^\pm\nonumber\\
&+\!\int_{-\infty}^{\infty}\!\!\!\!d\xi\left\{\!\!\!\begin{array}{ll}\int_0^t d\tau \\ \\ -\int_t^T d\tau\end{array}\!\!\!\!\!\!\!\!\! e^{-i(l+2k_R)(\xi-x)-(l(l+2k)+4ik_Rk_I)y+8ik_R(k_R^2-3k_I^2)\tau-4il(l^2+3kl+3k^2)(\tau-t)} H\mathcal F_1^\pm\mbox{\large $\Big]$}\right.\nonumber\\
&-\frac{1}{2\pi}\int_{0}^{-2k_R} dl \int_{-\infty}^{\infty}d\xi\int_y^\infty d\eta\, e^{-i(l+2k_R)(\xi-x)+(l(l+2k)+4ik_Rk_I)(\eta-y)} q \mathcal M_1^\pm.\label{M1}
\end{align}}

Let $l^\prime=l+2k_R$ in the integrals on the RHS of equation \eqref{M1}. Consequently,
\begin{align*}
l(l+2k)&= l^\prime(l^\prime-2\bar k)-4ik_Rk_I\\
il&= il^\prime-2ik_R\\
-4il(l^2+3kl+3k^2)&=-4il^\prime({l^\prime}^2-3\bar k l^\prime +3\bar k^2)+8ik_R(k_R^2-3k_I^2).
\end{align*}
Dropping the prime, equation \eqref{M1} becomes
{\small\begin{align}\label{M1b}
&\mathcal M_1^\pm=1+\frac{1}{2\pi}\mbox{\large $\Big($}\int_{-\infty}^{2k_R}\!\!\! dl +\!\! \int_{0}^{\infty}\!\!\! dl\mbox{\large $\Big)$}\mbox{\large $\Big[$}\int_{-\infty}^{\infty}d\xi\int_0^y d\eta\, e^{-il(\xi-x)+l(l-2\bar k)(\eta-y)} q\mathcal M_1^\pm\nonumber\\
&+\int_{-\infty}^{\infty}d\xi\left\{\begin{array}{ll}\int_0^t d\tau \\ \\ -\int_t^T d\tau\end{array}e^{-il(\xi-x)-l(l-2\bar k)y-4il(l^2-3\bar kl+3\bar k^2)(\tau-t)} H(\xi,\tau,k,l-2k_R)\mathcal F_1^\pm\right.\mbox{\large $\Big]$}\nonumber\\
&-\frac{1}{2\pi}\int_{2k_R}^{0} dl \int_{-\infty}^{\infty}d\xi\int_y^\infty d\eta\, e^{-il(\xi-x)+l(l-2\bar k)(\eta-y)} q\mathcal M_1^\pm.
\end{align}}
Employing the transformation $k\mapsto -\bar k$ in the definition \eqref{H} of $H$, we find
\begin{equation*}
H(\xi,\tau,k,l-2k_R)\mapsto H(\xi,\tau,-\bar k,l+2k_R)=H(\xi,\tau,k,l).
\end{equation*}
Introducing the notation
\begin{equation*}
\check{\mathcal M}_1^\pm=\mathcal M_1^\pm (x,y,t,-k_R,k_I),
\end{equation*}
equation \eqref{M1b} becomes:
{\small\begin{align*}
&\check{\mathcal M}_1^\pm=1+\frac{1}{2\pi}\mbox{\large $\Big($}\int_{-\infty}^{-2k_R}\!\!\! dl +\!\! \int_{0}^{\infty}\!\!\! dl\mbox{\large $\Big)$}\int_{-\infty}^{\infty}d\xi\int_0^y d\eta\, e^{-il(\xi-x)+l(l+2 k)(\eta-y)} q\check{\mathcal M}_1^\pm\nonumber\\
&+\frac{1}{2\pi}\mbox{\large $\Big($}\!\int_{-\infty}^{-2k_R}\!\!\!\!\!dl +\!\! \int_{0}^{\infty} \!\!\!\!dl\mbox{\large $\Big)$}\!\!\int_{-\infty}^{\infty}\!\!\!\!d\xi\left\{\!\!\!\!\begin{array}{ll}\int_0^t d\tau \\ \\ -\int_t^T d\tau\end{array}\!\!\!\!\!\!\!e^{-il(\xi-x)-l(l+2 k)y-4il(l^2+3kl+3k^2)(\tau-t)} H(\xi,\tau,k,l)\check{\mathcal F}_1^\pm\right.\nonumber\\
&-\frac{1}{2\pi}\int_{-2k_R}^{0} dl \int_{-\infty}^{\infty}d\xi\int_y^\infty d\eta\, e^{-il(\xi-x)+l(l+2 k)(\eta-y)} q\check{\mathcal M}_1^\pm.
\end{align*}}
Thus, it follows that this equation is identical with equation \eqref{direct2ii} hence, assuming that the solution is unique, it follows that
\begin{equation*}
\mathcal M_1^\pm(x,y,t,-k_R, k_I)=\mu_2^\pm(x,y,t,k_R,k_I),
\end{equation*}
or, equivalently,
\begin{equation*}
\mathcal M_1^\pm(x,y,t, k_R,k_I)=\mu_2^\pm(x,y,t,-k_R, k_I).
\end{equation*}

Therefore, using the definition \eqref{M1pm}, we obtain equation \eqref{kbar1finalprop} for the d-bar derivative of $\mu_1^\pm$. In a similar way, we can derive the expression \eqref{kbar2finalprop} for the d-bar derivative of $\mu_2^\pm$.
\QED
\end{PROOF}
 
In order to show that $\mu$ is continuous across the imaginary $k$-axis, we evaluate the expressions \eqref{direct1ii} and \eqref{direct2ii} at $k_R=0$ and then subtract the resulting equations:
{\small\begin{align*}
&\mbox{\large $\big($}\mu_1^\pm-\mu_2^\pm\mbox{\large $\big)$}\mbox{\large $\Big|$}_{k_R=0}=\frac{1}{2\pi}\int_{-\infty}^{\infty}dl\int_{-\infty}^{\infty}d\xi\int_0^yd\eta\ e^{-il(\xi-x)+l(l+2ik_I)(\eta-y)}q\mbox{\large $\big($}\mu_1^\pm-\mu_2^\pm\mbox{\large $\big)$}\mbox{\large $\Big|$}_{k_R=0}\nonumber\\
&+\frac{1}{2\pi}\int_{-\infty}^{\infty}\!\!dl\int_{-\infty}^{\infty}\!\!d\xi\int_0^t\!\!d\tau\, e^{-il(\xi-x)-l(l+2ik_I)y-4il(l^2+3ik_Il-3k_I^2)(\tau-t)}q \mbox{\large $\Big[$}H\mbox{\large $\big($}\phi_1^\pm-\phi_2^\pm\mbox{\large $\big)$}\mbox{\large $\Big]$}\mbox{\large $\Big|$}_{k_R=0}.
\end{align*}}
This equation evaluated at $y=0$ implies
{\small\begin{align*}
&\mbox{\large $\big($}\phi_1^\pm-\phi_2^\pm\mbox{\large $\big)$}\mbox{\large $\Big|$}_{k_R=0}=\nonumber\\
&=\frac{1}{2\pi}\int_{-\infty}^{\infty}\!\!dl\int_{-\infty}^{\infty}\!\!d\xi\int_0^t\!\!d\tau\, e^{-il(\xi-x)-l(l+2ik_I)y-4il(l^2+3ik_Il-3k_I^2)(\tau-t)}q \mbox{\large $\Big[$}H\mbox{\large $\big($}\phi_1^\pm-\phi_2^\pm\mbox{\large $\big)$}\mbox{\large $\Big]$}\mbox{\large $\Big|$}_{k_R=0};
\end{align*}}
hence, assuming uniqueness we find $\left(\phi_1^\pm-\phi_2^\pm\right)\!\Big|_{k_R=0}=0$. 

Then, it follows that $\mu_1^\pm-\mu_2^\pm$ satisfies a homogeneous integral equation, hence under the assumption of uniqueness, 
\begin{equation}\label{imagcont}
\mu_1^\pm(x,y,t,0,k_I)=\mu_2^\pm(x,y,t,0,k_I).
\end{equation}
Next, we compute the jump of $\mu$ across the real $k$-axis.

\begin{prop}(The jumps across the real $k$-axis)\label{kpiidiscprop}
Define the functions $\mu_1^\pm$ and $\mu_2^\pm$, by equation \eqref{spectralfunctionsii}, \eqref{direct1ii} and \eqref{direct2ii}. These functions satisfy the following discontinuity relations across the real $k$-axis:
\begin{align}
\Delta\mu_1(x,y,t,k_R)&=\int_{-\infty}^{0}d\lambda\chi_2(k_R,\lambda)\left(e_{\lambda}\breve \mu_2^+\right)(x,y,t,k_R,\lambda)\nonumber\\
&+\int_{-2k_R}^{\infty}d\lambda\chi_1(k_R,\lambda)\left(e_{\lambda}\breve \mu_1^+\right)(x,y,t,k_R,\lambda), \quad k_R\leq0,\label{discontinuity1prop}\\
\Delta\mu_2(x,y,t,k_R)&=\int_{-\infty}^{-2k_R}d\lambda\psi_2(k_R,\lambda)\left(e_{\lambda}\breve \mu_2^+\right)(x,y,t,k_R,\lambda)\nonumber\\
&+\int_{0}^{\infty}d\lambda\psi_1(k_R,\lambda)\left(e_{\lambda}\breve \mu_1^+\right)(x,y,t,k_R,\lambda),\quad k_R\geq0,\label{discontinuity2prop}
\end{align}
where $\Delta\mu_{1,2}$ are defined by equation \eqref{Deltaj}, and
\begin{align}
\breve \mu_j(x,y,t,k_R,\lambda)&:=\mu_j(x,y,t,-k_R-\dfrac{\lambda}{2},\dfrac{i\lambda}{2}),\quad j=1,2,\quad \lambda \in \mathbb R,\label{mubreve}\\
e_{\lambda}(x,y,t,k_R,\lambda)&:=e^{-i(\lambda+2k_R)x-\lambda(\lambda+2k_R)y-4i\lambda(\lambda^2+3k_R\lambda+3k_R^2)t-8ik_R^3t}.\label{elambda}
\end{align}
The functions $\chi_{2}$, $\chi_{1}$, $\psi_{2}$ and $\psi_1$ are the solutions of the following systems of Volterra integral equations:
\begin{subequations}
\begin{align}\label{chi2eqfinal}
&\chi_2(k_R,\lambda)-\int_{-\infty}^{\lambda}dl\chi_2(k_R,l) r_2(k_R, l, -2k_R-\lambda)\nonumber\\
&-\int_{-2k_R-\lambda}^\infty dl\chi_1(k_R, l)r_1(k_R,\lambda,-2k_R-\lambda)=p_1(k_R,-2k_R-\lambda), \quad k_R\leq0,\ \lambda\leq0,\\
&\chi_1(k_R,\lambda)-\int_{-\infty}^{-2k_R-\lambda}dl \chi_2(k_R,l)r_2(k_R, l, -2k_R-\lambda)\nonumber\\
&-\int_{\lambda}^{\infty}dl \chi_1(k_R, l)r_1(k_R,l, -2k_R-\lambda)=p_1(k_R,-2k_R-\lambda),\quad k_R\leq0,\ \lambda\geq-2k_R \label{chi1eqfinal}
\end{align}
\end{subequations}
and
\begin{subequations}
\begin{align}\label{psi2eqfinal}
&\psi_2(k_R,\lambda)+\int_{-\infty}^{\lambda}dl\psi_2(k_R,l) r_2(k_R, l, -2k_R-\lambda)\nonumber\\
&+\!\int_0^{-2k_R-\lambda} \!\!\!dl\psi_1(k_R, l)r_1(k_R,\lambda,-2k_R-\lambda)\!=p_1(k_R,-2k_R-\lambda), \ k_R\geq0,\, \lambda\leq-2k_R,\\
&\psi_1(k_R,\lambda)+\int_{-\infty}^{-2k_R-\lambda}dl \psi_2(k_R,l)r_2(k_R, l, -2k_R-\lambda)\nonumber\\
&+\int_{0}^{\lambda}dl \psi_1(k_R, l)r_1(k_R,l, -2k_R-\lambda)=p_1(k_R,-2k_R-\lambda),\quad k_R\geq0,\ \lambda\geq0,\label{psi1eqfinal}
\end{align}
\end{subequations}
with $r_2$, $r_1$ and $p_1$ defined by
\begin{align}
&r_2(k_R,\lambda, l)=\frac{1}{2\pi} \int_{-\infty}^{\infty}d\xi\int_{0}^\infty\!\! d\eta\, e^{-il\xi+l(l+2k_R)\eta} q_0(\xi,\eta)\left(e_{\lambda}\breve \mu_2^+\right) (\xi,\eta,0,k_R,\lambda),\label{r2}\\
&r_1(k_R,\lambda, l)=\frac{1}{2\pi}\int_{-\infty}^{\infty}\!\!\!\!\!\!d\xi\int_{0}^\infty\!\!\!\!\!\!d\eta\, e^{-il\xi+l(l+2k_R)\eta} q_0(\xi,\eta)\left(e_{\lambda}\breve \mu_1^+\right) (\xi,\eta,0,k_R,\lambda), \label{r1}\\
&p_1(k_R,l)=\frac{1}{2\pi}\int_{-\infty}^{\infty}d\xi\int_0^Td\tau\, e^{-il\xi-4il(l^2+3k_Rl+3k_R^2)\tau}H(\xi,\tau,k_R,l)\phi^-_1(\xi,\tau,k_R,0)\label{p1},
\end{align}
and $h$ given by equation \eqref{H}.
\end{prop}

\begin{PROOF} 
%The proof of the above proposition is quite complicated and involves several non-trivial steps; in particular, the starting point \eqref{mft} for our derivation below and the transformation \eqref{trans} on $k$ are not obvious at first glance.
   
Evaluating equations \eqref{direct1ii} at $k_I=0$ and subtracting the resulting expressions yields the following equation for $\Delta \mu_1$:
{\small\begin{align*}
&\Delta \mu_1= \frac{1}{2\pi}\mbox{\large $\Big($}\int_{-\infty}^{0}dl+\int_{-2k_R}^{\infty}dl\mbox{\large $\Big)$}\int_{-\infty}^{\infty}d\xi\int_0^y d\eta \ e^{-il(\xi-x)+l(l+2k_R)(\eta-y)}q\Delta \mu_1\nonumber\\
&+\frac{1}{2\pi}\mbox{\large $\Big($}\int_{-\infty}^{0}\!\!\!\!\!\!dl+\int_{-2k_R}^{\infty}\!\!\!\!\!dl\mbox{\large $\Big)$}\!\!\int_{-\infty}^{\infty}\!\!\!d\xi\int_0^t\!\!d\tau\,e^{-il(\xi-x)-l(l+2k_R)y-4il(l^2+3k_Rl+3k_R^2)(\tau-t)}\mbox{\large $\big($}H\phi^+_1\mbox{\large $\big)$}\mbox{\large $\big|_{k_I=0}$}\nonumber\\
&+\frac{1}{2\pi}\mbox{\large $\Big($}\int_{-\infty}^{0}\!\!\!\!\!\!dl+\int_{-2k_R}^{\infty}\!\!\!\!\!\!dl\mbox{\large $\Big)$}\!\!\int_{-\infty}^{\infty}\!\!\!d\xi\int_t^T\!\!d\tau\,e^{-il(\xi-x)-l(l+2k_R)y-4il(l^2+3k_Rl+3k_R^2)(\tau-t)}\mbox{\large $\big($}H\phi^-_1\mbox{\large $\big)$}\mbox{\large $\big|_{k_I=0}$}\nonumber\\
&-\frac{1}{2\pi}\int_{0}^{-2k_R}dl\int_{-\infty}^{\infty}d\xi\int_{y}^{\infty}d\eta\ e^{-il(\xi-x)+l(l+2k_R)(\eta-y)}q\Delta\mu_1, \quad k_R\leq 0.
\end{align*}}

Furthermore, using the identity \eqref{tT}, we can rewrite this expression in the following form:
{\small\begin{align}
&\Delta \mu_1= \frac{1}{2\pi}\mbox{\large $\Big($}\int_{-\infty}^{0}dl+\int_{-2k_R}^{\infty}dl\mbox{\large $\Big)$}\int_{-\infty}^{\infty}d\xi\int_0^y d\eta \ e^{-il(\xi-x)+l(l+2k_R)(\eta-y)}q\Delta \mu_1\nonumber\\
&+\frac{1}{2\pi}\mbox{\large $\Big($}\int_{-\infty}^{0}\!\!\!\!\!dl+\int_{-2k_R}^{\infty}\!\!\!\!dl\mbox{\large $\Big)$}\!\int_{-\infty}^{\infty}\!\!\!d\xi\int_0^t\!\!d\tau\,e^{-il(\xi-x)-l(l+2k_R)y-4il(l^2+3k_Rl+3k_R^2)(\tau-t)}H\mbox{\large $\big|_{k_I=0}$}\ \Delta\phi_1\nonumber\\
&-\frac{1}{2\pi}\int_{0}^{-2k_R}dl\int_{-\infty}^{\infty}d\xi\int_{y}^{\infty}d\eta\ e^{-il(\xi-x)+l(l+2k_R)(\eta-y)}q\Delta\mu_1\nonumber\\
&+\frac{1}{2\pi}\mbox{\large $\Big($}\int_{-\infty}^{0}\!\!\!\!\!\!dl+\int_{-2k_R}^{\infty}\!\!\!\!\!\!dl\mbox{\large $\Big)$}\!\!\int_{-\infty}^{\infty}\!\!\!d\xi\int_0^T\!\!d\tau\,e^{-il(\xi-x)-l(l+2k_R)y-4il(l^2+3k_Rl+3k_R^2)(\tau-t)}\mbox{\large $\big($}H\phi^-_1\mbox{\large $\big)$}\!\mbox{\large $\big|_{k_I=0}$}.\label{Delta1}
\end{align}}
In what follows, by analysing the Lax pair, we will construct an equation similar to equation \eqref{Delta1}. This equation will be satified by an appropriate combination of the functions $\mu_1^+$ and $\mu_2^+$. Therefore, assuming uniqueness, the solution $\Delta\mu_1$ of equation \eqref{Delta1} equals the above combination.

Consider the Fourier transform pair \eqref{ftpair}, but now let $l=s-2k_R$. This transformation will allow us, after applying a suitable change of variables, to split the exponential terms into two parts: one which will match the exponentials on the RHS of equation \eqref{Delta1} and one which will be \textit{independent} of $l$:

\begin{subequations}
\begin{equation}\label{mft}
\hat \mu(l,y,t,k_R,k_I)\Big|_{l=s-2k_R}\!\!=\int^{\infty}_{-\infty}dx\, e^{-i(s-2k_R)x}\Big(\mu(x,y,t,k_R,k_I)-1\Big), \ s \in \mathbb R, \ k_R\in \mathbb R
\end{equation}
and 
\begin{equation}\label{mift}
 \mu(x,y,t,k_R,k_I)=1+\frac{1}{2\pi}\int^{\infty}_{-\infty}ds\, e^{i(s-2k_R)x}\hat \mu(l,y,t,k_R,k_I)\Big|_{l=s-2k_R}, \quad k_R\in \mathbb R.
\end{equation}
\end{subequations}
We now follow a similar approach to the one presented in section \ref{kpiidp}: applying the transform \eqref{mft} on the Lax equation \eqref{lax1intro}, we find
\begin{equation}\label{m integrating factor y}
\bigg(\hat \mu_y+(s-2k_R)\left(s+2ik_I\right)\hat \mu\bigg)\Big|_{l=s-2k_R}=\int_{-\infty}^{\infty}dx\, e^{-i(s-2k_R)x}\, q\mu.
\end{equation}
Using the change of variables
\begin{equation}\label{trans}
\left.\begin{array}{ll}k_R &\mapsto -k_R-\dfrac{\lambda}{2}  \\ k_I&\mapsto \dfrac{i\lambda}{2}\end{array}\right\} \Rightarrow k\mapsto -k_R-\lambda,\quad \lambda \in \mathbb R,
\end{equation}
equation \eqref{m integrating factor y} becomes:
\begin{equation*}
\bigg(\hat {\breve\mu}_y+(s+2k_R+\lambda)\left(s-\lambda\right)\hat {\breve\mu}\bigg)\Big|_{l=s+2k_R+\lambda}=\int_{-\infty}^{\infty}dx\, e^{-i(s+2k_R+\lambda)x} q {\breve\mu}.
\end{equation*}
with $\breve \mu$ defined by equation \eqref{mubreve}. Hence,
\begin{equation*}
\bigg(e^{(s+2k_R+\lambda)(s-\lambda)y}\ \hat {\breve\mu}\Big|_{l=s+2k_R+\lambda}\bigg)_y=e^{(s+2k_R+\lambda)(s-\lambda)y} \widehat{q\breve \mu}\Big|_{l=s+2k_R+\lambda}, 
\end{equation*}
where 
\begin{equation*}
\widehat{q\breve \mu}\Big|_{l=s+2k_R+\lambda}=\int_{-\infty}^{\infty}dx \, e^{-i(s+2k_R+\lambda)x}q\breve \mu.
\end{equation*}
We solve this equation by integrating with respect to $y$ either from $y=0$ or from $y=\infty$:
{\small\begin{align}\label{iffinal}
\hat{\breve \mu}\mbox{\large $\big|$}_{l=s+2k_R+\lambda}\!\!\!=\!\left\{\!\!\!\!\begin{array}{ll} \int_0^y\! d\eta\, e^{(s+2k_R+\lambda)(s-\lambda)(\eta-y)}\widehat{q\breve \mu}\mbox{\large $\big|$}_{l=s+2k_R+\lambda}\!\!+\!e^{-(s+2k_R+\lambda)(s-\lambda)y}\, \hat{\breve\phi}\mbox{\large $\big|$}_{l=s+2k_R+\lambda} \\  \\ -\int_y^\infty d\eta\, e^{(s+2k_R+\lambda)(s-\lambda)(\eta-y)}\widehat{q\breve \mu}\mbox{\large $\Big|$}_{l=s+2k_R+\lambda}, \end{array}\right.
\end{align}}
where we have assumed that $\hat \mu \rightarrow 0$ as $y\rightarrow \infty$.

We are interested in obtaining an expression for the difference $\Delta \mu_1$, thus we assume that $k_R\leq0$. Also, since $\mu_1$ must be bounded for $k_R\leq0$, the exponential involved in equation \eqref{iffinal} must be bounded and hence, we need to control the sign of $(s+2k_R+\lambda)(s-\lambda)(\eta-y)$. We consider two cases for $\lambda$, namely $\lambda\leq-k_R$ or $\lambda\geq-k_R$. According to the definition \eqref{spectralfunctionsii}, the distinction between $\mu_1$ and $\mu_2$ is made on the basis of the sign of the real part of $k$. However, the transformation \eqref{trans} implies that $k=-k_R-\lambda \in \mathbb R$, hence $\mu=\mu_2$ in the first case, whereas $\mu=\mu_1$ in the second case.

Inverting the expression \eqref{iffinal} (by making use of equation \eqref{mift}), recalling the definition \eqref{mubreve} and relabeling $s$ to $l$, we find the following equations:
{\small\begin{align*}
&\breve \mu_2=1+\frac{1}{2\pi}\mbox{\large $\Big($}\int_{-\infty}^{\lambda}\!dl+\int_{-2k_R-\lambda}^{\infty}\!\!\!dl\mbox{\large $\Big)$}\!\int_{-\infty}^{\infty}\!\!\!d\xi\int_{0}^{y}\!d\eta\, e^{-i(l+2k_R+\lambda)(\xi-x)+(l+2k_R+\lambda)(l-\lambda)(\eta-y)}q\breve \mu_2\nonumber\\
&+\frac{1}{2\pi}\mbox{\large $\Big($}\int_{-\infty}^{\lambda}dl+\int_{-2k_R-\lambda}^{\infty}dl\mbox{\large $\Big)$}\int_{-\infty}^{\infty}d\xi\, e^{-i(l+2k_R+\lambda)(\xi-x)-(l+2k_R+\lambda)(l-\lambda)y}\ \breve \phi_2\nonumber\\
&-\frac{1}{2\pi}\int^{-2k_R-\lambda}_{\lambda}dl\int_{-\infty}^{\infty}d\xi\int_{y}^{\infty}d\eta\, e^{-i(l+2k_R+\lambda)(\xi-x)+(l+2k_R+\lambda)(l-\lambda)(\eta-y)}q\breve \mu_2,
\end{align*}}
which is valid for $\lambda\leq-k_R,\, k_R\leq0,$ and
{\small\begin{align*}
&\breve \mu_1=1+\frac{1}{2\pi}\mbox{\large $\Big($}\!\int_{-\infty}^{-2k_R-\lambda}\!dl+\int_{\lambda}^{\infty}\!\!\!dl\mbox{\large $\Big)$}\!\int_{-\infty}^{\infty}\!\!\!d\xi\int_{0}^{y}\!\!d\eta\, e^{-i(l+2k_R+\lambda)(\xi-x)+(l+2k_R+\lambda)(l-\lambda)(\eta-y)}q\breve \mu_1\nonumber\\
&+\frac{1}{2\pi}\mbox{\large $\Big($}\int_{-\infty}^{-2k_R-\lambda}dl+\int_{\lambda}^{\infty}dl\mbox{\large $\Big)$}\int_{-\infty}^{\infty}d\xi\, e^{-i(l+2k_R+\lambda)(\xi-x)-(l+2k_R+\lambda)(l-\lambda)y}\ \breve \phi_1\nonumber\\
&-\frac{1}{2\pi}\int^{\lambda}_{-2k_R-\lambda}dl\int_{-\infty}^{\infty}d\xi\int_{y}^{\infty}d\eta\, e^{-i(l+2k_R+\lambda)(\xi-x)+(l+2k_R+\lambda)(l-\lambda)(\eta-y)}q\breve \mu_1,
\end{align*}}
which is valid for $\lambda\geq-k_R, \, k_R\leq0.$

Following the approach used in section \ref{kpiidp}, we can employ the second Lax equation \eqref{lax2intro} in order to eliminate the dependence in $y$ from the second term on the RHS of the above expressions. This yields the following equations:
{\small\begin{align}\label{mu2afinal}
&e_{\lambda}\, \breve \mu_2^+=e_{\lambda}+\frac{1}{2\pi}\mbox{\large $\Big($}\int_{-\infty}^{\lambda}dl+\int_{-2k_R-\lambda}^{\infty}dl\mbox{\large $\Big)$}\int_{-\infty}^{\infty}d\xi\int_{0}^{y}d\eta\, e^{-il(\xi-x)+l(l+2k_R)(\eta-y)}q\,e_\lambda\,\breve \mu_2^+\nonumber\\
&+\frac{1}{2\pi}\mbox{\large $\Big($}\int_{-\infty}^{\lambda}\!\!\!\!dl+\!\!\int_{-2k_R-\lambda}^{\infty}\!\!\!\!\!\!dl\mbox{\large $\Big)$}\!\int_{-\infty}^{\infty}\!\!\!d\xi\!\int_0^t\!\!\!d\tau\, e^{-il(\xi-x)-l(l+2k_R)y-4il(l^2+3k_Rl+3k_R^2)(\tau-t)}\breve H\,e_\lambda|_{\eta=0}\ \breve \phi_2^+\nonumber\\
&-\frac{1}{2\pi}\int^{-2k_R-\lambda}_{\lambda}dl\int_{-\infty}^{\infty}d\xi\int_{y}^{\infty}d\eta\,e^{-il(\xi-x)+l(l+2k_R)(\eta-y)}q\,e_\lambda\,\breve \mu_2^+,
\end{align}}
which is valid for $\lambda\leq-k_R, \, k_R\leq0,$ and
{\small\begin{align}\label{mu2bfinal}
&e_{\lambda}\, \breve \mu_1^+=e_{\lambda}+\frac{1}{2\pi}\mbox{\large $\Big($}\int_{-\infty}^{-2k_R-\lambda}dl+\int_{\lambda}^{\infty}dl\mbox{\large $\Big)$}\int_{-\infty}^{\infty}d\xi\int_{0}^{y}d\eta\, e^{-il(\xi-x)+l(l+2k_R)(\eta-y)}q\,e_\lambda\,\breve \mu_1^+\nonumber\\
&+\frac{1}{2\pi}\mbox{\large $\Big($}\int_{-\infty}^{-2k_R-\lambda}\!\!\!\!\!\!\!dl+\!\!\int_{\lambda}^{\infty}\!\!\!dl\mbox{\large $\Big)$}\!\!\int_{-\infty}^{\infty}\!\!\!d\xi\!\int_0^t\!\!\!d\tau\, e^{-il(\xi-x)-l(l+2k_R)y-4il(l^2+3k_Rl+3k_R^2)(\tau-t)}\breve H\,e_\lambda|_{\eta=0}\, \breve \phi_1^+\nonumber\\
&-\frac{1}{2\pi}\int^{\lambda}_{-2k_R-\lambda}dl\int_{-\infty}^{\infty}d\xi\int_{y}^{\infty}d\eta\,e^{-il(\xi-x)+l(l+2k_R)(\eta-y)}q\,e_\lambda\,\breve \mu_1^+,
\end{align}}
which is valid for $\lambda\geq-k_R, \, k_R\leq0$, where $e_{\lambda}$ is defined by equation \eqref{elambda}, $\breve H$ is given by
\begin{equation}\label{Hmod}
\breve H(x,t,k_R,l):=H(x,t,-k_R-\lambda, l+2k_R+\lambda)
\end{equation}
and $H$ is defined by equation \eqref{H}.

In order for $\breve \mu_{1,2}$ to satisfy an equation similar to equation \eqref{Delta1}, we need to split the integrals with respect to $l$. However, in order to ensure that the resulting integrals remain bounded, we impose two additional conditions: in equation \eqref{mu2afinal}, where $\lambda\leq-k_R$, we actually require the stronger condition that $\lambda\leq0$, and in equation \eqref{mu2bfinal}, where $\lambda\geq-k_R$, we demand the stronger condition that $\lambda\geq-2k_R$. Under these new restrictions, the expressions \eqref{mu2afinal} and \eqref{mu2bfinal} become
{\small\begin{align}\label{mu2a}
&e_{\lambda}\breve \mu_2^+=e_{\lambda}-R_2\nonumber\\
&+\frac{1}{2\pi}\mbox{\large $\Big($}\int_{-\infty}^{0}dl+\int_{-2k_R}^{\infty}dl\mbox{\large $\Big)$}\int_{-\infty}^{\infty}d\xi\int_{0}^{y}d\eta\, e^{-il(\xi-x)+l(l+2k_R)(\eta-y)}q\,e_{\lambda}\,\breve \mu_2^+\nonumber\\
&+\frac{1}{2\pi}\mbox{\large $\Big($}\int_{-\infty}^{0}\!\!\!dl+\!\int_{-2k_R}^{\infty}\!\!\!dl\mbox{\large $\Big)$}\!\!\int_{-\infty}^{\infty}\!\!\!d\xi\int_0^t \!\!d\tau\, e^{-il(\xi-x)-l(l+2k_R)y-4il(l^2+3k_Rl+3k_R^2)(\tau-t)}\,\breve H\,e_\lambda|_{\eta=0}\,\breve \phi_2^+\nonumber\\
&-\frac{1}{2\pi}\int_{0}^{-2k_R}dl\int_{-\infty}^{\infty}d\xi\int_{y}^{\infty}d\eta\,e^{-il(\xi-x)+l(l+2k_R)(\eta-y)}q\,e_{\lambda}\,\breve \mu_2^+,
\end{align}}
which is valid for $\lambda\leq0, \, k_R\leq0,$ and
{\small\begin{align}\label{mu2b}
&e_{\lambda}\breve \mu_1^+=e_{\lambda}-R_1\nonumber\\
&+\frac{1}{2\pi}\mbox{\large $\Big($}\int_{-\infty}^{0}dl+\int_{-2k_R}^{\infty}dl\mbox{\large $\Big)$}\int_{-\infty}^{\infty}d\xi\int_{0}^{y}d\eta\, e^{-il(\xi-x)+l(l+2k_R)(\eta-y)}q\,e_{\lambda}\,\breve \mu_1^+\nonumber\\
&+\frac{1}{2\pi}\mbox{\large $\Big($}\int_{-\infty}^{0}\!\!\!dl+\!\int_{-2k_R}^{\infty}\!\!\!dl\mbox{\large $\Big)$}\!\!\int_{-\infty}^{\infty}\!\!\!d\xi\int_0^t \!\!d\tau\, e^{-il(\xi-x)-l(l+2k_R)y-4il(l^2+3k_Rl+3k_R^2)(\tau-t)}\,\breve H\,e_\lambda|_{\eta=0}\,\breve \phi_1^+\nonumber\\
&-\frac{1}{2\pi}\int_{0}^{-2k_R}dl\int_{-\infty}^{\infty}d\xi\int_{y}^{\infty}d\eta\,e^{-il(\xi-x)+l(l+2k_R)(\eta-y)}q\,e_{\lambda}\,\breve \mu_1^+,
\end{align}}
which is valid for $\lambda\geq-2k_R, \, k_R\leq0$.

The quantities $R_2$ and $R_1$ are given by the following expressions:
{\small\begin{align}
&R_2(x,y,t,k_R,\lambda)=\frac{1}{2\pi}\mbox{\large $\Big($}\int_{\lambda}^0 \!\!dl+\int^{-2k_R-\lambda}_{-2k_R}\!\!\!dl\mbox{\large $\Big)$}\!\int_{-\infty}^{\infty}d\xi\int_{0}^{\infty}d\eta\, e^{-il(\xi-x)+l(l+2k_R)(\eta-y)}q\,e_{\lambda}\,\breve \mu_2^+\nonumber\\
&+\frac{1}{2\pi}\mbox{\large $\Big($}\!\int^{0}_{\lambda}\!\!dl+\!\!\int^{-2k_R-\lambda}_{-2k_R}\!\!\!\!\!dl\mbox{\large $\Big)$}\!\!\int_{-\infty}^{\infty}\!\!\!\!d\xi\int_0^t \!\!\!\!d\tau\, e^{-il(\xi-x)-l(l+2k_R)y-4il(l^2+3k_Rl+3k_R^2)(\tau-t)}\breve H\,e_\lambda|_{\eta=0}\,\breve \phi_2^+,\label{R1}
\end{align}}
which is defined for $\lambda\leq0, \, k_R\leq0$, and
{\small\begin{align}
&R_1(x,y,t,k_R,\lambda)=\frac{1}{2\pi}\mbox{\large $\Big($}\int^{0}_{-2k_R-\lambda}\!\!\!dl+\int^{\lambda}_{-2k_R}\!\!\!dl\mbox{\large $\Big)$}\!\int_{-\infty}^{\infty}\!\!d\xi\int_{0}^{\infty}\!d\eta\, e^{-il(\xi-x)+l(l+2k_R)(\eta-y)}q\,e_{\lambda}\,\breve \mu_1^+\nonumber\\
&+\frac{1}{2\pi}\mbox{\large $\Big($}\!\int^{0}_{-2k_R-\lambda}\!\!\!\!\!\!dl+\!\!\int^{\lambda}_{-2k_R}\!\!\!\!\!dl\mbox{\large $\Big)$}\!\!\int_{-\infty}^{\infty}\!\!\!\!d\xi\!\int_0^t\!\!\! d\tau\, e^{-il(\xi-x)-l(l+2k_R)y-4il(l^2+3k_Rl+3k_R^2)(\tau-t)}\breve H e_\lambda|_{\eta=0}\,\breve \phi_1^+\label{R2},
\end{align}}
which is defined for $\lambda\geq-2k_R, \, k_R\leq0$.

We need to verify that the exponentials involved in the expressions \eqref{mu2a}, \eqref{mu2b}, \eqref{R1} and \eqref{R2} are bounded. This is straightforward for the exponentials $e^{-l(l+2k_R)y}$, $e^{l(l+2k_R)(\eta-y)}$ and $e_\lambda$ in \eqref{mu2a} and \eqref{mu2b}.

Concerning $R_2$, we have  $$\lambda\leq l\leq0\leq-2k_R \quad \mathrm{or}\quad 0\leq-2k_R\leq l\leq-2k_R-\lambda.$$ Thus, the exponentials $e^{-l(l+2k_R)y}$ and $e^{(l-\lambda)(l+2k_R+\lambda)\eta}$, which are involved in the definition of $R_2$, are bounded.

Concerning $R_1$, we have $$-2k_R-\lambda\leq l\leq0\leq-2k_R \quad \mathrm{or}\quad 0\leq-2k_R\leq l\leq\lambda.$$ Thus, the exponentials $e^{-l(l+2k_R)y}$ and $e^{(l-\lambda)(l+2k_R+\lambda)\eta}$, which are involved in the definition of $R_1$, are also bounded.

The expressions for $R_2$ and $R_1$ can be rewritten in the following form:
\begin{align}
R_2&=\frac{1}{2\pi}\left(\int^{0}_{\lambda}dl+\int^{-2k_R-\lambda}_{-2k_R}dl\right)e^{ilx-l(l+2k_R)y+4il(l^2+3k_Rl+3k_R^2)t}\nonumber\\
&\times\Bigg\{\int_{-\infty}^{\infty}d\xi\int_{0}^{\infty}d\eta\,e^{-il\xi+l(l+2k_R)\eta-4il(l^2+3k_Rl+3k_R^2)t}q\,e_{\lambda}\,\breve \mu_2^+\nonumber\\
&+\int_{-\infty}^{\infty}d\xi\int_0^t d\tau\, e^{-il\xi-4il(l^2+3k_Rl+3k_R^2)\tau}\,\breve H\,e_\lambda|_{\eta=0}\,\breve \phi_2^+\Bigg\}\label{R1mod}
\end{align}
and 
\begin{align}
R_1&=\frac{1}{2\pi}\left(\int^{0}_{-2k_R-\lambda}dl+\int^{\lambda}_{-2k_R}dl\right)e^{ilx-l(l+2k_R)y+4il(l^2+3k_Rl+3k_R^2)t}\nonumber\\
&\times\Bigg\{\int_{-\infty}^{\infty}d\xi\int_{0}^{\infty}d\eta\,e^{-il\xi+l(l+2k_R)\eta-4il(l^2+3k_Rl+3k_R^2)t}\,q\,e_{\lambda}\,\breve \mu_1^+\nonumber\\
&+\int_{-\infty}^{\infty}d\xi\int_0^t d\tau\, e^{-il\xi-4il(l^2+3k_Rl+3k_R^2)\tau}\,\breve H\,e_\lambda|_{\eta=0}\,\breve \phi_1^+\Bigg\}\label{R2mod}.
\end{align}
Applying the transformation \eqref{trans} to the identity \eqref{GRjump}, and then replacing $l$ by $l+\lambda+2k_R$, we find the following identity:
\begin{align*}
&\int_{-\infty}^{\infty}d\xi\int_{0}^\infty d\eta\, e^{-il\xi+l(l+2k_R)\eta}\, q_0\left(e_\lambda \breve \mu_j^+\right)\!\Big|_{t=0}\nonumber\\
&=\int_{-\infty}^{\infty}d\xi\int_{0}^\infty d\eta\, e^{-il\xi+l(l+2k_R)\eta-4il(l^2+3k_Rl+3k_R^2)t}\, q\,e_\lambda\breve \mu_j^+ \nonumber\\
&+\int_{0}^{t}d\tau\int_{-\infty}^{\infty}d\xi\, e^{-il\xi-4il(l^2+3k_Rl+3k_R^2)\tau}H\,e_\lambda|_{\eta=0}\, \breve \phi_j^+, \quad j=1,2.
\end{align*}
Note that if $l$ satisfies the restrictions imposed in the domain of $R_1$ and $R_2$ then the exponential $e^{(l-\lambda)(l+2k_R+\lambda)\eta}$ is bounded, hence the global relation \eqref{GRjump} is satisfied. Thus, employing the global relation in the expressions \eqref{R1mod} and \eqref{R2mod}, we find
\begin{equation}\label{R1final}
R_2(x,y,t,k_R,\lambda)=\left(\int^{0}_{\lambda}dl+\int_{-2k_R-\lambda}^{-2k_R}dl\right)E(x,y,t,k_R,l)r_2(k_R,\lambda, l),\quad \lambda\leq 0,
\end{equation}
and
\begin{equation}\label{R2final}
R_1(x,y,t,k_R,\lambda)=\left(\int^{0}_{-2k_R-\lambda}dl+\int_{\lambda}^{-2k_R}\!\!dl\right) E(x,y,t,k_R,l)r_1(k_R, \lambda, l),\ \lambda\geq -2k_R,
\end{equation}
where
\begin{equation}
E(x,y,t,k,l):=e^{ilx-l(l+2k)y+4il(l^2+3kl+3k^2)t},\label{E}
\end{equation}
and $r_2$ and $r_1$ are defined by equations \eqref{r2} and \eqref{r1}.

In summary, we have obtained two integral equations for $\breve \mu_2^+$ and $\breve \mu_1^+$, namely equations \eqref{mu2a} and \eqref{mu2b}, with the quantities $R_2$ and $R_1$ given by equations \eqref{R1final} and \eqref{R2final} respectively.

Next, we multiply equation \eqref{mu2a} by $\chi_2$ and integrate with respect to $\lambda$ from $-\infty$ to 0; similarly, we multiply equation \eqref{mu2b} by $\chi_1$ and integrate with respect to $\lambda$ from $-2k_R$ to $\infty$ (the precise form of the functions $\chi_1$ and $\chi_2$ will be determined later). Adding the two resulting expressions and defining $\mathcal S_1$ by
\begin{equation}\label{S1}
\mathcal S_1(x,y,t,k_R)=\int_{-\infty}^{0}d\lambda\chi_2(k_R,\lambda)\,e_{\lambda}\,\breve \mu_2^++\int_{-2k_R}^{\infty}d\lambda\chi_1(k_R,\lambda)\,e_{\lambda}\,\breve \mu_1^+,
\end{equation}
we find
\begin{align}\label{inteqS1}
&\mathcal S_1=\frac{1}{2\pi}\left(\int_{-\infty}^{0}dl+\int_{-2k_R}^{\infty}dl\right)\int_{-\infty}^{\infty}d\xi\int_{0}^{y}d\eta\, e^{-il(\xi-x)+l(l+2k_R)(\eta-y)}q\, \mathcal S_1\nonumber\\
&+\frac{1}{2\pi}\left(\int_{-\infty}^{0}\!\!\!\!dl+\!\!\!\int_{-2k_R}^{\infty}\!\!\!dl\right)\!\!\int_{-\infty}^{\infty}\!\!\!\!\!d\xi\int_0^t \!\!\!d\tau\, e^{-il(\xi-x)-l(l+2k_R)y-4il(l^2+3k_Rl+3k_R^2)(\tau-t)}\breve H\mathcal { S}_1|_{\eta=0}\nonumber\\
&-\frac{1}{2\pi}\int_{0}^{-2k_R}dl\int_{-\infty}^{\infty}d\xi\int_{y}^{\infty}d\eta\,e^{-il(\xi-x)+l(l+2k_R)(\eta-y)}q\, \mathcal S_1\nonumber\\
&+\mathcal F_1(x,y,t,k_R),
\end{align}
where $\breve H$ is defined by equation \eqref{Hmod}, and 
\begin{equation*}
\mathcal F_1(x,y,t,k_R)=\int_{-\infty}^{0}d\lambda\chi_2(k_R,\lambda)\left(e_{\lambda}-R_2\right)+\int_{-2k_R}^{\infty}d\lambda\chi_1(k_R,\lambda)\left(e_{\lambda}-R_1\right).
\end{equation*}
In order to identify equation \eqref{inteqS1} with equation \eqref{Delta1}, the functions $\chi_1$ and $\chi_2$ must be defined in such way that equations \eqref{inteqS1} and \eqref{Delta1} have the same forcing, i.e. so that
\begin{equation*}
\mathcal F_1(x,y,t,k_R)=\left(\int_{-\infty}^{0}dl+\int_{-2k_R}^{\infty}dl\right)\, E(x,y,t,k_R,l)p_1(k_R,l),
\end{equation*}
where $E$ is defined by equation \eqref{E} and $p_1$ is given by equation \eqref{p1}. This implies that $\chi_1$ and $\chi_2$ must satisfy the following condition:
{\small\begin{align}\label{comparison1}
&\int_{-\infty}^{0}d\lambda\chi_2(k_R,\lambda)\,e_\lambda(x,y,t,k_R,\lambda)-\int_{-\infty}^{0}d\lambda \chi_2(k_R,\lambda)\int^0_\lambda dl\, E(x,y,t,k_R,l)r_2(k_R, \lambda, l)\nonumber\\
&-\int_{-\infty}^{0}d\lambda \chi_2(k_R,\lambda)\int^{-2k_R-\lambda}_{-2k_R} dl \, E(x,y,t,k_R, l)r_2(k_R, \lambda, l)\nonumber\\
&+\int_{-2k_R}^\infty d\lambda \chi_1(k_R,\lambda)\,e_\lambda(x,y,t,k_R,\lambda)\nonumber\\
&-\int_{-2k_R}^{\infty}d\lambda \chi_1(k_R, \lambda) \int^{0}_{-2k_R-\lambda}dl\, E(x,y,t,k_R,l)r_1(k_R,\lambda, l)\nonumber\\
&-\int_{-2k_R}^{\infty}d\lambda \chi_1(k_R, \lambda) \int^{\lambda}_{-2k_R}dl\, E(x,y,t,k_R,l)r_1(k_R,\lambda, l)=\nonumber\\
&=\int_{-\infty}^{0}d\lambda\, E(x,y,t,k_R,\lambda)p_1(k_R,\lambda)+\int_{-2k_R}^{\infty}d\lambda\, E(x,y,t,k_R,\lambda)p_1(k_R,\lambda),
\end{align}}
where $l$ has been relabeled as $\lambda$ on the RHS.

Our aim is to obtain appropriate integral equations for $\chi_1$ and $\chi_2$. In this respect, we note that we can interchange the order of integration with respect to $\lambda$ and $l$:
\begin{equation*}
\int_{-\infty}^0d\lambda\int_\lambda^0 dl=\int_{-\infty}^{0}dl\int_{-\infty}^{l}d\lambda,
\end{equation*}
\begin{equation*}
\int_{-\infty}^{0}d\lambda\int^{-2k_R-\lambda}_{-2k_R}dl=\int_{-2k_R}^{\infty}dl\int_{-\infty}^{-2k_R-l}d\lambda,
\end{equation*}
\begin{equation*}
\int_{-2k_R}^{\infty}d\lambda\int^0_{-2k_R-\lambda} dl=\int_{-\infty}^{0}dl\int_{-2k_R-l}^{\infty}d\lambda,
\end{equation*}
and
\begin{equation*}
\int_{-2k_R}^{\infty}d\lambda\int^{\lambda}_{-2k_R}dl=\int_{-2k_R}^{\infty}dl\int_{l}^\infty d\lambda.
\end{equation*}
\begin{figure}[ht]
\begin{center}
\resizebox{7cm}{!}{\input{int1_2.pstex_t}}
\end{center}
%\caption{Interchanging the order of integration between $\lambda$ and $l$.}
\label{int1_2}
\end{figure}
\begin{figure}[ht]
\begin{center}
\resizebox{9cm}{!}{\input{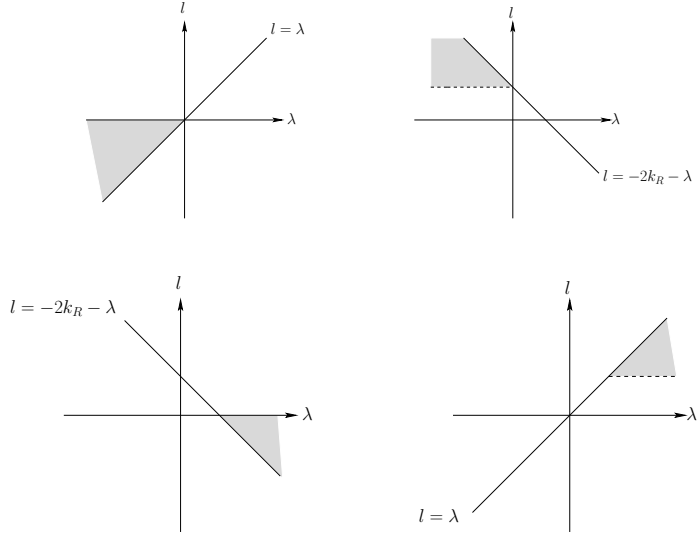}}
\end{center}
\caption{Interchanging the order of integration between $\lambda$ and $l$.}
\label{int3_4}
\end{figure}
Therefore, equation \eqref{comparison1} becomes:
{\small\begin{align*}
&\int_{-\infty}^0dl\chi_2(k_R,l)\,e_l(x,y,t,k_R,l)-\int_{-\infty}^{0}dl\, E(x,y,t,k_R,l)\int_{-\infty}^{l}d\lambda \chi_2(k_R,\lambda)r_2(k_R, \lambda, l)\nonumber\\
&-\int_{-2k_R}^{\infty}dl\, E(x,y,t,k_R, l)\int_{-\infty}^{-2k_R-l}d\lambda \chi_2(k_R,\lambda) r_2(k_R, \lambda, l)\nonumber\\
&+\int_{-2k_R}^\infty\!\! dl \chi_1(k_R,l)\, e_l(x,y,t,k_R,l)-\int_{-\infty}^{0}\!\!dl \, E(x,y,t,k_R,l)\int_{-2k_R-l}^{\infty}d\lambda \chi_1(k_R, \lambda)r_1(k_R,\lambda, l)\nonumber\\
&-\int_{-2k_R}^{\infty}dl\, E(x,y,t,k_R,l)\int_{l}^\infty d\lambda\chi_1(k_R, \lambda)r_1(k_R,\lambda, l)=\nonumber\\
&=\int_{-\infty}^{0}dl\, E(x,y,t,k_R,l)p_1(k_R,l)+\int_{-2k_R}^{\infty}dl\, E(x,y,t,k_R,l)p_1(k_R,l).
\end{align*}}
In addition, it follows from the definitions \eqref{elambda} and \eqref{E} of $e_\lambda$ and $E$, that
\begin{equation}\label{Ee}
 E(x,y,t,k_R, -\lambda-2k_R)=e_\lambda(x,y,t,k_R,\lambda).
\end{equation}
Hence the transformation $l \mapsto -2k_R-l$ in the terms involving the exponential $E$ yields
\begin{align}\label{comparison3}
&\int_{-\infty}^0dl\, e_l(x,y,t,k_R,l)\Bigg\{\chi_2(k_R,l)-\int_{-\infty}^{l}d\lambda \chi_2(k_R,\lambda) r_2(k_R, \lambda, -2k_R-l)\nonumber\\
&-\int_{-2k_R-l}^\infty d\lambda\chi_1(k_R, \lambda)r_1(k_R,\lambda,-2k_R-l)-p_1(k_R,-2k_R-l)\Bigg\}\nonumber\\
&+\int_{-2k_R}^{\infty}dl\, e_l(x,y,t,k_R,l)\Bigg\{\chi_1(k_R,l)-\int_{-\infty}^{-2k_R-l}d\lambda \chi_2(k_R,\lambda)r_2(k_R, \lambda, -2k_R-l)\nonumber\\
&-\int_{l}^{\infty}d\lambda \chi_1(k_R, \lambda)r_1(k_R,\lambda, -2k_R-l)-p_1(k_R,-2k_R-l)\Bigg\}=0.
\end{align}
This condition can be identically satisfied provided that we set both curly brackets equal to zero, i.e. provided that
\begin{align*}%\label{chi1eq1}
&\chi_2(k_R,\lambda)-\int_{-\infty}^{\lambda}dl\chi_2(k_R,l) r_2(k_R, l, -2k_R-\lambda)\nonumber\\
&-\int_{-2k_R-\lambda}^\infty dl\chi_1(k_R, l)r_1(k_R,\lambda,-2k_R-\lambda)=p_1(k_R,-2k_R-\lambda), \quad \lambda\leq0
\end{align*}
and
\begin{align*}%\label{chi2eq1}
&\chi_1(k_R,\lambda)-\int_{-\infty}^{-2k_R-\lambda}dl \chi_2(k_R,l)r_2(k_R, l, -2k_R-\lambda)\nonumber\\
&-\int_{\lambda}^{\infty}dl \chi_1(k_R, l)r_1(k_R,l, -2k_R-\lambda)=p_1(k_R,-2k_R-\lambda),\quad \lambda\geq-2k_R,
\end{align*}
which are equations \eqref{chi2eqfinal} and \eqref{chi1eqfinal} respectively. 

Note that the restrictions on $\lambda$ originate from equation \eqref{comparison3}; furthermore, these restrictions ensure that $r_2$ and $r_1$ are well defined, see equations \eqref{r2} and \eqref{r1}.

Equations \eqref{psi2eqfinal} and \eqref{psi1eqfinal} can be derived in a similar way. 
\QED
\end{PROOF}

\subsection{The spectral functions}\label{kpiiir}

Consider the initial-boundary value problem for the KPII equation \eqref{kpiiintro}, with the initial and the boundary values defined by equations \eqref{icintro}-\eqref{lkpbc2intro}. This problem is well posed provided that one of the two boundary values, \eqref{lkpbc1intro} or \eqref{lkpbc2intro}, is known. It is then possible to obtain the unknown boundary value in terms of the initial condition and the given boundary condition via the global relation \eqref{GRjump}.

Now, we will define a map from the initial condition $q_0$ and from the boundary values $g$ and $h$ to the spectral functions of propositions \ref{kpiid-barprop} and \ref{kpiidiscprop}:
\begin{align}\label{map}
\{q_0(x,y),\ g(x,t),\ h(x,t)\}\mapsto\{\alpha_{1,2}^-(k_R,k_I), \beta_{1,2}^\pm (k_R,k_I), \chi_{1,2}(k_R,\lambda),\psi_{1,2}(k_R,\lambda)\},
\end{align}
where $(x,y,t) \in \Omega$, $k_R\leq0$ for $\{\alpha^-_1, \beta^\pm_1\}$ and $k_R\geq0$ for $\{\alpha^-_2, \beta^\pm_2\}$, $k_I\geq0$ for  $\beta_{1,2}^+$, $k_I\leq0$ for $\{\alpha_{1,2}^-, \beta_{1,2}^-\}$ and $\lambda \in \mathbb R$.

\begin{enumerate}
\item $ q_0 \mapsto \{\rho_1^+(x,y,k_R,k_I),\rho_2^+(x,y,k_R,k_I)\}$, for $(x,y,t)\in \Omega$ and $k_I\geq0$, \newline where the functions $\rho_1^+$ and $\rho_2^+$ are defined in terms of $q_0$ via the following linear integral equations:
\begin{subequations}\label{maprhoa}
\begin{align}\label{mapeq1}
\rho^+_1&=1+\frac{1}{2\pi}\left(\int_{-\infty}^{0}dl+\int_{-2k_R}^{\infty}dl\right)\int_{-\infty}^{\infty}d\xi \int_{0}^{y}d\eta\, e^{-il(\xi-x)+l(l+2k)(\eta-y)}q_0\rho_1^{+}\nonumber\\ 
&-\frac{1}{2\pi}\int_{0}^{-2k_R}dl\int_{-\infty}^{\infty}d\xi\int_{y}^{\infty}d\eta \, e^{-il(\xi-x)+l(l+2k)(\eta-y)}q_0\rho_1^{+}, \quad k_R\leq0,
\end{align}
\begin{align}\label{mapeq2}
\rho^+_2&=1+\frac{1}{2\pi}\left(\int_{-\infty}^{-2k_R}dl+\int_{0}^{\infty}dl\right)\int_{-\infty}^{\infty}d\xi \int_{0}^{y}d\eta\, e^{-il(\xi-x)+l(l+2k)(\eta-y)}q_0\rho_2^{+}\nonumber\\ 
&-\frac{1}{2\pi}\int_{-2k_R}^{0}dl\int_{-\infty}^{\infty}d\xi\int_{y}^{\infty}d\eta \, e^{-il(\xi-x)+l(l+2k)(\eta-y)}q_0\rho_2^{+}, \quad k_R\geq0.
\end{align}
\end{subequations}

\item {\small$\{q_0,g,h\} \mapsto \{\rho_1^-(x,y,k_R,k_I),\rho_2^-(x,y,k_R,k_I)\},\{\phi_1^-(x,t,k_R,k_I), \phi_2^-(x,t,k_R,k_I)\}$},\\
for $(x,y,t) \in \Omega$ and $k_I\leq0$,\\
where the functions $\rho_1^-$ and $\rho_2^-$ are defined in terms of $\{q_0,g,h\}$ via the following linear integral equations:
\begin{subequations}\label{maprhob}
\begin{align}\label{mapeq3}
\rho^-_1&=1+\frac{1}{2\pi}\left(\int_{-\infty}^{0}dl+\int_{-2k_R}^{\infty}dl\right)\int_{-\infty}^{\infty}d\xi \int_{0}^{y}d\eta\, e^{-il(\xi-x)+l(l+2k)(\eta-y)}q_0\rho_1^{-}\nonumber\\ 
&-\frac{1}{2\pi}\left(\int_{-\infty}^{0}\!\!\!dl+\!\int_{-2k_R}^{\infty}\!\!\!dl\right)\!\!\int_{-\infty}^{\infty}\!\!\!\!d\xi\!\int_{0}^{T}\!\!\!\!\!d\tau \, e^{-il(\xi-x)-l(l+2k)y-4il(l^2+3kl+3k^2)\tau} H\phi_1^{-}\nonumber\\
&-\frac{1}{2\pi}\int_{0}^{-2k_R}dl\int_{-\infty}^{\infty}d\xi\int_{y}^{\infty}d\eta\,e^{-il(\xi-x)+l(l+2k)(\eta-y)}q_0\rho_1^{-}, \quad k_R\leq0,
\end{align}
\begin{align}\label{mapeq4}
\rho^-_2&=1+\frac{1}{2\pi}\left(\int_{-\infty}^{-2k_R}dl+\int_{0}^{\infty}dl\right)\int_{-\infty}^{\infty}d\xi \int_{0}^{y}d\eta\, e^{-il(\xi-x)+l(l+2k)(\eta-y)}q_0\rho_2^{-}\nonumber\\ 
&-\frac{1}{2\pi}\left(\int_{-\infty}^{-2k_R}\!\!\!dl+\!\int_{0}^{\infty}\!\!\!dl\right)\!\!\int_{-\infty}^{\infty}\!\!\!\!d\xi\!\int_{0}^{T}\!\!\!\!\!d\tau\, e^{-il(\xi-x)-l(l+2k)y-4il(l^2+3kl+3k^2)\tau} H\phi_2^{-}\nonumber\\
&-\frac{1}{2\pi}\int_{-2k_R}^{0}dl\int_{-\infty}^{\infty}d\xi\int_{y}^{\infty}d\eta \, e^{-il(\xi-x)+l(l+2k)(\eta-y)}q_0\rho_2^{-}, \quad k_R\geq0,
\end{align}
\end{subequations}
and
\begin{subequations}\label{mapphi}
\begin{align}\label{mapeq5}
\phi^-_1&=1-\frac{1}{2\pi}\mbox{\large $\Big($}\int_{-\infty}^{0}\!\!\!dl+\!\!\int_{-2k_R}^{\infty}\!\!\!\!\!dl\mbox{\large $\Big)$}\!\!\int_{-\infty}^{\infty}\!\!\!\!d\xi\!\int_{t}^{T}\!\!\!d\tau  e^{-il(\xi-x)-4il(l^2+3kl+3k^2)(\tau-t)} H\phi_1^{-}\nonumber\\
&+\frac{1}{2\pi}\int_{0}^{-2k_R}dl\int_{-\infty}^{\infty}d\xi\int_{0}^{t}\!\!\!d\tau  e^{-il(\xi-x)-4il(l^2+3kl+3k^2)(\tau-t)} H\phi_1^{-},\nonumber\\
&-\frac{1}{2\pi}\int_{0}^{-2k_R}dl\int_{-\infty}^{\infty}d\xi\int_{0}^{\infty}d\eta \, e^{-il(\xi-x)+l(l+2k)\eta+4il(l^2+3kl+3k^2)t}q_0\rho_1^{-},
\end{align}
\begin{align}\label{mapeq6}
\phi^-_2&=1-\frac{1}{2\pi}\mbox{\large $\Big($}\int_{-\infty}^{-2k_R}\!\!\!dl+\!\!\int_{0}^{\infty}\!\!\!\!\!dl\mbox{\large $\Big)$}\!\!\int_{-\infty}^{\infty}\!\!\!\!d\xi\!\int_{t}^{T}\!\!\!d\tau  e^{-il(\xi-x)-4il(l^2+3kl+3k^2)(\tau-t)} H\phi_2^{-}\nonumber\\
&+\frac{1}{2\pi}\int_{-2k_R}^{0}dl\int_{-\infty}^{\infty}d\xi\int_{0}^{t}\!\!\!d\tau  e^{-il(\xi-x)-4il(l^2+3kl+3k^2)(\tau-t)} H\phi_2^{-},\nonumber\\
&-\frac{1}{2\pi}\int_{-2k_R}^{0}dl\int_{-\infty}^{\infty}d\xi\int_{0}^{\infty}d\eta \, e^{-il(\xi-x)+l(l+2k)\eta+4il(l^2+3kl+3k^2)t}q_0\rho_2^{-},
\end{align}
\end{subequations}
with $H$ defined by equation \eqref{H}.

The motivation for the above definitions emanates from proposition \ref{kpiidpprop}. Indeed, evaluating equations \eqref{direct1ii} and \eqref{direct2ii} at $t=0$ and letting $\mu_{1,2}^\pm\big|_{t=0}=\rho_{1,2}^\pm$, we recover equations \eqref{mapeq1}-\eqref{mapeq4}. Furthermore, applying the inverse Fourier transform \eqref{ift} for $\phi_{1,2}^-$ on equations \eqref{phi1-hat} and \eqref{phi2-hat}, we find equations \eqref{mapeq5} and \eqref{mapeq6}.

\item $\{q_0,g,h\}\mapsto \{\alpha_{1,2}^-(k_R,k_I),\beta_{1,2}^\pm(k_R,k_I)\}$,\\
where the functions $\alpha^-_{1,2}$ and $\beta_{1,2}^\pm$ are defined by the equations:
\begin{subequations}\label{mapsp}
\begin{align}
\alpha_1^-(k_R,k_I)&=\frac{1}{2\pi}\int_{-\infty}^{\infty}d\xi\int_0^T d\tau\, e^{2ik_R\xi+8ik_R(k_R^2-3k_I^2)\tau} H(\xi,\tau,k,-2k_R)\phi_1^-\label{mapeq7},\\
\alpha_2^-(k_R,k_I)&=\frac{1}{2\pi}\int_{-\infty}^{\infty}d\xi\int_0^Td\tau\, e^{2ik_R\xi+8ik_R(k_R^2-3k_I^2)\tau}H(\xi,\tau,k,-2k_R)\phi_2^-\label{mapeq9},\\
\beta_1^{\pm}(k_R,k_I)&=\frac{1}{2\pi} \int_{-\infty}^{\infty}d\xi\int_0^\infty d\eta\, e^{2ik_R\xi-4ik_Rk_I\eta}q_0\rho_1^\pm\label{mapeq8},\\
\beta_2^{\pm}(k_R,k_I)&=\frac{1}{2\pi} \int_{-\infty}^{\infty}d\xi\int_0^\infty d\eta\ e^{2ik_R\xi-4ik_Rk_I\eta}q_0\rho_2^\pm.\label{mapeq10}
\end{align}
\end{subequations}
\item $\{q_0,g,h\}\mapsto\{r_{1,2}(k_R,\lambda,l), p_1(k_R,l)\}$, with $l\in \mathbb R$,

where the functions $r_{1,2}$ and $p_1$ are defined by the following equations:
\begin{subequations}\label{maprp}
\begin{align}
&r_2(k_R,\lambda, l)=\frac{1}{2\pi} \int_{-\infty}^{\infty}d\xi\int_{0}^\infty d\eta\, e^{-il\xi+l(l+2k_R)\eta} q_0(\xi,\eta)\left(e_{\lambda}\breve \rho_2^+\right) (\xi,\eta,k_R,\lambda)\label{mapeq11},\\
&r_1(k_R,\lambda, l)=\frac{1}{2\pi}\int_{-\infty}^{\infty}d\xi\int_{0}^\infty d\eta\, e^{-il\xi+l(l+2k_R)\eta} q_0(\xi,\eta)\left(e_{\lambda}\breve \rho_1^+\right) (\xi,\eta,k_R,\lambda), \label{mapeq12}\\
&p_1(k_R,l)=\frac{1}{2\pi}\int_{-\infty}^{\infty}\!\!d\xi\int_0^T\!\!d\tau\, e^{-il\xi-4il(l^2+3k_Rl+3k_R^2)\tau}H(\xi,\tau,k_R,l)\phi^-_1(\xi,\tau,k_R,0)\label{mapeq13},
\end{align}
\end{subequations}
with $e_\lambda$ given by equation \eqref{elambda} and 
\begin{equation*}
\breve\rho_{1,2}^+(x,y,k_R,\lambda)=\rho^+_{1,2}(x,y,-k_R-\dfrac{\lambda}{2},\dfrac{i\lambda}{2}). 
\end{equation*}

\item $\{r_{1,2},\, p_1\}\mapsto \{\chi_{1,2}(k_R,\lambda), \psi_{1,2}(k_R,\lambda)\}$,\\
where the functions $\chi_{1,2}$ and $\psi_{1,2}$ are defined via the following integral equations:
{\small\begin{subequations}
\begin{align}\label{mapeq14}
&\chi_2(k_R,\lambda)-\int_{-\infty}^{\lambda}dl\chi_2(k_R,l) r_2(k_R, l, -2k_R-\lambda)\nonumber\\
&-\int_{-2k_R-\lambda}^\infty dl\chi_1(k_R, l)r_1(k_R,\lambda,-2k_R-\lambda)=p_1(k_R,-2k_R-\lambda), \quad k_R\leq0,\ \lambda\leq0,\\
&\chi_1(k_R,\lambda)-\int_{-\infty}^{-2k_R-\lambda}dl \chi_2(k_R,l)r_2(k_R, l, -2k_R-\lambda)\nonumber\\
&-\int_{\lambda}^{\infty}dl \chi_1(k_R, l)r_1(k_R,l, -2k_R-\lambda)=p_1(k_R,-2k_R-\lambda),\quad k_R\leq0,\ \lambda\geq-2k_R \label{mapeq15}
\end{align}
\end{subequations}
and
\begin{subequations}
\begin{align}\label{mapeq14p}
&\psi_2(k_R,\lambda)+\int_{-\infty}^{\lambda}dl\psi_2(k_R,l) r_2(k_R, l, -2k_R-\lambda)\nonumber\\
&+\!\int_0^{-2k_R-\lambda}\!\!\! dl\psi_1(k_R, l)r_1(k_R,\lambda,-2k_R-\lambda)=p_1(k_R,-2k_R-\lambda), \ k_R\geq0,\ \lambda\leq-2k_R,\\
&\psi_1(k_R,\lambda)+\int_{-\infty}^{-2k_R-\lambda}dl \psi_2(k_R,l)r_2(k_R, l, -2k_R-\lambda)\nonumber\\
&+\int_{0}^{\lambda}dl \psi_1(k_R, l)r_1(k_R,l, -2k_R-\lambda)=p_1(k_R,-2k_R-\lambda),\quad k_R\geq0,\ \lambda\geq0,\label{mapeq15p}
\end{align}
\end{subequations}}
\end{enumerate}
In summary, steps 1-5 provide the map \eqref{map} from $\{q_0,g,h\}$ to the spectral functions. 

As mentioned in the beginning of section \ref{invprob}, the function $\mu$ satisfied Pompeiu's formula \eqref{Pompeiu}. We notice that, under the definition \eqref{E} of the exponential $E$, equations \eqref{kbar1finalprop} and \eqref{kbar2finalprop} can be written as 
\begin{equation*}
\frac{\partial \mu_1^\pm}{\partial \bar k }(x,y,t,k_R,k_I)= E(x,y,t,k,-2k_R)\, \gamma_1^{\pm}(k_R,k_I)\,\mu_2^\pm(x,y,t,-k_R, k_I), \ k_R\leq 0
\end{equation*}
and
\begin{equation*}
\frac{\partial \mu_2^\pm}{\partial \bar k }(x,y,t,k_R,k_I)= -E(x,y,t,k,-2k_R)\, \gamma_2^{\pm}(k_R,k_I)\,\mu_1^\pm(x,y,t,-k_R, k_I),\ k_R \geq 0.
\end{equation*}
Using these expressions for the d-bar derivatives, as well as equations \eqref{discontinuity1prop} and \eqref{discontinuity2prop} for the discontinuities across the real $k$-axis, Pompeiu's formula \eqref{Pomp} yield the following expressions:
\begin{align}\label{Pompeiu2}
&\mu(x,y,t,k_R,k_I)=\nonumber\\
1&+\frac{1}{2i\pi} \int_{-\infty}^0 \frac{d\nu_R}{\nu_R-k}\Bigg[\int_{-\infty}^{0}d\lambda\,\chi_2\left(e_{\lambda}\breve \mu_2^+\right)+\int_{-2\nu_R}^{\infty}d\lambda\,\chi_1\left(e_{\lambda}\breve \mu_1^+\right)\Bigg](x,y,t,\nu_R,\lambda) \nonumber\\
&+\frac{1}{2i\pi}\int_{0}^{\infty} \frac{d\nu_R }{\nu_R-k}\,\Bigg[\int_{-\infty}^{-2\nu_R}d\lambda\,\psi_2\left(e_{\lambda}\breve \mu_2^+\right)+\int_{0}^{\infty}\!d\lambda\,\psi_1\left(e_{\lambda}\breve \mu_1^+\right)\!\Bigg](x,y,t,\nu_R,\lambda)\nonumber\\
&+\frac{1}{\pi}\int_{0}^{\infty}d\nu_R\int_{0}^{\infty}\frac{d\nu_I}{\nu-k}\,\gamma_2^{+}(\nu_R,\nu_I)E(x,y,t,\nu,-2\nu_R) \mu_1^+(x,y,t,-\nu_R, \nu_I)\nonumber\\
&-\frac{1}{\pi}\int_{-\infty}^{0}d\nu_R\int_{0}^{\infty}\frac{d\nu_I}{\nu-k}\,\gamma_1^{+}(\nu_R,\nu_I) E(x,y,t,\nu,-2\nu_R) \mu_2^+(x,y,t,-\nu_R, \nu_I)\nonumber\\
&+\frac{1}{\pi}\int_{0}^{\infty}d\nu_R\int_{-\infty}^{0} \frac{d\nu_I}{\nu-k}\,\gamma_2^{-}(\nu_R,\nu_I)E(x,y,t,\nu,-2\nu_R) \mu_1^-(x,y,t,-\nu_R, \nu_I)\nonumber\\
&-\frac{1}{\pi}\int_{-\infty}^{0}d\nu_R\int_{-\infty}^{0} \frac{d\nu_I}{\nu-k}\,\gamma_1^{-}(\nu_R,\nu_I)E(x,y,t,\nu,-2\nu_R) \mu_2^-(x,y,t,-\nu_R, \nu_I).
\end{align}
Using the transformation 
\begin{equation}\label{transhat}
\nu_R+\dfrac{\lambda}{2}=\hat \nu_R\quad \mathrm{and}\quad \dfrac{\lambda}{2}=\hat \nu_I, 
\end{equation}
the second and the third terms on the RHS of equation \eqref{Pompeiu2} become
\begin{align*}
&\frac{1}{i\pi} \int_{-\infty}^0d\hat \nu_R \int_{\hat \nu_R}^{0}d\hat \nu_I\,\frac{\chi_2(\hat \nu_R-\hat \nu_I,2\hat \nu_I)}{\hat \nu_R-\hat \nu_I-k}\,e_{\lambda}(x,y,t,\hat \nu_R-\hat \nu_I,2\hat \nu_I)\mu_2^+(x,y,t,-\hat \nu_R,i\hat \nu_I)\nonumber\\
+\,&\frac{1}{i\pi} \int_{0}^{\infty} d\hat \nu_R\int_{\hat \nu_R}^{\infty}d\hat\nu_I\,\frac{\chi_1(\hat \nu_R-\hat \nu_I,2\hat \nu_I)}{\hat \nu_R-\hat \nu_I-k}\,e_{\lambda}(x,y,t,\hat \nu_R-\hat \nu_I,2\hat \nu_I)\mu_1^+(x,y,t,-\hat \nu_R,i\hat \nu_I).
\end{align*}
\begin{figure}[ht]
\begin{center}
\resizebox{7cm}{!}{\input{jump1change1.pstex_t}}
\end{center}
\label{j1c1}
\end{figure}
\begin{figure}[ht]
\begin{center}
\resizebox{7cm}{!}{\input{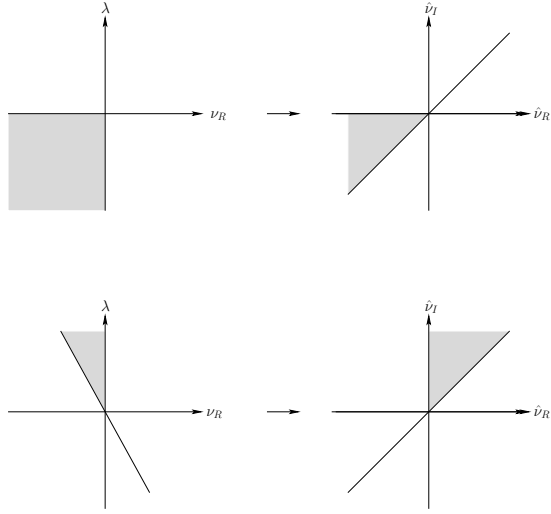}}
\end{center}
\caption{The regions of integration under the transformation \eqref{transhat}.}
\label{j2c2}
\end{figure}
Moreover, letting $\hat \nu_I=-i\tilde \nu_I$, we find
{\small\begin{align*}
&\frac{1}{\pi} \int_{-\infty}^0d\hat \nu_R \int_{0}^{i\hat \nu_R}d\tilde \nu_I\,\frac{\chi_2(\hat \nu_R+i\tilde \nu_I,-2i\tilde \nu_I)}{\hat \nu_R+i\tilde \nu_I-k}\,e_{\lambda}(x,y,t,\hat \nu_R+i\tilde \nu_I,-2i\tilde \nu_I)\mu_2^+(x,y,t,-\hat \nu_R,\tilde \nu_I)\nonumber\\
-\,&\frac{1}{\pi} \int_{0}^{\infty} d\hat \nu_R\int_{i\hat \nu_R}^{i\infty}d\tilde\nu_I\,\frac{\chi_1(\hat \nu_R+i\tilde \nu_I,-2i\tilde \nu_I)}{\hat \nu_R+i\tilde \nu_I-k}\,e_{\lambda}(x,y,t,\hat \nu_R+i \tilde \nu_I,-2i\tilde \nu_I)\mu_1^+(x,y,t,-\hat \nu_R,\tilde \nu_I).
\end{align*}}
The definitions \eqref{elambda} and \eqref{E} imply
\begin{equation*}
e_{\lambda}(x,y,t,\hat \nu_R+i\tilde \nu_I,-2i\tilde \nu_I)=E(x,y,t,\hat \nu_R+i\tilde \nu_I,-2\hat \nu_R).
\end{equation*}
Replacing the tilde by the hat and letting $\hat \nu=\nu_R+i \nu_I$, we can write the second and the third term as follows:
\begin{align}\label{jump1mod}
&\frac{1}{\pi} \int_{-\infty}^0d\hat \nu_R \int_{0}^{i\hat \nu_R}d\hat \nu_I\,\frac{1}{\hat \nu-k}\,\chi_2(\hat \nu,-2i\hat \nu_I)E(x,y,t,\hat \nu,-2\hat\nu_R)\mu_2^+(x,y,t,-\hat \nu_R,\hat \nu_I)\nonumber\\
-\,&\frac{1}{\pi} \int_{0}^{\infty} d\hat \nu_R\int_{i\hat \nu_R}^{i\infty}d\hat\nu_I\,\frac{1}{\hat \nu-k}\,\chi_1(\hat \nu,-2i\hat \nu_I)E(x,y,t,\hat \nu,-2\hat\nu_R)\mu_1^+(x,y,t,-\hat \nu_R,\hat \nu_I).
\end{align}
The same transformations applied on the fourth and the fifth term on the RHS of equation \eqref{Pompeiu2} yield
\begin{align}\label{jump2mod}
&\frac{1}{\pi} \int_{-\infty}^0d\hat \nu_R \int_{i\hat \nu_R}^{-i\infty}d\hat \nu_I\,\frac{1}{\hat \nu-k}\,\psi_2(\hat \nu,-2i\hat \nu_I)E(x,y,t,\hat \nu,-2\hat\nu_R)\mu_2^+(x,y,t,-\hat \nu_R,\hat \nu_I)\nonumber\\
-\,&\frac{1}{\pi} \int_{0}^{\infty} d\hat \nu_R\int_{0}^{i\hat \nu_R}d\hat\nu_I\,\frac{1}{\hat \nu-k}\,\psi_1(\hat \nu,-2i\hat \nu_I)E(x,y,t,\hat \nu,-2\hat\nu_R)\mu_1^+(x,y,t,-\hat \nu_R,\hat \nu_I).
\end{align}
Dropping the hats and inserting equations \eqref{jump1mod} and \eqref{jump2mod} into equation \eqref{Pompeiu2}, we conclude that
\begin{align}\label{Pompeiufinal}
&\mu(x,y,t,k_R,k_I)=\nonumber\\
1&+\frac{1}{\pi} \int_{-\infty}^0d \nu_R \int_{0}^{i\nu_R}\frac{d\nu_I}{ \nu-k}\,\chi_2(\nu,-2i \nu_I)E(x,y,t,\nu,-2\nu_R)\mu_2^+(x,y,t,- \nu_R, \nu_I)\nonumber\\
&-\frac{1}{\pi} \int_{0}^{\infty} d\nu_R\int_{i\nu_R}^{i\infty}\frac{d\nu_I}{ \nu-k}\,\chi_1(\nu,-2i\nu_I)E(x,y,t,\nu,-2\nu_R)\mu_1^+(x,y,t,-\nu_R,\nu_I)\nonumber\\
&+\frac{1}{\pi} \int_{-\infty}^0d \nu_R \int_{i \nu_R}^{-i\infty}\frac{d\nu_I}{ \nu-k}\,\psi_2( \nu,-2i \nu_I)E(x,y,t,\nu,-2\nu_R)\mu_2^+(x,y,t,- \nu_R, \nu_I)\nonumber\\
&-\frac{1}{\pi} \int_{0}^{\infty} d \nu_R\int_{0}^{i \nu_R}\frac{d\nu_I}{ \nu-k}\,\psi_1( \nu,-2i \nu_I)E(x,y,t,\nu,-2\nu_R)\mu_1^+(x,y,t,- \nu_R, \nu_I)\nonumber\\
&+\frac{1}{\pi}\int_{0}^{\infty}d\nu_R\int_{0}^{\infty}\frac{d\nu_I}{\nu-k}\,\gamma_2^{+}(\nu_R,\nu_I)E(x,y,t,\nu,-2\nu_R)\mu_1^+(x,y,t,-\nu_R, \nu_I)\nonumber\\
&-\frac{1}{\pi}\int_{-\infty}^{0}d\nu_R\int_{0}^{\infty}\frac{d\nu_I}{\nu-k}\,\gamma_1^{+}(\nu_R,\nu_I) E(x,y,t,\nu,-2\nu_R) \mu_2^+(x,y,t,-\nu_R, \nu_I)\nonumber\\
&+\frac{1}{\pi}\int_{0}^{\infty}d\nu_R\int_{-\infty}^{0} \frac{d\nu_I}{\nu-k}\,\gamma_2^{-}(\nu_R,\nu_I)E(x,y,t,\nu,-2\nu_R) \mu_1^-(x,y,t,-\nu_R, \nu_I)\nonumber\\
&-\frac{1}{\pi}\int_{-\infty}^{0}d\nu_R\int_{-\infty}^{0} \frac{d\nu_I}{\nu-k}\,\gamma_1^{-}(\nu_R,\nu_I)E(x,y,t,\nu,-2\nu_R) \mu_2^-(x,y,t,-\nu_R, \nu_I).
\end{align}

\begin{prop}\label{qfinal}
Assume that $q(x,y,t)$ satisfies an IBV problem for the KPII equation \eqref{kpiiintro} with given initial condition $q(x,y,0)=q_0(x,y)$. Let
\begin{equation*}
q(x,0,t)=g(x,t),\quad q_y(x,0,t)=h(x,t).
\end{equation*}
Assume that the following global relation is valid:
{\small\begin{align*}
&\int_{-\infty}^{\infty}d\xi\int_{0}^\infty d\eta\, e^{-il\xi+l(l+2k)\eta} q_0(\xi,\eta)\rho^-_{j}-\int_{-\infty}^{\infty}\!\!d\xi\int_{0}^{t}d\tau\, e^{-il\xi-4il(l^2+3kl+3k^2)\tau}H(\xi,\tau,k,l)\phi^-_{j}\nonumber\\
&=\int_{-\infty}^{\infty}d\xi\int_{0}^\infty\!\!\!\! d\eta\, e^{-il\xi+l(l+2k)\eta-4il(l^2+3kl+3k^2)t} q(\xi,\eta,t)\mu^-_{j},\ j=1,2, \,  k \in \mathbb C,\, l(l+2k_R)\leq0,
\end{align*}}
where $\mu^-_j$ are bounded $\forall k \in \mathbb C$, $\rho_j^-$, $\phi^-_j$ are defined by equations \eqref{maprhob}, \eqref{mapphi} and $H$ is given by \eqref{H}, i.e. by
\begin{equation*}
H(x,t,k,l)=3\bigg[g_x(x,t)-2i(l+k)g(x,t)-\partial_{x}^{-1}h(x,t)\bigg].
\end{equation*}

Given $q_0(x,y)$, $g(x,t)$ and $h(x,t)$, define the functions $\gamma^\pm_{1,2}(k_R,k_I)$ by
\begin{align*}
\gamma_1^+
(k_R,k_I)&=\beta_1^+(k_R,k_I), &k_R\leq0,\, k_I\geq0,\\
\gamma_1^-
(k_R,k_I)&=\beta_1^-(k_R,k_I)-\alpha^-_1(k_R,k_I), &k_R\leq0,\, k_I\leq0,\\
\gamma_2^+
(k_R,k_I)&=\beta_2^+(k_R,k_I), &k_R\geq0,\, k_I\geq0,\\
\gamma_2^-
(k_R,k_I)&=\beta_2^-(k_R,k_I)-\alpha^-_2(k_R,k_I), &k_R\geq0,\, k_I\leq0
\end{align*}
where the functions $\alpha^-_{1,2}(k_R,k_I),\ \beta_{1,2}^\pm(k_R,k_I)$ are given by
\begin{align*}
\alpha_1^-(k_R,k_I)&=\frac{1}{2\pi}\int_{-\infty}^{\infty}\!d\xi\int_0^T\!\!\!d\tau\, e^{2ik_R\xi+8ik_R(k_R^2-3k_I^2)\tau} H(\xi,\tau,k,-2k_R)\phi_1^-(\xi,\tau,k_R,k_I),\\
\beta_1^{\pm}(k_R,k_I)&=\frac{1}{2\pi} \int_{-\infty}^{\infty}d\xi\int_0^\infty d\eta\, e^{2ik_R\xi-4ik_Rk_I\eta}q_0(\xi,\eta)\rho_{1}^\pm(\xi, \eta,k_R,k_I),\\
\alpha_2^-(k_R,k_I)&=\frac{1}{2\pi}\int_{-\infty}^{\infty}\!d\xi\int_0^T\!\!\!d\tau\, e^{2ik_R\xi+8ik_R(k_R^2-3k_I^2)\tau}H(\xi,\tau,k,-2k_R)\phi_2^-(\xi,\tau,k_R,k_I),\\
\beta_2^{\pm}(k_R,k_I)&=\frac{1}{2\pi} \int_{-\infty}^{\infty}d\xi\int_0^\infty d\eta\ e^{2ik_R\xi-4ik_Rk_I\eta}q_0(\xi,\eta)\rho_{2}^\pm(\xi, \eta,k_R,k_I).
\end{align*}
Assume that the following systems of Volterra integral equations have a unique solution for the functions $\chi_{1,2}(k_R,\lambda)$ and $\psi_{1,2}(k_R,\lambda)$ respectively:
\begin{align*}
&\chi_2(k_R,\lambda)-\int_{-\infty}^{\lambda}dl\chi_2(k_R,l) r_2(k_R, l, -2k_R-\lambda)\nonumber\\
&-\int_{-2k_R-\lambda}^\infty dl\chi_1(k_R, l)r_1(k_R,\lambda,-2k_R-\lambda)=p_1(k_R,-2k_R-\lambda), \quad k_R\leq0,\ \lambda\leq0,\\
&\chi_1(k_R,\lambda)-\int_{-\infty}^{-2k_R-\lambda}dl \chi_2(k_R,l)r_2(k_R, l, -2k_R-\lambda)\nonumber\\
&-\int_{\lambda}^{\infty}dl \chi_1(k_R, l)r_1(k_R,l, -2k_R-\lambda)=p_1(k_R,-2k_R-\lambda),\quad k_R\leq0,\ \lambda\geq-2k_R 
\end{align*}
and
\begin{align*}
&\psi_2(k_R,\lambda)+\int_{-\infty}^{\lambda}dl\psi_2(k_R,l) r_2(k_R, l, -2k_R-\lambda)\nonumber\\
&+\!\int_0^{-2k_R-\lambda}\!\!\! dl\psi_1(k_R, l)r_1(k_R,\lambda,-2k_R-\lambda)\!=p_1(k_R,-2k_R-\lambda), \ k_R\geq0,\, \lambda\leq-2k_R,\\
&\psi_1(k_R,\lambda)+\int_{-\infty}^{-2k_R-\lambda}dl \psi_2(k_R,l)r_2(k_R, l, -2k_R-\lambda)\nonumber\\
&+\int_{0}^{\lambda}dl \psi_1(k_R, l)r_1(k_R,l, -2k_R-\lambda)=p_1(k_R,-2k_R-\lambda),\quad k_R\geq0,\ \lambda\geq0,
\end{align*}
where
\begin{align*}
&r_1(k_R,\lambda, l)=\frac{1}{2\pi}\int_{-\infty}^{\infty}\!\!\!\!\!\!d\xi\int_{0}^\infty\!\!\!\!\!\!d\eta\, e^{-il\xi+l(l+2k_R)\eta} q_0(\xi,\eta)\left(e_{\lambda}\breve \rho_1^+\right) (\xi,\eta,k_R,\lambda),\, \lambda \geq-2k_R,\\
&r_2(k_R,\lambda, l)=\frac{1}{2\pi} \int_{-\infty}^{\infty}d\xi\int_{0}^\infty\!\! d\eta\, e^{-il\xi+l(l+2k_R)\eta} q_0(\xi,\eta)\left(e_{\lambda}\breve \rho_2^+\right) (\xi,\eta,k_R,\lambda),\ \lambda \geq 0,\\
&p_1(k_R,l)=\frac{1}{2\pi}\int_{-\infty}^{\infty}d\xi\int_0^Td\tau\, e^{-il\xi-4il(l^2+3k_Rl+3k_R^2)\tau}H(\xi,\tau,k_R,l)\phi^-_1(\xi,\tau,k_R,0),
\end{align*}
and $e_\lambda$ is given by equation \eqref{Ee}.

Given $\{\chi_{1,2},\, \psi_{1,2},\, \gamma_{1,2}^\pm\}$, define $\mu$ as the solution of equation \eqref{Pompeiufinal}, where 
\begin{align*}
\mu(x,y,t,k_R,k_I) = \left\lbrace \begin{array}{ll} \mu_{1}^{+}(x,y,t,k_R,k_I), &\quad k_R\leq0,\ k_I\geq0, \\ \\ \mu_{1}^{-}(x,y,t,k_R,k_I), &\quad k_R\leq0,\ k_I\leq0,  \\ \\ \mu_{2}^{-}(x,y,t,k_R,k_I), &\quad k_R\geq0,\ k_I\leq0, \\ \\ \mu_{2}^{+}(x,y,t,k_R,k_I), &\quad k_R\geq0,\ k_I\geq0, \end{array} \right.
\end{align*}
and assume that equation \eqref{Pompeiufinal} has a unique solution.

Define $q(x,y,t)$ by
\begin{equation*}
q(x,y,t)=-2i\lim_{k\rightarrow \infty}\bigg[k\mu_x(x,y,t,k_R,k_I)\bigg].
\end{equation*}
Then $q(x,y,t)$ solves the KPII equation \eqref{kpiiintro}.
\end{prop}

\textbf{Remark 3.1} By differentiating equation \eqref{Pompeiufinal} with respect to $x$ we can write an integral representation for the solution $q$:
\begin{align*}
\frac{i\pi}{2}\,q&=-\int_{-\infty}^0\!\!d\nu_R \int_{0}^{i\nu_R}\!d\nu_I\,\chi_2(\nu,-2i \nu_I)\Bigg(E(x,y,t,\nu,-2\nu_R)\mu_2^+(x,y,t,- \nu_R, \nu_I)\Bigg)_x\nonumber\\
&+ \int_{0}^{\infty} d\nu_R\int_{i\nu_R}^{i\infty}d\nu_I\,\chi_1(\nu,-2i\nu_I)\Bigg(E(x,y,t,\nu,-2\nu_R)\mu_1^+(x,y,t,-\nu_R,\nu_I)\Bigg)_x\nonumber\\
&- \int_{-\infty}^0d \nu_R \int_{i \nu_R}^{-i\infty}d \nu_I\,\psi_2( \nu,-2i \nu_I)\Bigg(E(x,y,t,\nu,-2\nu_R)\mu_2^+(x,y,t,- \nu_R, \nu_I)\Bigg)_x\nonumber\\
&+ \int_{0}^{\infty} d \nu_R\int_{0}^{i \nu_R}d\nu_I\,\psi_1( \nu,-2i \nu_I)\Bigg(E(x,y,t,\nu,-2\nu_R)\mu_1^+(x,y,t,- \nu_R, \nu_I)\Bigg)_x\nonumber\\
&-\int_{0}^{\infty}d\nu_R\int_{0}^{\infty}d\nu_I\,\gamma_2^{+}(\nu_R,\nu_I)\Bigg(E(x,y,t,\nu,-2\nu_R) \mu_1^+(x,y,t,-\nu_R, \nu_I)\Bigg)_x\nonumber\\
&+\int_{-\infty}^{0}d\nu_R\int_{0}^{\infty}d\nu_I\,\gamma_1^{+}(\nu_R,\nu_I) \Bigg(E(x,y,t,\nu,-2\nu_R) \mu_2^+(x,y,t,-\nu_R, \nu_I)\Bigg)_x\nonumber\\
&-\int_{0}^{\infty}d\nu_R\int_{-\infty}^{0} d\nu_I\,\gamma_2^{-}(\nu_R,\nu_I)\Bigg(E(x,y,t,\nu,-2\nu_R) \mu_1^-(x,y,t,-\nu_R, \nu_I)\Bigg)_x\nonumber\\
&+\int_{-\infty}^{0}d\nu_R\int_{-\infty}^{0} d\nu_I\,\gamma_1^{-}(\nu_R,\nu_I)\Bigg(E(x,y,t,\nu,-2\nu_R) \mu_2^-(x,y,t,-\nu_R, \nu_I)\Bigg)_x.
\end{align*}
\vskip 3mm
\textbf{Remark 3.2} (The case of zero initial condition)
\vskip 2mm
Consider the special case $q_0(x,y)=0$. Then, equations \eqref{chi2eqfinal}-\eqref{psi1eqfinal} imply that the functions $\chi_{1,2}$ and $\psi_{1,2}$ are given explicitly by
\begin{equation*}
 \chi_2(k_R,\lambda)=\chi_1(k_R,\lambda)=\psi_2(k_R,\lambda)=\psi_1(k_R,\lambda)=p_1(k_R,-2k_R-\lambda),
\end{equation*}
with $p_1$ defined by equation \eqref{p1}. 

Also, equations \eqref{beta1prop} and \eqref{beta2prop} yield
\begin{equation*}
\beta_1^\pm(k_R,k_I)=\beta_2^\pm(k_R,k_I)=0,
\end{equation*}
thus, according to equations \eqref{gamma1pprop}-\eqref{gamma2mprop}:
\begin{equation*}
\gamma_1^+(k_R,k_I)=\gamma_2^+(k_R,k_I)=0
\end{equation*}
and 
\begin{equation*}
\gamma_1^-(k_R,k_I)=-\alpha_1^-(k_R,k_I), \quad \gamma_2^-(k_R,k_I)=-\alpha_2^-(k_R,k_I),
\end{equation*}
where $\alpha_{1,2}^-$ are defined by equations \eqref{alpha1-prop} and \eqref{alpha2-prop}.

Hence, equation \eqref{Pompeiufinal} for $\mu$ becomes significantly simpler, i.e. 
\begin{align*}
&\mu(x,y,t,k_R,k_I)=\nonumber\\
1&+\frac{1}{\pi} \int_{-\infty}^0d \nu_R \int_{0}^{i\nu_R}\frac{d\nu_I}{ \nu-k}\,p_1(\nu,-2\nu_R)E(x,y,t,\nu,-2\nu_R)\mu_2^+(x,y,t,- \nu_R, \nu_I)\nonumber\\
&-\frac{1}{\pi} \int_{0}^{\infty} d\nu_R\int_{i\nu_R}^{i\infty}\frac{d\nu_I}{ \nu-k}\,p_1(\nu,-2\nu_R)E(x,y,t,\nu,-2\nu_R)\mu_1^+(x,y,t,-\nu_R,\nu_I)\nonumber\\
&+\frac{1}{\pi} \int_{-\infty}^0d \nu_R \int_{i \nu_R}^{-i\infty}\frac{d\nu_I}{ \nu-k}\,p_1(\nu,-2\nu_R)E(x,y,t,\nu,-2\nu_R)\mu_2^+(x,y,t,- \nu_R, \nu_I)\nonumber\\
&-\frac{1}{\pi} \int_{0}^{\infty} d \nu_R\int_{0}^{i \nu_R}\frac{d\nu_I}{ \nu-k}\,p_1(\nu,-2\nu_R)E(x,y,t,\nu,-2\nu_R)\mu_1^+(x,y,t,- \nu_R, \nu_I)\nonumber\\
&-\frac{1}{\pi}\int_{0}^{\infty}d\nu_R\int_{-\infty}^{0} \frac{d\nu_I}{\nu-k}\,\alpha_2^{-}(\nu_R,\nu_I)E(x,y,t,\nu,-2\nu_R) \mu_1^-(x,y,t,-\nu_R, \nu_I)\nonumber\\
&+\frac{1}{\pi}\int_{-\infty}^{0}d\nu_R\int_{-\infty}^{0} \frac{d\nu_I}{\nu-k}\,\alpha_1^{-}(\nu_R,\nu_I)E(x,y,t,\nu,-2\nu_R) \mu_2^-(x,y,t,-\nu_R, \nu_I),
\end{align*}
where $E(x,y,t,k,l)=e^{ilx-l(l+2k)y+4il(l^2+3kl+3k^2)t}$, the functions $p_1$ and $\alpha_{1,2}^-$ are given by
\begin{align*}
p_1(k_R,l)&=\frac{1}{2\pi}\int_{-\infty}^{\infty}d\xi\int_0^Td\tau\, e^{-il\xi-4il(l^2+3k_Rl+3k_R^2)\tau}H(\xi,\tau,k_R,l)\phi^-_1(\xi,\tau,k_R,0),\\
\alpha_1^-(k_R,k_I)&=\frac{1}{2\pi}\int_{-\infty}^{\infty}\!d\xi\int_0^T\!\!\!d\tau\, e^{2ik_R\xi+8ik_R(k_R^2-3k_I^2)\tau} H(\xi,\tau,k,-2k_R)\phi_1^-(\xi,\tau,k_R,k_I),\\
\alpha_2^-(k_R,k_I)&=\frac{1}{2\pi}\int_{-\infty}^{\infty}\!d\xi\int_0^T\!\!\!d\tau\, e^{2ik_R\xi+8ik_R(k_R^2-3k_I^2)\tau}H(\xi,\tau,k,-2k_R)\phi_2^-(\xi,\tau,k_R,k_I)
\end{align*}
and $\phi_{1,2}^-$ are defined by equations \eqref{mapphi}.
\subsection{The linear limit}\label{kpiill}

Suppose now that $q$ is small, i.e. let 
\begin{equation*}
q(x,y,t)=\varepsilon u(x,y,t)+O(\varepsilon^2), \quad \varepsilon\rightarrow 0.
\end{equation*}
Then, the $O(\varepsilon)$ term of the nonlinear KPII equation \eqref{kpiiintro} yields a linear equation which, as expected, is the linearised KP equation:
\begin{align*}
u_t+u_{xxx}+3\partial^{-1}_x u_{yy}=0.
\end{align*}
The following integral representation of the solution of an IBV problem was derived in section \ref{lkp}:
{\small\begin{align}\label{u}
&u=\frac{2}{\pi^2}\mbox{\large $\Big($}\int_{0}^{\infty}d\nu_R-\int_{-\infty}^{0}d\nu_R\mbox{\large $\Big)$}\int_{-\infty}^{\infty}d\nu_I\int_{-\infty}^{\infty}d\xi \int_{0}^{\infty} d\eta\,  E(x-\xi,y-\eta,t,\nu,-2\nu_R)\nu_R\, u_0\nonumber\\
&-\frac{2}{\pi^2}\int_{0}^\infty\!\!\!d\nu_R\mbox{\large $\Big($}\int_{-\infty}^0\!\!\!\!d\nu_I+\!\int_{0}^{i \infty}\!\!\!\!d\nu_I\mbox{\large $\Big)$}\!\!\int_{-\infty}^{\infty}\!\!\!\!\!d\xi\int_{0}^{T}\!\!\!\!d\tau\, E(x-\xi,y,t-\tau,\nu,-2\nu_R)\,3\nu_R\mbox{\large $\big($}2\nu_I v\!-\partial_\xi^{-1}w\mbox{\large $\big)$}\nonumber\\
&+\frac{2}{\pi^2}\int_{0}^\infty\!\!\!d\nu_R\mbox{\large $\Big($}\!\int_{-\infty}^0\!\!\!\!d\nu_I+\!\int_{0}^{-i \infty}\!\!\!\!\!\!\!d\nu_I\!\mbox{\large $\Big)$}\!\!\int_{-\infty}^{\infty}\!\!\!\!\!d\xi\int_{0}^{T}\!\!\!\!d\tau\, E(x-\xi,y,t-\tau,\nu,-2\nu_R)\,3\nu_R\mbox{\large $\big($}2\nu_I v\!-\partial_\xi^{-1}w\mbox{\large $\big)$}\!,
\end{align}}
\vskip-5mm
\noindent where \vskip-5mm
\begin{equation*}
u_0(x,y)=u(x,y,0),\quad v(x,t)=u(x,0,t),\quad w(x,t)=u_y(x,0,t).
\end{equation*}
For the nonlinear problem, equations \eqref{direct1ii} and \eqref{direct2ii} suggest that, in the linear limit, $\mu=1+O(\varepsilon)$. 

Moreover, equations \eqref{r2}, \eqref{r1} and \eqref{p1} imply that
\begin{equation*}
r_2=O(\varepsilon),\quad r_1=O(\varepsilon), \quad p_1=O(\varepsilon),
\end{equation*}
and then, equations \eqref{chi2eqfinal} and \eqref{chi1eqfinal} suggest that $\chi_2=O(\varepsilon)$ and $\chi_1=O(\varepsilon)$. 

In fact, we can write
\begin{equation*}
\mu=1+\varepsilon M+O(\varepsilon^2), \ p_1=\varepsilon P_1+O(\varepsilon^2), \ \chi_2=\varepsilon X_2+O(\varepsilon^2), \ \chi_1=\varepsilon X_1+O(\varepsilon^2),
\end{equation*}
for some functions $M$, $P_1$, $X_2$ and $X_1$, independent of $\varepsilon$. Then, equation \eqref{p1} yields
{\small\begin{equation*}
P_1(k_R,l)=\frac{1}{2\pi}\int_{-\infty}^{\infty}d\xi\int_0^Td\tau\, E(-\xi,0,-\tau,k_R,l)\, 3\mbox{\large $\big[$}v_\xi(\xi,\tau)-2i(l+k_R)v(\xi,\tau)-\partial_\xi^{-1}w(\xi,\tau)\mbox{\large $\big]$}.
\end{equation*}}
Integrating by parts, we find
\begin{equation*}
P_1(k_R,l)=\frac{1}{2\pi}\int_{-\infty}^{\infty}d\xi\int_0^Td\tau\, E(-\xi,0,-\tau,k_R,l)\, 3\bigg[-i(l+2k_R)v(\xi,\tau)-\partial_\xi^{-1}w(\xi,\tau)\bigg],
\end{equation*}
therefore, equations \eqref{chi2eqfinal} and \eqref{chi1eqfinal} yield
\begin{equation}\label{X12}
X_2(k_R,\lambda)=X_1(k_R,\lambda)=P_1(k_R,-2k_R-\lambda).
\end{equation}
A similar analysis for $\psi_1$ and $\psi_2$ implies
\begin{equation*}
\Psi_2(k_R,\lambda)=\Psi_1(k_R,\lambda)=P_1(k_R,-2k_R-\lambda).
\end{equation*}
Furthermore, equations \eqref{alpha1-prop}-\eqref{beta2prop} suggest that $A_{1,2}^\pm=A$ and $B_{1,2}^\pm=B$, where
\begin{align}
B(k_R,k_I)&=\frac{1}{2\pi}\int_{-\infty}^{\infty}d\xi\int_{0}^{\infty}d\eta\, E(-\xi,-\eta,0,k,-2k_R) u_0(\xi,\eta), \label{Beta}\\
A(k_R,k_I)&=\frac{1}{2\pi}\int_{-\infty}^{\infty}\!\!\!\!d\xi\int_0^T\!\!\!\!d\tau\,E(-\xi,0,-\tau,k,-2k_R)\, 3\bigg[2k_Iv(\xi,\tau)-\partial_\xi^{-1}w(\xi,\tau)\bigg].\label{Alpha}
\end{align}
Hence it follows from the definitions \eqref{gamma1pprop} and \eqref{gamma2pprop} that $\gamma_1^+=\gamma_2^+$ and $\gamma_1^-=\gamma_2^-$.

Then, the $O(\varepsilon)$ term of equation \eqref{Pompeiufinal} yields
\begin{align}\label{Mu}
M&=\frac{1}{\pi} \left(\int_{-\infty}^0\!\!\!d\nu_R \int_{0}^{-i\infty}\!\!\!\!d\nu_I-\int_{0}^{\infty}\!\!\!\!d\nu_R\int_{0}^{i\infty}\!\!\!d\nu_I\right)\frac{1}{ \nu-k}\,E(x,y,t,\nu,-2\nu_R)P_1(\nu,-2\nu_R)\nonumber\\
&+\frac{1}{\pi}\left(\int_{0}^{\infty}\!d\nu_R\int_{-\infty}^{\infty}d\nu_I-\int_{-\infty}^0\!d\nu_R \int_{-\infty}^{\infty}d\nu_I\right)\frac{1}{\nu-k}\,E(x,y,t,\nu,-2\nu_R)B(\nu_R,\nu_I)\nonumber\\
&-\frac{1}{\pi}\left(\int_{0}^{\infty}\!d\nu_R\int_{-\infty}^0d\nu_I-\int_{-\infty}^0\!d\nu_R \int_{-\infty}^0d\nu_I\right)\frac{1}{\nu-k}\,E(x,y,t,\nu,-2\nu_R)A(\nu_R,\nu_I).
\end{align}
Employing the formulae \eqref{X12}, \eqref{Beta} and \eqref{Alpha} into equation \eqref{Mu}, and using the definition
\begin{equation*}
u(x,y,t)=-2i\lim_{k\rightarrow \infty}\bigg[k M_x(x,y,t,k_R,k_I)\bigg],
\end{equation*}
we recover the expression \eqref{u}.

\section{Conclusion}
It is shown in \cite{F2002b} that there exists a much simpler derivation for the integral representation of $q$ for the linearised KPII equation than the one presented in section \ref{lkp}. However, this simpler approach cannot be ``nonlinearised''. On the other hand, the approach presented in section \ref{lkp} contains all necessary conceptual steps needed for the solution of KPII. However, the analysis of KPII is technically quite involved. In particular, in order to formulate proposition \eqref{qfinal} , it is necessary to compute both the d-bar derivatives $\partial \mu_j^\pm / \partial \bar k$ and the jumps $(\mu_j^+-\mu_j^-),\ j=1,2$; the latter computation is quite complicated, as was also the case for the analogous computation is DSII \cite{FDS2009}.

We note that an additional novelty of the IBV problems in 2+1, as opposed to IBV problems in 1+1, is that now the global relation plays an even more important role. Namely, it is necessary in order to express the jumps of $\mu$ across the real $k$-axis in terms of the initial and the boundary values.

Several problems remain open, in particular:
\begin{enumerate}
\item For a well posed problem, either $g$ or $h$ are prescribed as boundary conditions; on the other hand, the solution obtained in proposition \ref{qfinal} depends on \textit{both} $g$ and $h$. Thus, in order for this solution to be effective it is necessary to use the global relation to eliminate the unknown boundary value. For the linearised KPII, this is achieved in \cite{F2002b}; the analogous problem for the KPII remains open. In this connection we note that until recently, for problems in 1+1, it was not possible to express the spectral functions directly in terms of the given boundary conditions, but it was necessary to first obtain the unknown boundary values , i.e. to determine the so-called Dirichlet to Neumann map (see \cite{MonvelFShep2003}-\cite{FTreharne2008}). However, in a recent breakthrough, it has been shown that for the nonlinear Schr\"odinger (NLS) equation on the half-line it is possible to express the spectral functions directly in terms of the given initial and boundary conditions. The question of whether this new approach can be extended to equations in 2+1 remains open.

\item For the initial-value problem of the KPII, the linear integral equation analogous to \eqref{Pompeiufinal} admits a unique solution for $\mu$ for real initial conditions. This is due to the existence of a so-called ``vanishing lemma'', which is based on the theory of generalised analytic functions of Vekua. The question of whether there exists an analogous result for equation \eqref{Pompeiufinal} remains open.

\item The spectral functions are defined in terms of linear integral equations. The question of existence and uniqueness for these equations remains open. Actually, it is expected that these equations do possess homogeneous solutions which will give rise to coherent structures.

\item In spite of the fact that the representation of $q$ involves both $g$ and $h$, it should still be possible to obtain effective formulae for the large $t$ asymptotics of the solution. (for equations in 1+1 this has been achieved by the Deift-Zhou approach \cite{DeiftZhou1992}, \cite{DeiftZhou1993})

\item The new method is particularly effective for a class of boundary conditions called linearisable. In this case, it is possible to express the spectral functions directly in terms of the given initial and boundary conditions using only \textit{algebraic} manipulations, see for example \cite{F2002a}-\cite{FLenells2010b}. The question of identifying lineasable boundary conditions for the KPII remains open. 

\item For KPII the formalism presented here involves the crucial assumption that several linear integral equations have a unique solution. In spite of the fact that these equations are of Fredholm type, it is not difficult to establish uniqueness under the assumption of sufficiently ``small data''. However, the elimination of the ``small norm'' assumption is a formidable task (see also remark 3 above). 
\end{enumerate}

\section*{Acknowledgments}
ASF is grateful to Guggenheim Foundation for partial support and to V. Zakharov for suggesting the investigation of KPII on the half-plane.\vskip 3mm
\noindent D. Mantzavinos was supported by an EPSRC Doctoral Training Grant.
\end{document}